\documentclass[a4paper, 11pt]{article}
\usepackage[a4paper, margin=1in]{geometry}
\usepackage[T1]{fontenc} 
\usepackage{lmodern} 
\usepackage{microtype} 
\usepackage{amsmath} 
\usepackage{amssymb} 
\usepackage{amsthm} 
\usepackage{amsfonts} 
\usepackage{booktabs} 
\usepackage{multirow} 
\usepackage{mathrsfs} 
\usepackage[title,titletoc]{appendix} 
\usepackage{graphicx} 
\usepackage{subcaption} 
\usepackage{tikz}
\usetikzlibrary{positioning,arrows.meta,calc,decorations.pathreplacing}

\usepackage[
	backend=biber,
	style=trad-abbrv,
	citestyle=numeric-comp,
	sorting=none, 
	maxbibnames=99, 
	isbn=false,
	url=false,
	doi=true,
	eprint=false,
	date=year
]{biblatex}
\allowdisplaybreaks

\numberwithin{equation}{section}

\usepackage{hyperref} 
\hypersetup{
	colorlinks=true,
	linkcolor=blue,
	citecolor=blue,
	urlcolor=blue
}
\usepackage{cleveref} 

\theoremstyle{plain}
\newtheorem{theorem}{Theorem}[section]
\newtheorem{lemma}[theorem]{Lemma}
\newtheorem{proposition}[theorem]{Proposition}
\newtheorem{corollary}[theorem]{Corollary}

\theoremstyle{definition}
\newtheorem{definition}[theorem]{Definition}
\newtheorem{assumption}[theorem]{Assumption}
\newtheorem{algorithm}{Algorithm}[section]
\newtheorem{example}{Example}[section]

\theoremstyle{remark}
\newtheorem{remark}{Remark}[section]

\newcommand{\fixthmcref}[1]{%
	\AddToHook{env/#1/begin}{\crefalias{theorem}{#1}}%
}
\fixthmcref{lemma}
\fixthmcref{proposition}
\fixthmcref{corollary}
\fixthmcref{definition}
\fixthmcref{assumption}

\newcommand{\setcrefnames}[5]{\crefname{#1}{#2}{#3}\Crefname{#1}{#4}{#5}}
\setcrefnames{assumption}{assumption}{assumptions}{Assumption}{Assumptions}
\setcrefnames{algstep}{step}{steps}{Step}{Steps}
\setcrefnames{claim}{part}{parts}{Part}{Parts}
\setcrefnames{condition}{condition}{conditions}{Condition}{Conditions}

\newcommand{\dd}{\,\mathrm d}
\newcommand{\R}{\mathbb R}
\newcommand{\Pk}{\mathbb P}
\newcommand{\Qk}{\mathbb Q}
\newcommand{\GE}{\mathcal G_{\mathrm E}}
\newcommand{\vect}[1]{\mathbf{#1}}
\newcommand{\mat}[1]{\mathsf{#1}}
\newcommand{\Ark}{\mat A}
\newcommand{\brk}{\vect b}

\newcommand{\stateT}{\mat T}
\newcommand{\bx}{\vect x}

\newcommand{\bn}{\vect n}
\newcommand{\bfv}{\vect f}
\newcommand{\DG}{\mathcal F}
\newcommand{\loc}{\mathcal G}
\newcommand{\Th}{\mathcal T_h}
\newcommand{\Vh}{\mathbb V_h}
\newcommand{\ZK}{\mathbb Z^{k_K}(K)}
\newcommand{\FF}{\boldsymbol{\mathscr F}}
\newcommand{\calM}{\mathcal M}
\newcommand{\norm}[1]{\lVert #1\rVert}
\newcommand{\abs}[1]{\lvert #1\rvert}
\newcommand{\avg}[1]{\langle #1\rangle}
\newcommand{\Ah}{\mathcal A_h}

\DeclareMathOperator{\diam}{diam}

\hypersetup{
	pdftitle={PACL: Patch-based Admissibility and Conservative Limiting framework for bound-preserving compact single-step DG methods for hyperbolic conservation laws}, 
	pdfauthor={Yue Wu; Zheng Sun}
}
\usepackage{mathtools}
\mathtoolsset{showonlyrefs} 

\title{PACL: Patch-based Admissibility and Conservative Limiting framework for bound-preserving compact single-step DG methods for hyperbolic conservation laws}

\author{
	Yue~Wu\thanks{Division of Applied Mathematics, Brown University, Providence, RI 02912, USA. Email: \href{mailto:yue_wu3@brown.edu}{\texttt{yue\_wu3@brown.edu}}. ORCID: \href{https://orcid.org/0009-0006-3898-8204}{\texttt{0009-0006-3898-8204}}.}
	\and Zheng~Sun\thanks{Corresponding author. Department of Mathematics, The University of Alabama, Tuscaloosa, AL 35487, USA. Email: \href{mailto:zsun30@ua.edu}{\texttt{zsun30@ua.edu}}. ORCID: \href{https://orcid.org/0000-0003-3763-3015}{\texttt{0000-0003-3763-3015}}.}
}

\date{}
\begin{document}

\maketitle
\vspace{-1.25em}
\begin{abstract} We develop a Patch-based Admissibility and Conservative Limiting (PACL) framework for compact single-step discontinuous Galerkin (DG) discretizations of hyperbolic conservation laws. The key idea is to establish admissibility of weighted averages over overlapping cell patches after a complete high-order time step. Continuous patch weights cancel interior numerical flux contributions, yielding admissibility estimates for weighted patch averages under suitable time-step conditions. A conservative patch correction uses these averages as reference states to obtain admissible physical cell averages, and a subsequent cellwise scaling enforces whole-cell or nodal admissibility in the original DG space. The framework requires neither a strong-stability-preserving decomposition nor auxiliary first-order discretizations, and all limiting is performed at the end of each time step. For scalar conservation laws, we prove a patchwise weak maximum principle for compact Runge--Kutta DG (cRKDG), ADER-DG with converged or fixed-Picard predictors, and Lax--Wendroff DG. We also prove that, under the stated local mesh and scaling assumptions, a single application of the scalar limiter preserves approximation order at active bounds, including flat extrema. For the compressible Euler equations, we establish patch admissibility for cRKDG and fixed-Picard ADER-DG under time-step conditions that maintain positive predictor density while allowing nonpositive predictor pressure. The analysis covers both exact spatial integration and compatible quadrature. Numerical experiments demonstrate high-order accuracy and bound preservation, including problems with near-vacuum states and strong shocks.
\end{abstract}

\begin{paragraph}{Keywords:}
	Bound-preserving methods; discontinuous Galerkin methods; compact Runge--Kutta DG; ADER-DG; hyperbolic conservation laws; compressible Euler equations
\end{paragraph}

\begin{paragraph}{MSC 2020 codes:}
	Primary 65M60; Secondary 35L65, 65M12, 76M10.
\end{paragraph}

\section{Introduction}\label{sec:introduction}

Discontinuous Galerkin (DG) methods combine high-order accuracy, geometric flexibility, and compact spatial discretizations for hyperbolic conservation laws \cite{10.1023/A:1012873910884}. Several fully discrete formulations retain this locality in the temporal discretization. The compact Runge--Kutta DG (cRKDG) method uses a cell-local operator at the inner Runge--Kutta (RK) stages and restores conservative inter-element coupling in the completed update \cite{10.1137/23M158629X}. A class of ADER-DG methods constructs a local space-time predictor within each cell, followed by a conservative corrector \cite{10.1016/j.jcp.2010.10.024}. Lax--Wendroff DG (LWDG) methods instead approximate the time-averaged flux through a Taylor expansion and Cauchy--Kowalevski (CK) substitutions \cite{10.1016/j.cma.2004.11.007}. Despite these differences, the methods considered here share a common structure: the intermediate temporal computations are local to each cell, while intercell coupling enters only through conservative numerical fluxes at the end of each time step. We refer to this class as compact single-step DG methods.

\subsection{A brief literature review}

A central challenge for these high-order discretizations is to preserve the physical bounds of the underlying conservation laws. Scalar solutions may satisfy a maximum principle, while admissible states of the compressible Euler equations require positive density and internal energy. High-order polynomial approximations can violate these constraints even when the exact solution remains admissible. Bound-preserving methods therefore seek to enforce a prescribed convex admissible set without sacrificing conservation or high-order accuracy. Below, we briefly review two major classes of bound- and invariant-domain-preserving high-order methods. For a comprehensive review, see \cite{10.48550/arXiv.2512.09116}.

The Zhang--Shu framework for bound-preserving high-order DG methods \cite{10.1016/j.jcp.2009.12.030,10.1016/j.jcp.2010.08.016,10.1098/rspa.2011.0153} follows a \emph{weak-to-strong} construction. A convex cell-average decomposition and bound-preserving first-order building blocks ensure admissibility of the physical cell average after a forward Euler update, under a suitable CFL condition and admissibility constraints on the input polynomial. We refer to this as \emph{weak admissibility}, or equivalently, a \emph{weak invariant-domain-preserving} property. A scaling toward the admissible cell average then enforces the corresponding \emph{strong admissibility} property at prescribed nodes or throughout the cell. Subsequent developments address various admissibility constraints, mesh geometries, and improved CFL conditions \cite{10.1007/s00211-011-0443-7,10.1007/s10915-011-9472-8,10.1016/j.jcp.2018.03.004,10.1137/21M1458247,10.1137/23M1549365,10.1016/j.jcp.2025.113769}. High-order time integration commonly uses strong-stability-preserving (SSP) RK methods, with scaling at the relevant intermediate stages \cite{10.1137/S003614450036757X}. A related stagewise bound-preserving construction has been developed for cRKDG \cite{10.1016/j.jcp.2025.114071}.

Another major family limits high-order solutions or numerical fluxes toward those of a compatible admissibility-preserving low-order discretization. This strategy encompasses flux-corrected DG methods, subcell blending, parametrized flux limiting, and convex limiting, including monolithic formulations \cite{10.1016/bs.hna.2016.08.002, 10.1016/j.jcp.2016.12.031,10.1016/j.cma.2018.11.036,10.1016/j.camwa.2021.02.012,10.1016/j.cma.2021.113876,10.1016/j.jcp.2022.111850,10.1016/j.compfluid.2022.105627,10.1007/s42967-023-00321-6}. Related approaches have been developed for Lax--Wendroff DG and high-order flux-reconstruction methods \cite{10.1007/s10915-016-0291-9,10.1007/s10915-024-02482-9,10.1007/s42967-024-00463-1,10.1007/s10915-025-03118-2}. In particular, the cRKFR method in \cite{10.1007/s10915-025-03118-2} applies low-order-flux-based admissibility limiting to the completed time-averaged update, followed by cellwise scaling.

Although these approaches differ in their limiting mechanisms and temporal implementations, their admissibility arguments rely on suitable low-order building blocks. Several studies have also examined structural properties of completed high-order DG updates and used end-of-step corrections to enforce the desired constraints. A completed-step weak admissibility property has been established for third-order Lax--Wendroff DG methods, with bound preservation enforced through a single scaling limiter \cite{10.1016/j.jcp.2022.111591}. For entropy stability, weak cell entropy inequalities have been established for cRKDG and more general compact single-step DG methods, with strong cell entropy inequalities enforced through a final correction \cite{10.48550/arXiv.2610.18837,10.48550/arXiv.2610.30535}. TVD and TVB properties have also been obtained without requiring TVD or SSP time discretization, by confining the nonlinear correction to the final conservative update \cite{10.48550/arXiv.2606.26692}. This focus on the completed step motivates us to seek a weak admissibility property that supplies the reference states needed for conservative bound preservation.

\subsection{Our contribution}

In this work, we develop a Patch-based Admissibility and Conservative Limiting (PACL) framework based on the aforementioned principles. The key is to seek completed-step weak admissibility over overlapping cell patches, rather than over individual cells. Each mesh vertex defines a patch of incident cells, weighted by the corresponding continuous partition-of-unity hat function. 
Testing the completed DG update with this function and summing over the patch cancels the conservative interface-flux contributions. 
The resulting evolution identity involves only the old DG polynomials and the cell-local time-averaged volume fluxes. 
By estimating these fluxes, we prove that the normalized weighted patch averages remain in the admissible set under a suitable time-step condition. 
We call this property \emph{patchwise weak admissibility}. 
These averages also admit an interpretation on an overlapping dual mesh, with the conserved quantities in each primal cell apportioned among neighboring vertex patches through the partition-of-unity weights.

To recover strong admissibility, we introduce \emph{conservative patch limiting}, consisting of patchwise correction, cell-average recombination, and cellwise Zhang--Shu scaling.
This construction is related to the conservative redistribution framework of Guermond \cite{10.48550/arXiv.2606.22839}, which provides a general approach to conservative postprocessing through local corrections and weighted recombination.
The distinctive feature of PACL is that the completed-step weak-admissibility analysis establishes the feasibility of the patch correction a priori: under the stated CFL conditions, the weighted patch averages of the high-order DG update are admissible on every prescribed patch.
This property simplifies the conservative redistribution to a contraction of the weighted cell states toward their admissible patch average using a common coefficient, preserving the assigned patch mass by construction.
The patch corrections are computed independently from the same uncorrected solution.
A single recombination then yields admissible physical cell averages, with every conserved quantity preserved over the union of the active patches.
Finally, Zhang--Shu scaling preserves the corrected cell averages and enforces admissibility on the prescribed enforcement set in the original DG space.

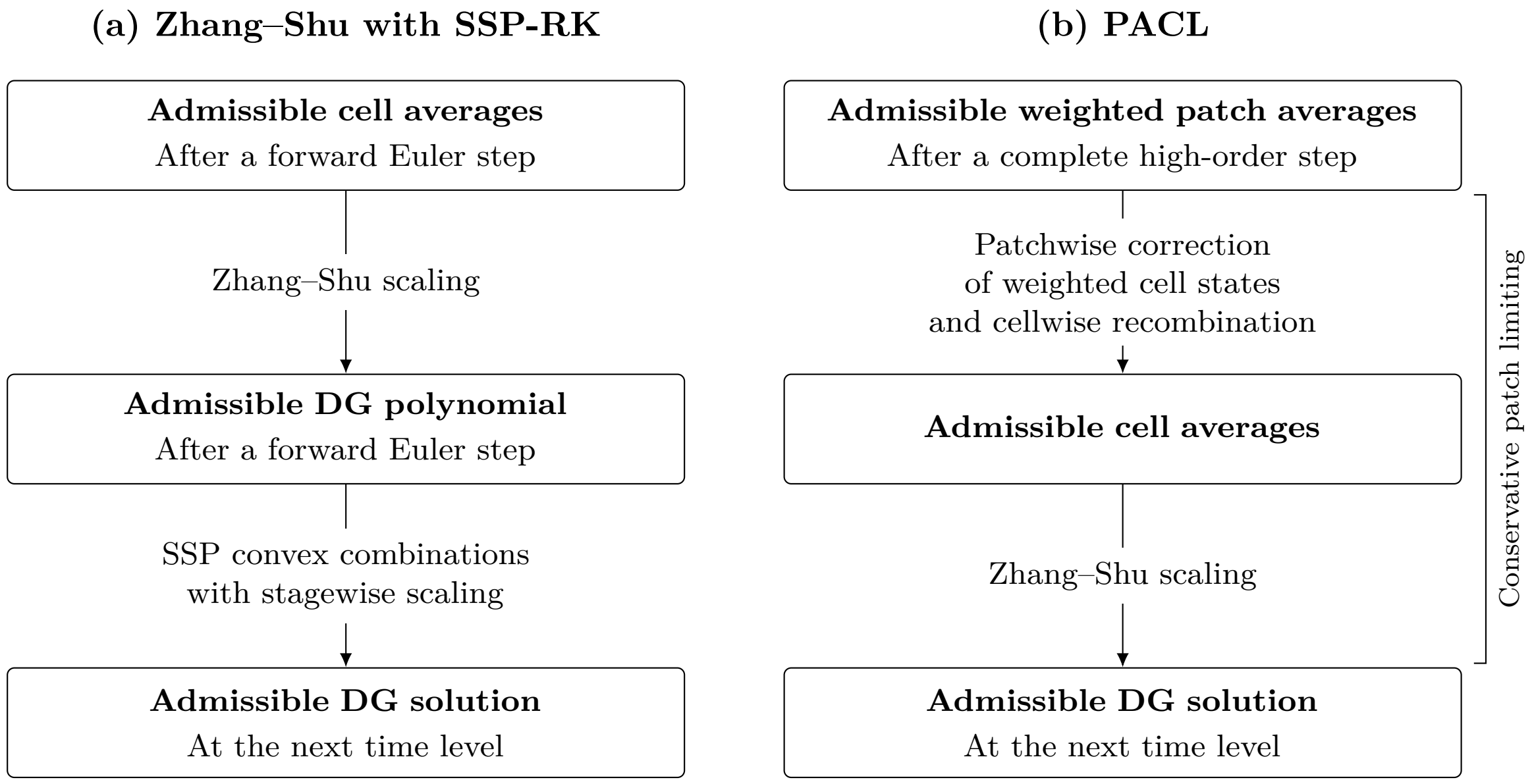
\begin{figure}[!ht]
\centering
\begin{tikzpicture}[
	font=\footnotesize,
	state/.style={
		draw, thin, rounded corners=2pt,
		text width=.405\linewidth,
		minimum height=1.10cm,
		inner sep=5pt, align=center
	},
	flow/.style={-{Latex[length=4pt,width=3.5pt]},thin},
	operation/.style={
		fill=white, align=center,
		text width=.38\linewidth,
		inner xsep=2pt, inner ysep=4pt
	},
	heading/.style={font=\small\bfseries,align=center}
]

\node[state] (L1) at (0,0) {
	\textbf{Admissible cell averages}\\[2pt]
	After a forward Euler step
};
\node[state] (R1) at (.49\linewidth,0) {
	\textbf{Admissible weighted patch averages}\\[2pt]
	After a complete high-order step
};

\node[state] (L2) at (0,-2.95) {
	\textbf{Admissible DG polynomial}\\[2pt]
	After a forward Euler step
};
\node[state] (R2) at (.49\linewidth,-2.95) {
	\textbf{Admissible cell averages}
};

\node[state] (L3) at (0,-5.90) {
	\textbf{Admissible DG solution}\\[2pt]
	At the next time level
};
\node[state] (R3) at (.49\linewidth,-5.90) {
	\textbf{Admissible DG solution}\\[2pt]
	At the next time level
};

\node[heading,anchor=south] at ($(L1.north)+(0,.22cm)$)
	{(a) Zhang--Shu with SSP-RK};
\node[heading,anchor=south] at ($(R1.north)+(0,.22cm)$)
	{(b) PACL};

\draw[flow] (L1.south) -- (L2.north)
	node[midway,operation] {Zhang--Shu scaling};
\draw[flow] (L2.south) -- (L3.north)
	node[midway,operation] {
		SSP convex combinations\\
		with stagewise scaling
	};

\draw[flow] (R1.south) -- (R2.north)
	node[midway,operation] {
		Patchwise correction of weighted cell states\\
and cellwise recombination
	};
\draw[flow] (R2.south) -- (R3.north)
	node[midway,operation] {
		Zhang--Shu scaling
	};

\coordinate (CPLtop) at ($(R1.south east)+(.24cm,-.04cm)$);
\coordinate (CPLbot) at ($(R3.north east)+(.24cm,.04cm)$);
\draw[thin]
	($(CPLtop)+(-.12cm,0)$) -- (CPLtop)
	-- (CPLbot) -- ($(CPLbot)+(-.12cm,0)$);
\node[rotate=90,font=\scriptsize,anchor=north,inner sep=2pt]
	at ($(CPLtop)!0.5!(CPLbot)+(.06cm,0)$)
	{Conservative patch limiting };
\end{tikzpicture}
\caption{Comparison of the Zhang--Shu construction with SSP-RK time
stepping and PACL. The initial bounds follow from bound-preserving
first-order building blocks and a convex cell-average decomposition
(left), or from continuous patch tests and local flux estimates (right).
Conservative patch limiting consists of the two indicated
recovery steps.}
\label{fig:pacl-vs-zhang-shu}
\end{figure}

PACL thus retains the weak-to-strong philosophy of the Zhang--Shu framework, but replaces its cellwise forward Euler property with a patchwise property of the completed high-order step. 
For compressible Euler equations, the forward Euler case recovers the classical characteristic-speed CFL scaling for the admissibility of weighted patch averages; see \Cref{rem:euler-forward-euler}.
The admissibility analysis requires neither an SSP decomposition nor inner-stage bound-preserving scaling. 
It therefore avoids the associated SSP order restrictions and the potential accuracy deterioration caused by repeated stagewise corrections. \Cref{fig:pacl-vs-zhang-shu} summarizes this comparison.
Unlike low-order-flux-based limiting, PACL obtains its admissible reference states from the high-order solution itself, without constructing an auxiliary admissibility-preserving low-order flux or update. 


We establish the framework first for scalar conservation laws. 
A weighted positivity inequality and a local $L^1$ Lipschitz estimate for the time-averaged volume flux yield an abstract patchwise weak maximum principle. 
The required flux estimates are verified for explicit cRKDG, converged and fixed-Picard ADER-DG, and LWDG. 
For cRKDG and ADER-DG, intermediate predictor values may leave the target interval, provided that they remain in a slightly enlarged interval on which the physical flux is well defined and Lipschitz. 
LWDG requires no predictor-range condition, although higher derivatives of the physical flux enter the estimates. 
We also prove that the conservative change in each physical cell average is bounded by the local error of the uncorrected approximation. 
Together with the accuracy estimate for Zhang--Shu scaling, this yields \Cref{thm:pacl-accuracy}: under its local mesh and scaling assumptions, one application of scalar PACL preserves the approximation order, including at active bounds and flat extrema.

For the compressible Euler equations, the analysis must additionally control the nonlinear internal-energy constraint. 
We use geometric quasilinearization (GQL) \cite{10.1137/21M1458247} to express this constraint through a family of linear tests. 
A cellwise linear transformation based on the old cell average removes the local density, velocity, and internal-energy scales before the temporal increment is estimated. 
The resulting conditions permit arbitrarily small absolute density and pressure, while requiring the old pointwise density to remain bounded below relative to its own cell average. 
For cRKDG and fixed-Picard ADER-DG, we establish positive predictor density and admissibility of the completed weighted patch averages, without requiring positive predictor pressure. 
Thus, the local predictor may leave the Euler admissible set through the pressure, while the completed patch averages remain admissible and provide the reference states for conservative patch limiting.

Both the scalar and Euler analyses cover exact spatial integration and compatible quadrature on the conforming meshes and approximation spaces specified below. 
With whole-cell constraints, the exact-integration formulation yields whole-cell admissibility. 
The quadrature formulation retains the physical cell and patch moments, and enforces admissibility at the prescribed sets. 
In either case, limiting is confined to the end of the time step, while the stated time-step conditions ensure predictor well-definedness.
An additional oscillation-control procedure can be incorporated, provided that it preserves the corrected physical cell averages and the local DG space. 
Numerical experiments with cRKDG and ADER-DG examine smooth accuracy, scalar discontinuities, and Euler flows with near-vacuum states and strong shocks.

The main contributions and advantages are summarized as follows.
\begin{itemize}
	\item We establish patchwise weak admissibility for completed high-order-in-time DG updates, without SSP decompositions or inner-stage bound-preserving corrections.
	\item We construct a direct conservative patch limiter, using admissible reference states supplied by the high-order update rather than an auxiliary low-order discretization.
	\item For scalar laws, we cover cRKDG, ADER-DG, and LWDG, and prove accuracy preservation of one limiter application at active bounds and flat extrema.
	\item For Euler, we establish patch admissibility for cRKDG and fixed-Picard ADER-DG, allowing nonpositive predictor pressure while ensuring positive predictor density; see \Cref{rem:predictor-pressure} for the treatment of numerical interface fluxes.
\end{itemize}

\subsection{Organization of the paper}

The paper is organized as follows.
\Cref{sec:preliminaries} introduces the compact methods and their common completed update.
\Cref{sec:framework} establishes scalar patchwise weak admissibility and verifies the local flux estimates, while \Cref{sec:euler} treats the compressible Euler equations. 
\Cref{sec:patch-limiting} presents conservative patch limiting and its scalar accuracy estimate. 
\Cref{sec:quad} develops the compatible-quadrature formulation, and \Cref{sec:numerical} presents numerical experiments.

\section{Preliminaries and compact discretizations}\label{sec:preliminaries}

PACL applies to compact single-step DG methods through their common conservative completed update \eqref{eq:corrector}. The time discretization determines a cell-local time-averaged volume flux, while a conservative numerical flux couples adjacent cells. We formulate this common update and identify the local flux maps for cRKDG, ADER-DG, and  LWDG.

\subsection{Notations and the general compact update}\label{subsec:mesh-space}

Consider the scalar conservation law
\begin{equation*}
	\partial_t u + \nabla\cdot\bfv(u) = 0,
	\quad \bfv(u) = (f_1(u), \ldots, f_d(u)),
	\quad \bx\in\Omega\subset\R^d.
\end{equation*}
Let $\{\Th\}$ be a family of conforming mesh discretizations of $\Omega$. The cells are bounded and connected, with Lipschitz boundaries. Write $h_K=\diam(K)$, let $\abs{K}$ and $\abs{\partial K}$ denote the cell volume and boundary measure, and let $\bn_K$ be the outward unit normal. In one dimension, boundary integration is the sum of the two endpoint values and $\abs{\partial K}=2$.

We assume each cell is the image of a reference simplex or box $\widehat K$ under a $C^1$ bijective map $T_K$ with a $C^1$ inverse. The reference cells belong to a fixed finite collection, and the maps satisfy the following cellwise shape regularity condition,
\begin{equation}\label{eq:mesh-regularity}
	\norm{DT_K}_{L^\infty(\widehat K)} \leq C_{\mathrm{map}} h_K,
	\quad
	\norm{(DT_K)^{-1}}_{L^\infty(\widehat K)}\leq C_{\mathrm{map}}h_K^{-1},
\end{equation}
where $DT_K$ denotes the Jacobian matrix of $T_K$ and
$C_{\mathrm{map}}$ is independent of $K$ and the mesh. The boundary conditions are discussed in \Cref{rem:boundaries}.

Fix $k_{\max}$ and let $k_K$ denote the local polynomial degree on cell $K$. The patchwise admissibility analysis assumes $1\leq k_K\leq k_{\max}$ so that the DG test space contains the restrictions of the vertex hats, while the limiting algorithm will also treat degree-zero cells as described in \Cref{rem:degree-zero-cells}. On a reference simplex \(\widehat{K}\), $\widehat{\mathbb Z}^{k_K} = \Pk^{k_K}(\widehat{K})$ consists of polynomials of total degree at most $k_K$. On a reference box \(\widehat{K}\), $\widehat{\mathbb Z}^{k_K} = \Qk^{k_K}(\widehat{K})$ consists of tensor-product polynomials of degree at most $k_K$ in each coordinate. Set the mapped function space
\[
	\ZK=\{\widehat v\circ T_K^{-1}:\widehat v\in \widehat{\mathbb Z}^{k_K}\}.
\]
These spaces are contained in $W^{1,\infty}(K)\cap C(\overline K)$ and contain constants. Affine simplices recover the physical spaces $\Pk^{k_K}(K)$, while tensor-product and curved elements use the corresponding mapped spaces. Define
\begin{equation*}
	\Vh = \{v_h\in L^2(\Omega):v_h|_K\in\ZK\text{ for every }K\in\Th\}.
\end{equation*}
We use
\[
	(v,w)_K = \int_Kvw\dd\bx, \quad
	\norm{v}_{1,K}=\int_K \abs{v} \dd\bx, \quad
	\norm{v}_{\infty,K}=\operatorname*{ess\,sup}_K \abs{v}.
\]

Write $u_K^n = u_h^n|_K$. On a shared face, superscripts $-$ and $+$ denote the traces from $K$ and its neighbor, respectively. 

A \emph{compact single-step DG method} is represented by a cell-local time-averaged volume flux $\FF_{K,\tau}$. Given $u_h^n\in\Vh$ and $\tau>0$, the raw completed update $u_h^\star$ is defined by
\begin{equation}\label{eq:corrector}
	(u_h^\star-u_h^n,\varphi)_K
	=\tau\left[\int_K\FF_{K,\tau}(u_h^n|_K)\cdot\nabla\varphi\dd\bx - \sum_{e\subset\partial K}\int_e\widehat{\mathscr F}_{K,e}\,\varphi^-\dd s\right],
	\quad \forall\varphi\in\ZK.
\end{equation}
The time-averaged numerical flux is conservative on every shared face,
\begin{equation}\label{eq:corrector-face}
	\widehat{\mathscr F}_{K,e}+\widehat{\mathscr F}_{K',e}=0,
	\quad e=\partial K\cap\partial K',
\end{equation}
and the local volume flux preserves constants,
\begin{equation}\label{eq:constantflux}
	\FF_{K,\tau}(c)=\bfv(c).
\end{equation}
Equation \eqref{eq:corrector} separates the local time evolution from the conservative interface coupling. 
The patch analysis in \Cref{sec:framework} eliminates the interface contribution and uses estimates of the local volume-flux map $\FF_{K,\tau}$. These estimates are established in \Cref{subsec:flux-estimates} for each time discretization. Spatial integration is exact in \Cref{sec:preliminaries,sec:framework,sec:euler,sec:patch-limiting}, while \Cref{sec:quad} develops the compatible quadrature formulation.

\subsection{cRKDG}\label{subsec:crk-predictor}

In cRKDG \cite{10.1137/23M158629X}, cell-local inner stages are followed by an inter-element conservative DG update. For $f_{\bn}(v)=\bfv(v)\cdot\bn$, define the cell-local operator by
\begin{equation*}
	(\loc_K(v),\varphi)_K
	=\int_K\bfv(v)\cdot\nabla\varphi\dd\bx - \int_{\partial K}f_{\bn_K}(v)\varphi\dd s,
	\quad \forall\varphi\in\ZK.
\end{equation*}
The operator $\loc_K$ uses only the polynomial in $K$ and its interior trace, and the divergence theorem gives $\loc_K(c)=0$ for a constant $c$.

For $v_h\in\Vh$, define the conservative DG operator $\DG(v_h)\in\Vh$ by
\begin{equation}\label{eq:DG-operator}
	\begin{aligned}
		(\DG(v_h),\varphi)_K
		&=
		\int_K\bfv(v_h)\cdot\nabla\varphi\dd\bx -
		\sum_{e\subset\partial K}
		\int_e\widehat f(v_h^-,v_h^+;\bn_K)\,
		\varphi^-\dd s,
		\quad \forall\varphi\in\ZK.
	\end{aligned}
\end{equation}
Here the numerical normal flux is consistent and conservative:
\[
	\widehat f(c,c;\bn)=\bfv(c)\cdot\bn,
	\quad
	\widehat f(a,b;\bn)
	=-\widehat f(b,a;-\bn).
\]
Replacing the exterior trace by the interior trace in \eqref{eq:DG-operator} recovers $\loc_K$.

Let $\Ark=(a_{\ell j})_{s \times s}$ be strictly lower triangular and $\brk=(b_1,\ldots,b_s)^\top$ with $\sum_\ell b_\ell=1$ be the RK coefficients. No SSP assumption and no sign restriction on the coefficients is imposed. The resulting admissibility condition therefore imposes no SSP-related restriction on the temporal order. The cell-local stages and completed conservative step are
	\begin{align}
		U_{1,K}(v)&=v,\\
		U_{\ell,K}(v)&=v+\tau\sum_{j<\ell}a_{\ell j}\loc_K(U_{j,K}(v)),\quad \ell=2,\ldots,s,\label{eq:rkstage}\\
		u_h^\star&=u_h^n+\tau\sum_{\ell=1}^s b_\ell\DG(U_{\ell,h}).
	\end{align}
Here $U_{\ell,h}|_K=U_{\ell,K}(u_K^n)$, and the face flux in \eqref{eq:corrector} is the same RK-weighted sum of the stage numerical fluxes. No face communication is needed to construct the inner stages. The corresponding time-averaged volume flux is
\begin{equation}\label{eq:rkflux}
	\FF_{K,\tau}(v)=\sum_{\ell=1}^s b_\ell\bfv(U_{\ell,K}(v)).
\end{equation}
For a constant input, every stage remains constant and $\sum_\ell b_\ell=1$ gives \eqref{eq:constantflux}.

\subsection{ADER-DG with local space-time predictors}\label{sec:ader}

ADER uses a local space-time Galerkin predictor \cite{10.1016/j.jcp.2010.10.024}, rather than RK stages. For a spatial flux field $\vect W$, write its linear weak form, which we regard as a divergence discretization of \(\vect W\), as
\begin{equation}\label{eq:linear-spatial-discrete-div}
	(\mathscr R_K[\vect W],\varphi)_K
	=\int_K\vect W\cdot\nabla\varphi\dd\bx
	-\int_{\partial K}(\vect W\cdot\bn_K)\varphi\dd s,
	\quad \forall\varphi\in\ZK.
\end{equation}
This map contains only weak assembly and the local mass solve. The cRKDG operator above satisfies $\loc_K(v)=\mathscr R_K[\bfv(v)]$. ADER uses the divergence discretization directly in its space-time equations. Define the local space-time polynomial space
\[
	\mathbb{X}_K^r=\Pk^r([0,1];\ZK),
	\quad s=(t-t^n)/\tau,
\]
and use the norms
\[
	\norm{q}_{\mathbb{X},1,K}=\sup_{0\leq s\leq 1}\norm{q(s)}_{1,K},
	\quad
	\norm{q}_{\mathbb{X},\infty,K}=\sup_{0\leq s\leq 1}\norm{q(s)}_{\infty,K}.
\]

\begin{assumption}[Time discretization of the ADER predictor]
	\label{ass:ader-time}
	Fix an integer $r\ge0$ and time quadrature nodes $s_i\in[0,1]$ with positive weights $\vartheta_i$ satisfying $\sum_i \vartheta_i = 1$. The rule is exact through degree $\max\{r,2r-1,0\}$. The temporal pairing
	\[
		\int_0^1(q,\partial_s\phi)_K\dd s
	\]
	is evaluated exactly for all $q,\phi\in \mathbb{X}_K^r$. The same time quadrature rule is used to evaluate the nonlinear term in the local predictor and the volume and interface fluxes in the conservative corrector.
\end{assumption}

A standard choice is $r+1$ Gauss points. Accuracy of a particular ADER scheme additionally requires its usual predictor and corrector quadrature accuracy. With initial polynomial $v$, the local predictor is defined by
\begin{equation}\label{eq:ader-weak}
	\begin{aligned}
		&(q_K(v)(1),\phi(1))_K
		-\int_0^1(q_K(v),\partial_s\phi)_K\dd s\\
		=& (v,\phi(0))_K
		+\tau\sum_i\vartheta_i
		(\mathscr R_K[\bfv(q_K(v)(s_i))],\phi(s_i))_K,  
		\quad \forall\phi\in \mathbb{X}_K^r.
	\end{aligned}
\end{equation}
The initial condition is imposed weakly in time. For sampled $\ZK$-valued $z$, introduce the fixed linear time-solve operator $\mathcal K_r$ by
\begin{equation}\label{eq:timeoperator}
	((\mathcal K_r z)(1),\phi(1))_K
	-\int_0^1(\mathcal K_r z,\partial_s\phi)_K\dd s
	=\sum_i\vartheta_i(z(s_i),\phi(s_i))_K,
	\quad \forall\phi\in \mathbb{X}_K^r.
\end{equation}

The time matrix is nonsingular. \Cref{lem:time-solve} proves this and gives its norm bounds. For a time-dependent flux field, the composition $\mathcal K_r\mathscr R_K$ means that $\mathscr R_K$ is applied at each time node before the time solve.

There are two ways to use this local predictor. The \emph{converged predictor} solves \eqref{eq:ader-weak}, i.e., the nonlinear fixed-point equation
\begin{equation}\label{eq:ader-fixed}
	q_K(v)=v+\tau\mathcal K_r\mathscr R_K[\bfv(q_K(v))],
\end{equation}
whereas the \emph{fixed-Picard predictor} stops after a prescribed number of Picard updates,
\begin{equation}\label{eq:ader-picard}
	q_K^{[0]}(v)=v,
	\quad
	q_K^{[j+1]}(v)=v+\tau\mathcal K_r\mathscr R_K[\bfv(q_K^{[j]}(v))],
	\quad j=0,\ldots,J-1.
\end{equation}
In either case the corrector volume flux is the positive time average
\begin{equation}\label{eq:ader-flux}
	\FF_{K,\tau}^{\mathrm{ADER}}(v)
	=\sum_i\vartheta_i\bfv(q(v)(s_i)),
\end{equation}
with $q=q_K$ for the converged predictor and $q=q_K^{[J]}$ for the fixed-Picard predictor.

\begin{remark}[Picard iterations and formal accuracy]\label{rem:ader-iteration-order}
	For smooth solutions, the deferred-correction interpretation relates the formal accuracy of ADER to the number of Picard updates~\cite[Theorem~4.5]{10.1007/s10915-020-01397-5}. For a target order $p$, a standard sufficient choice is temporal degree $r=p-1$ and $J=p$ predictor updates, together with spatial and quadrature accuracy of at least $p$~\cite[Section~3.2]{10.1007/s42967-023-00290-w}. The initialization $q_K^{[0]}=v$ is not counted as an update, and the conservative corrector is performed afterward. The third-order implementation below uses $r=2$ and $J=3$.
\end{remark}

In both cases, the face flux in \eqref{eq:corrector} is obtained from the same positive time average of a conservative two-state flux evaluated on the predictor traces. A constant input has zero spatial discrete divergence \eqref{eq:linear-spatial-discrete-div} and is fixed by each Picard update. It is also a solution of the converged predictor equation. The latter solution is selected by the local existence and uniqueness result in \Cref{prop:ader-converged}.

\subsection{LWDG}\label{subsec:lwdg}

The LWDG method \cite{10.1016/j.cma.2004.11.007} constructs the time-averaged flux from the old polynomial through a Taylor expansion and CK substitutions, without introducing nonlinear intermediate predictor states.

We use fixed-degree polynomial spaces on affine, shape-regular cells, so the spatial derivatives required by the CK recursion are well defined. For Taylor order $p\geq 1$, assume $\bfv\in C^p$ on a neighborhood of the old-state range, set $V_0=v$, and define the CK coefficients recursively by
\begin{equation}\label{eq:lw-jets}
	\mathtt A_j(v)=
	\left.\frac{\mathrm d^j}{\mathrm ds^j}
	\bfv\!\left(\sum_{\ell=0}^j\frac{s^\ell}{\ell!}V_\ell\right)
	\right|_{s=0},
	\quad
	V_{j+1}=-\nabla\cdot\mathtt A_j(v).
\end{equation}
The local time-averaged volume flux is
\begin{equation}\label{eq:lw-flux}
	\FF_{K,\tau}^{\mathrm LW}(v)
	=\sum_{j=0}^{p-1}\frac{\tau^j}{(j+1)!}\mathtt A_j(v).
\end{equation}
Only old-cell data enter these coefficients. For a constant input, $\mathtt A_0(c)=\bfv(c)$ and $\mathtt A_j(c)=0$ for $j\geq 1$, so constant consistency is automatic.

For example, $V_1=-\nabla\cdot\bfv(v)$ and $\mathtt A_1(v)=-\bfv'(v)\nabla\cdot\bfv(v)$. The first correction already involves a spatial derivative of the old polynomial. \Cref{lem:ck-lipschitz} records the corresponding scaling for all orders.

A conservative numerical flux completes \eqref{eq:corrector}. Unlike the two preceding constructions, this volume flux requires no nonlinear intermediate predictor state. Its analysis instead uses the higher flux smoothness stated in \Cref{lem:ck-lipschitz}.

\begin{remark}[Prediction in primitive variables]\label{rem:predictor-in-primitive-var}
	A local predictor may evolve primitive variables while the corrector updates the conserved variables. 
	Examples include the ADER finite volume and DG schemes of Zanotti and Dumbser~\cite{10.1186/s40668-015-0014-x} and the locally implicit LWDG scheme of Felton et al.~\cite{10.48550/arXiv.1806.06756}.
	The primitive-variable MUSCL--Hancock scheme predicts interface states at the half time step and uses them to update the conserved cell averages~\cite{10.1007/978-3-642-60543-7_6}.
\end{remark}

\section{Patchwise weak admissibility for scalar conservation laws}\label{sec:framework}

The scalar PACL analysis combines the patch evolution identity
\eqref{eq:weighted-identity} with the monotonicity estimate
\eqref{eq:comparison}.
Continuous patch hats cancel conservative interface-flux contributions,
as shown in \eqref{eq:cancel-flux}.
The weighted positivity bound implied by \eqref{eq:eta} controls
the old-state contribution, while the local Lipschitz estimate
\eqref{eq:weighted-lipschitz} controls the corresponding difference
in the time-averaged flux contribution.

\subsection{Continuous patch hats and cancellation}

\begin{assumption}[Embedded nonnegative partition of unity]\label{ass:pou}
	There is a finite family of $\{\psi_a\}_{a\in\Ah}$, which we call \emph{vertex hat functions}, such that
	\begin{equation*}
		\psi_a\in\Vh\cap C^0(\overline\Omega),\quad
		\psi_a\ge0,\quad
		\sum_{a\in\Ah}\psi_a=1.
	\end{equation*}
	Here $\mathcal{A}_h$ is the collection of patch indices --- typically vertices for conforming meshes. The hat traces agree on periodically identified faces.
\end{assumption}

\begin{remark}[Boundary conditions]\label{rem:boundaries}
	The patch admissibility analysis below is stated for periodic boundaries or for patches whose hat function vanishes on the physical boundary. \Cref{app:physical-boundaries} gives the corresponding admissibility conditions for boundary patches adjacent to symmetry boundaries, extrapolation outflow boundaries, or boundaries with prescribed exterior states.
\end{remark}

\begin{figure}[!ht]
	\centering
	\input{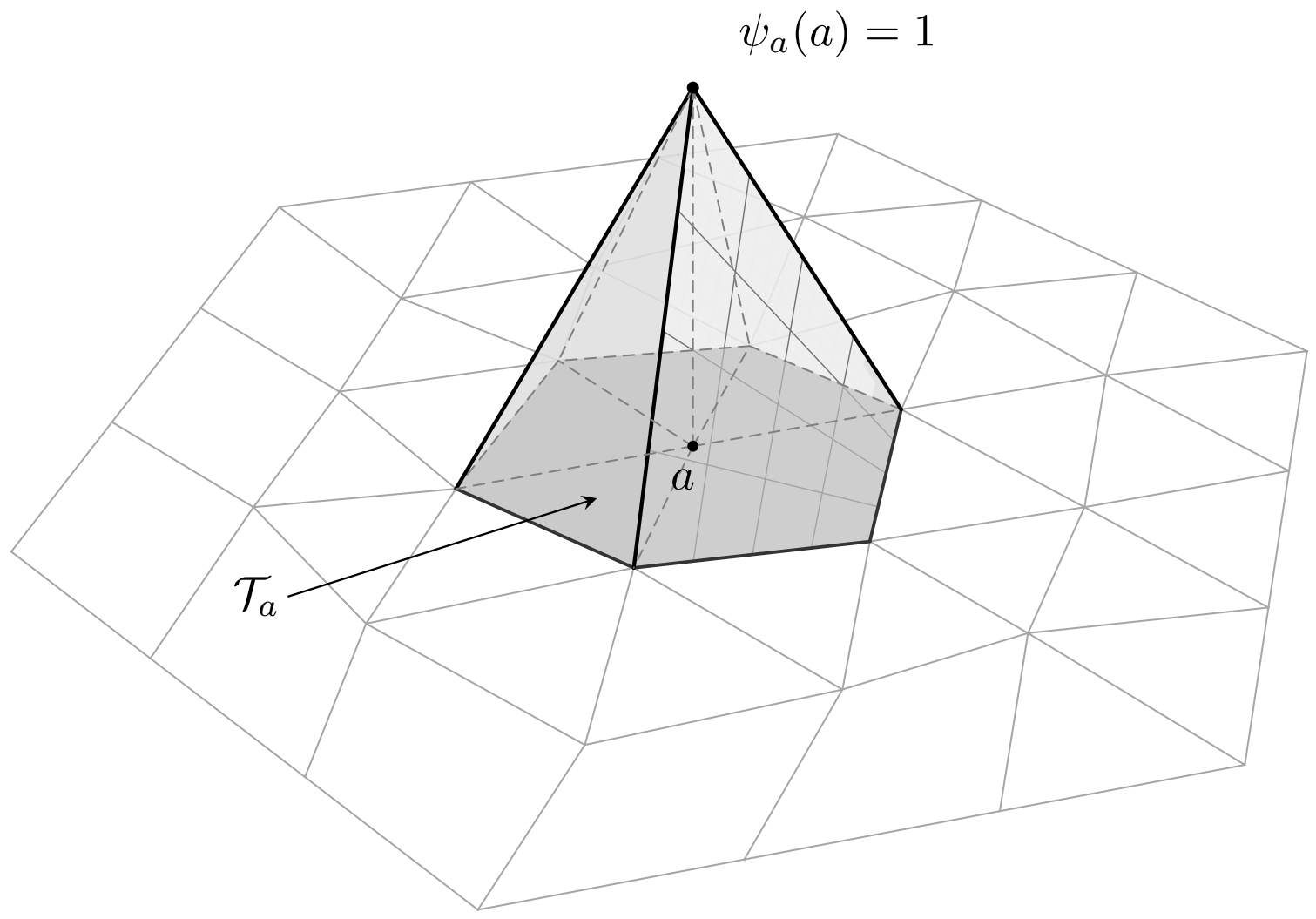}
	\caption{A vertex patch $\mathcal T_a$ and its continuous nodal weight $\psi_a$ on a mixed triangular--quadrilateral mesh. The shaded footprint consists of four triangles and one quadrilateral sharing $a$. On each triangle, $\psi_a$ is the barycentric coordinate associated with $a$. On the quadrilateral, $\psi_a\circ F_K=(1-\xi)(1-\eta)$, where $F_K$ is the bilinear map from $[0,1]^2$ with $F_K(0,0)=a$. The thin lines on the quadrilateral surface show parameter lines, rather than mesh edges. These pieces have matching linear traces on shared edges. The normalized patch average integrates $v_h$ against $\psi_a/\mathfrak m_a$, whose integral over the patch equals one.}
	\label{fig:patch-weighted-average}
\end{figure}

For each vertex $a$, define the patch and the associated masses by
\begin{equation}\label{eq:patch-masses}
	\mathcal T_a=\{K\in\Th:\psi_a|_K\not\equiv0\},
	\quad
	m_{aK}=\int_K\psi_a\dd\bx,
	\quad
	\mathfrak m_a=\sum_{K\in\mathcal T_a}m_{aK}>0.
\end{equation}
The weighted mass and normalized patch average are
\begin{equation}\label{eq:mass}
	\calM_a(v_h)=\sum_{K\in\mathcal T_a}\int_K\psi_av_h\dd\bx,
	\quad
	\avg{v_h}_a=\frac{\calM_a(v_h)}{\mathfrak m_a}.
\end{equation}
See \Cref{fig:patch-weighted-average} for an illustration on a 2D hybrid triangular-quadrilateral mesh.

\begin{remark}[Dual-cell interpretation]\label{rem:dual-cell}
	In a vertex-centered dual mesh, the indicator function of a dual cell selects the volume and the conserved quantity assigned to its vertex. 
	Here the hat function $\psi_a$ plays the same role through
	a continuous weight. At each point $\bx\in K$, the fraction assigned to vertex $a$ is $\psi_a(\bx)$, and these fractions sum to one. 
	Thus $m_{aK}$ is the volume contributed by $K$ to vertex $a$, and $\mathfrak m_a$ is the total effective volume associated with $a$. Weighted moments of the form $M_a(v_h)$ also arise as generalized
conserved quantities in local-conservation formulations of continuous
Galerkin methods~\cite{10.1016/j.compfluid.2017.06.018}.
\end{remark}

\begin{definition}[Patchwise weak admissibility]\label{def:patchwise-admissibility}
	For prescribed bounds $m<M$, a raw approximation $w_h\in\Vh$ is \emph{patchwise weakly admissible} with respect to $[m,M]$ if
	\[
		\avg{w_h}_a\in[m,M]\quad\text{for every patch }a.
	\]
	For Euler equations, the analogous condition is $\avg{\vect W_h}_a\in\GE$, where the convex admissible set $\GE$ is defined in \Cref{sec:euler}.
\end{definition}
This condition permits violations of the prescribed bounds by the raw polynomial and its cell averages. The conservative correction in \Cref{sec:patch-limiting} uses the admissible patch averages to recover admissible cell averages.

Because the weights form a partition of unity, we have
\begin{equation*}
	\sum_{a:K\in\mathcal T_a} m_{aK}=\abs{K},
	\quad
	\sum_a\calM_a(v_h)=\int_\Omega v_h\dd\bx.
\end{equation*}

\begin{lemma}[Patch evolution identity]\label{lem:cancellation}
	For the meshes and spaces of \Cref{subsec:mesh-space}, suppose \Cref{ass:pou} holds. Let the completed update \eqref{eq:corrector} have a finite conservative face flux satisfying \eqref{eq:corrector-face}. Then
	\begin{equation}\label{eq:weighted-identity}
		\calM_a(u_h^\star)=\calM_a(u_h^n)
		+\tau\sum_{K\in\mathcal T_a}
		\int_K\FF_{K,\tau}(u_K^n)\cdot\nabla\psi_a\dd\bx.
	\end{equation}
	If the local volume flux additionally satisfies \eqref{eq:constantflux}, then, for every constant state $c$,
	\begin{equation}\label{eq:constantweighted}
		\sum_{K\in\mathcal T_a}\int_K\FF_{K,\tau}(c)\cdot\nabla\psi_a\dd\bx=0.
	\end{equation}
\end{lemma}

\begin{proof}
	Choose $\varphi=\psi_a|_K$ in \eqref{eq:corrector} and sum over all cells in the patch. On an interior face, the two contributions equal
	\begin{equation}\label{eq:cancel-flux}
		-\tau\int_e\widehat{\mathscr F}_{K,e}(\psi_a^- -\psi_a^+)\dd s=0,
	\end{equation}
	by conservation of the numerical flux and continuity of $\psi_a$. Periodic boundary contributions cancel in the same manner. This proves \eqref{eq:weighted-identity}. For a constant state, use \eqref{eq:constantflux} and the divergence theorem to obtain
	\[
		\sum_K\int_K\bfv(c)\cdot\nabla\psi_a\dd\bx=0.
	\]
\end{proof}

For fixed old data and local volume-flux maps, \eqref{eq:weighted-identity} gives the same patch averages for every finite conservative interface flux. The admissibility analysis therefore concerns the cell-local volume fluxes, while the interface flux continues to determine the remaining components of the raw DG update.

\subsection{One-dimensional illustration on a uniform mesh}\label{subsec:1d-illustration}

Take a uniform periodic mesh $I_j=[x_{j-1/2},x_{j+1/2}]$ with common cell length $h$ and polynomial degree $k\geq 1$. Write $X_j(\xi)=x_{j-1/2}+h\xi$ and $v_j=v_h\circ X_j$ for $0\leq\xi\leq 1$. The hat at the shared endpoint $x_{j+1/2}$ is
\[
	\psi_{j+1/2}(x)=
	\begin{cases}
		(x-x_{j-1/2})/h, & x\in I_j,\\
		(x_{j+3/2}-x)/h, & x\in I_{j+1},\\
		0, & \text{otherwise}.
	\end{cases}
\]
Each weighted cell mass equals $h/2$, and the patch mass is $h$. Consequently its normalized average is simply
\begin{equation*}
	\avg{v_h}_{j+1/2}
	=\int_0^1\xi v_j(\xi)\dd\xi
	+\int_0^1(1-\xi)v_{j+1}(\xi)\dd\xi.
\end{equation*}
The two terms together return $c$ for a constant input $v_h=c$. For a reference-cell input $v$, define the averaged local flux 
\[
	F_{j,\tau}(v)=\frac{1}{h}\int_{I_j}
	\FF_{I_j,\tau}(v\circ X_j^{-1})(x)\dd x.
\]
Testing the two cell equations with the two pieces of the hat cancels the numerical flux at their shared endpoint. The two outer face terms vanish. Since the hat derivatives are $1/h$ and $-1/h$, the result is
\begin{equation}\label{eq:1d-hat-evolution}
	\avg{u_h^\star}_{j+1/2}
	=\int_0^1\xi v_j\dd\xi+\int_0^1(1-\xi)v_{j+1}\dd\xi
	+\frac\tau h\left(F_{j,\tau}(v_j)-F_{j+1,\tau}(v_{j+1})\right),
\end{equation}
where now $v_j=u_h^n\circ X_j$. Thus only the local averaged fluxes remain, irrespective of which method in \Cref{sec:preliminaries} produced them.

Equation~\eqref{eq:1d-hat-evolution} illustrates the comparison
argument used in the general proof below. Its right-hand side
expresses the completed patch average as a function of the two
old-cell polynomials. We seek a time-step condition under which
increasing either polynomial pointwise cannot decrease this average.
Since constant inputs are preserved, comparison with the constant
states $m$ and $M$ then gives the desired bounds whenever
$m\leq v_j,v_{j+1}\leq M$.

In the one-dimensional formula, increasing $v_j$ by a nonnegative
degree-$k$ polynomial $p\ge0$, while keeping $v_j+p$ within the prescribed bounds,
raises the weighted old-state term by $\int_0^1\xi p\,\dd\xi$.
Note that there exists a constant $\eta_k>0$, depending only on $k$,
such that every polynomial $p$ of degree at most $k$ that is
nonnegative on $[0,1]$ satisfies
\begin{equation*}
    \int_0^1 \xi p(\xi)\dd\xi
    \geq \eta_k\norm{p}_{L^1(0,1)}.
\end{equation*}
If the local flux maps are $L^1$-Lipschitz on polynomials taking
values in $[m,M]$, with constants $L_j(\tau)$ and $L_{j+1}(\tau)$,
the corresponding left-cell flux contribution can decrease by
at most $(\tau/h)L_j(\tau)\norm{p}_{L^1(0,1)}$.
Reflection gives the same weighted estimate with $1-\xi$,
so the right-cell contribution is treated similarly.
Thus, the condition
$\frac{\tau}{h}
    \max\{L_j(\tau),L_{j+1}(\tau)\}
    \leq \eta_k$ 
ensures the required monotonicity: the flux contribution cannot
decrease by more than the weighted old-state contribution increases.
The next subsection formulates these weighted positivity and
local flux requirements on general conforming meshes and combines
them to prove scalar patchwise weak admissibility.

\subsection{Scalar patchwise weak admissibility}

Fix target bounds $m<M$ contained in an open interval on which the physical flux is locally Lipschitz. For each active pair $(a,K)$, meaning $K\in\mathcal T_a$
or, equivalently, $\psi_a|_K\not\equiv0$, define the
weighted positivity constant
\begin{equation}\label{eq:eta}
	\eta_{aK}:=
	\inf_{\substack{0\ne p\in\ZK\\ p\ge0\text{ on }K}}
	\frac{\int_K\psi_ap\dd\bx}{\norm{p}_{1,K}}.
\end{equation}

\begin{proposition}[Positive patch weights and standard constructions]\label{prop:patch-positivity}
	Let the cells and local spaces be as in \Cref{subsec:mesh-space}. 
	Under \Cref{ass:pou}, the constant $\eta_{aK}$ in \eqref{eq:eta} is positive for every active pair $(a,K)$.

	For $k_K\geq1$, standard vertex hats satisfy \Cref{ass:pou}: barycentric coordinates on conforming simplicial meshes, multilinear vertex functions on affine tensor-product meshes, and their pullbacks on mapped cells with compatible face parameterizations and the corresponding mapped polynomial spaces.
	With the prescribed boundary identifications, these functions form a continuous nonnegative partition of unity.
	Under the degree bound and geometric regularity \eqref{eq:mesh-regularity}, there are constants $\eta_*>0$ and $C_\psi>0$, independent of $h$, such that
	\[
		\eta_{aK}\geq\eta_*,
		\quad
		h_K\norm{\nabla\psi_a}_{\infty,K}\leq C_\psi
	\]
	for every active pair.
\end{proposition}

\begin{proof}
	Consider $\mathcal P_K=\{p\in\ZK:p\ge0,\ \norm{p}_{1,K}=1\}$. This set is closed and bounded in a finite-dimensional space, hence compact. By the Weierstrass extreme value theorem, the continuous functional $p\mapsto\int_K\psi_a p$ attains its minimum on $\mathcal P_K$. Suppose its minimum were zero. Both factors in the integrand are nonnegative and continuous, so the minimizing polynomial would vanish wherever $\psi_a>0$. Since the pair is active, that set contains a nonempty open subset of $K$. A polynomial vanishing on an open set is identically zero, contrary to its normalization. For a mapped polynomial the same argument applies after pullback to the reference cell. 

	On an affine simplex the barycentric coordinates are nonnegative, sum to one, and agree on shared faces. They belong to $\Pk^{k_K}$ for $k_K\geq 1$. On the reference box the corresponding functions are $\prod_{r=1}^d\xi_r^{\epsilon_r}(1-\xi_r)^{1-\epsilon_r}$, $\epsilon_r\in\{0,1\}$, and belong to $\Qk^{k_K}$ for $k_K\geq 1$. Compatible element maps preserve nonnegativity, the partition identity, and matching traces.

	For the uniform estimate let $T_K:\widehat K\to K$ be the element map and $J_K=\abs{\det DT_K}$. With hats denoting pullbacks,
	\[
		\frac{\int_K\psi_ap}{\int_Kp}
		\geq\frac{J_{K,\min}}{J_{K,\max}}
		\frac{\int_{\widehat K}\widehat\psi_a\widehat p}
		{\int_{\widehat K}\widehat p}
		\geq\frac{J_{K,\min}}{J_{K,\max}}\widehat\eta_a.
	\]
	The reference constant $\widehat\eta_a$ is positive by the first part of the proof. The reference cells and degrees range over finite collections. Condition \eqref{eq:mesh-regularity} bounds both $J_{K,\max}/J_{K,\min}$ and $h_K\norm{(DT_K)^{-1}}_\infty$, which proves the two asserted controls, using $\nabla\psi_a\circ T_K=(DT_K)^{-\top}\widehat\nabla\widehat\psi_a$.
\end{proof}

\begin{remark}[Spatial degree and patch masses]
	These constructions embed continuous $\Pk^1$ or mapped $\Qk^1$ vertex functions in the DG test space. Thus, direct testing by vertex hats is available on cells of degree at least one. For cells of degree zero, \Cref{rem:degree-zero-cells} constructs the required weighted moments from the conservative update. On a connected conforming mesh, a continuous piecewise-constant test function is globally constant and does not provide the vertex patches used here.
	
	On standard affine cells the patch masses are given in \Cref{tab:patch-mass-standard}. For example, on a triangle the corrected cell average in \Cref{sec:patch-limiting} is the arithmetic average of three corrected weighted states.
	\begin{table}[htb]
		\centering
		\caption{Table of patch masses of standard affine cells.}
		\label{tab:patch-mass-standard}
		\begin{tabular}{lll}
			\toprule
			cell & local patch weight & $m_{aK}/\abs{K}$\\
			\midrule
			interval & linear endpoint hat & $1/2$\\
			$d$-simplex & barycentric coordinate & $1/(d+1)$\\
			$d$-box & multilinear vertex hat & $1/2^d$\\
			\bottomrule
		\end{tabular}
	\end{table}
\end{remark}


For later use define the scalar functional
\begin{equation}\label{eq:weighted-flux-functional}
	\mathscr F_{aK,\tau}(v)
	=\int_K\FF_{K,\tau}(v)\cdot\nabla\psi_a\dd\bx.
\end{equation}
The time discretization enters the patch argument only through a local Lipschitz estimate for this functional.

\begin{assumption}[Local weighted-flux Lipschitz continuity]\label{ass:weighted-lipschitz}
	For every active pair $(a,K)$ there is $L_{aK}(\tau)\ge0$ such that
	\begin{equation}\label{eq:weighted-lipschitz}
		\abs{\mathscr F_{aK,\tau}(v)-\mathscr F_{aK,\tau}(w)}
		\leq L_{aK}(\tau)\norm{v-w}_{1,K}
	\end{equation}
	for every $v,w\in\ZK$ satisfying $m\leq v,w\leq M$ on $K$. The local flux map is assumed to be well defined for every such comparison input, not only for the particular old solution.
\end{assumption}

A sufficient condition for \Cref{ass:weighted-lipschitz} is
\begin{equation}\label{eq:volume-lipschitz}
	\norm{\FF_{K,\tau}(v)-\FF_{K,\tau}(w)}_{1,K}
	\leq L_K^{\mathrm{vol}}(\tau)\norm{v-w}_{1,K},
\end{equation}
for then one may take $L_{aK}(\tau)=\norm{\nabla\psi_a}_{\infty,K}L_K^{\mathrm{vol}}(\tau)$.

\begin{theorem}[Scalar patchwise weak maximum principle]\label{thm:main}
	For the meshes and spaces of \Cref{subsec:mesh-space}, consider \eqref{eq:corrector} under \Cref{ass:pou,ass:weighted-lipschitz}. Require the finite conservative face flux \eqref{eq:corrector-face}, the volume-flux consistency \eqref{eq:constantflux}, and $\eta_{aK}>0$ for every active pair $(a, K)$. If the old polynomial satisfies $m\leq u_h^n\leq M$ cellwise and the CFL condition holds
	\begin{equation}\label{eq:scalar-patch-cfl}
		\tau L_{aK}(\tau)\leq\eta_{aK}
		\quad\text{for every active pair }(a,K),
	\end{equation}
	then the completed update is patchwise weakly admissible:
	\begin{equation}\label{eq:scalar-patch-conclusion}
		m\leq\avg{u_h^\star}_a\leq M
		\quad\text{for every }a.
	\end{equation}
	No intermediate-stage maximum principle is assumed.
\end{theorem}
\begin{proof}
	For a collection $\{v_K\}_{K\in\mathcal T_a}$ with $m\leq v_K\leq M$, define
	\[
		H_a(\{v_K\})=
		\frac{1}{\mathfrak m_a}\sum_{K\in\mathcal T_a}
		\left(\int_K\psi_av_K\dd\bx+\tau\mathscr F_{aK,\tau}(v_K)\right).
	\]
	Replace one input $v_K$ by $w_K\geq v_K$ and set $p=w_K-v_K\ge0$. Then
	\begin{equation}\label{eq:comparison}
		\begin{aligned}
			\mathfrak m_a\left(H_a(\ldots,w_K,\ldots)-H_a(\ldots,v_K,\ldots)\right)
			&=\int_K\psi_ap\dd\bx
			+\tau\left(\mathscr F_{aK,\tau}(w_K)-\mathscr F_{aK,\tau}(v_K)\right)\\
			&\geq(\eta_{aK}-\tau L_{aK}(\tau))\norm{p}_{1,K}\ge0.
		\end{aligned}
	\end{equation}
	Thus $H_a$ is monotone in the pointwise order of each cell polynomial. Constant preservation and \eqref{eq:constantweighted} give $H_a(\{c\})=c$. Comparing the old input collection with the constant collections $m$ and $M$, and using \Cref{lem:cancellation}, yields \eqref{eq:scalar-patch-conclusion}.
\end{proof}

\Cref{thm:main} separates the geometric conditions on the patch weights from the local flux estimate. \Cref{prop:patch-positivity} provides standard partitions with uniform geometric bounds. \Cref{subsec:flux-estimates} verifies \Cref{ass:weighted-lipschitz} for the methods in \Cref{sec:preliminaries}.

\subsection{Verification of the local flux estimates}\label{subsec:flux-estimates}

We now verify the local flux estimate required by \Cref{thm:main}. Common spatial discrete divergence bounds reduce the analysis to the dependence of each local predictor on its initial polynomial. The time discretization enters through the resulting flux majorant and predictor-legality condition.

\begin{assumption}[Scalar flux interval]\label{ass:scalar-data}
	Fix $m<M$ and $\widetilde m<m<M<\widetilde M$ such that $[\widetilde m,\widetilde M]$ lies in the domain of $\bfv$. The physical flux is Lipschitz on this enlarged interval with constant $\alpha$.
\end{assumption}
Set
\[
	R=(M-m)/2,
	\quad
	\varepsilon_I=\min\{m-\widetilde m,\widetilde M-M\}>0.
\]
The enlarged interval permits cRKDG stages and ADER predictors to leave $[m,M]$ while keeping their flux evaluations well defined. The Taylor-coefficient estimates instead use the stronger smoothness and affine-cell hypotheses stated in \Cref{lem:ck-lipschitz}.

\begin{lemma}[Linear discrete divergence and scalar flux estimates]\label{lem:GLip}
	For the meshes and polynomial spaces of \Cref{subsec:mesh-space}, let $\mathscr R_K$ be defined by \eqref{eq:linear-spatial-discrete-div}. There are constants $C_{R,K}$, $C_{k_K}$, and $C_{k_K}^{(\infty)}$, bounded uniformly in terms of the reference cells, $k_{\max}$, and $C_{\mathrm{map}}$, with the following properties.

	\begin{enumerate}
		\renewcommand{\labelenumi}{\normalfont(\roman{enumi})}
		\renewcommand{\theenumi}{(\roman{enumi})}
		\item\label[claim]{part:GLip-discrete-div} For every continuous flux field $\vect W\in C(\overline K)^d$,
		\begin{align}
			\norm{\mathscr R_K[\vect W]}_{1,K}
			&\leq\frac{C_{R,K}}{h_K}
			\left(\norm{\vect W}_{L^1(K)}
			+h_K\norm{\vect W\cdot\bn_K}_{L^1(\partial K)}\right),
			\label{eq:discrete-div-L1-detail}\\
			\norm{\mathscr R_K[\vect W]}_{\infty,K}
			&\leq\frac{C_{R,K}}{h_K}\norm{\vect W}_{L^\infty(K)}.
			\label{eq:discrete-div-infinity-detail}
		\end{align}

		\item\label[claim]{part:GLip-flux} If \Cref{ass:scalar-data} also holds, then, for $v,w\in\ZK$ taking values in $[\widetilde m,\widetilde M]$ on $\overline K$,
		\begin{equation}\label{eq:GLip}
			\begin{aligned}
				\norm{\mathscr R_K[\bfv(v)-\bfv(w)]}_{1,K}
				&\leq \frac{\alpha C_{k_K}}{h_K}\norm{v-w}_{1,K},\\
				\norm{\mathscr R_K[\bfv(v)-\bfv(w)]}_{\infty,K}
				&\leq \frac{\alpha C_{k_K}^{(\infty)}}{h_K}\norm{v-w}_{\infty,K}.
			\end{aligned}
		\end{equation}
	\end{enumerate}
\end{lemma}
\begin{proof} 
	First, \eqref{eq:mesh-regularity} gives $\abs{\det DT_K}\asymp h_K^d$, $\abs{K}\asymp h_K^d$, and $\abs{\partial K}\leq Ch_K^{d-1}$. Here, we write $A\asymp B$ if $cB\leq A\leq CB$ for positive
constants $c$ and $C$ independent of the cell and the mesh size. Norm equivalence on the fixed reference polynomial spaces, followed by change of variables and the chain rule, therefore yields, for every $z\in\ZK$,
	\begin{equation}\label{eq:local-inverse-trace}
		\begin{aligned}
			\abs{K}^{-1/r}\norm{z}_{r,K} \leq \norm{z}_{\infty,K}&\leq C\abs{K}^{-1/r}\norm{z}_{r,K},\quad r=1,2,\\
			h_K\norm{\nabla z}_{\infty,K}&\leq C\norm{z}_{\infty,K},\quad
			h_K\norm{z}_{1,\partial K}\leq C\norm{z}_{1,K}.
		\end{aligned}
	\end{equation}
	The constants are uniform because the reference cells and degrees range over finite collections and the element maps satisfy the same scaled bounds. Let $\Pi_K$ be the physical $L^2$ projection onto $\ZK$. Its $L^2$ contraction property and \eqref{eq:local-inverse-trace} give, for $g\in L^\infty(K)$,
	\[
		\norm{\Pi_Kg}_{\infty,K}
		\leq C\abs{K}^{-1/2}\norm{\Pi_Kg}_{2,K}
		\leq C\abs{K}^{-1/2}\norm{g}_{2,K}
		\leq C\norm{g}_{\infty,K}.
	\]
	The gradient estimate and continuity up to the boundary then imply
	\begin{equation}\label{eq:projection-bound}
		h_K\norm{\nabla\Pi_Kg}_{\infty,K}
		+\norm{\Pi_Kg}_{\infty,\partial K}\leq C\norm{g}_{\infty,K}.
	\end{equation}
	For $\norm{g}_{\infty,K}\leq 1$, orthogonality and $\mathscr R_K[\vect W]\in\ZK$ give
	\[
		\begin{aligned}
			\abs{(\mathscr R_K[\vect W],g)_K}
			&=\abs{(\mathscr R_K[\vect W],\Pi_Kg)_K}\\
			&\leq\norm{\vect W}_{L^1(K)}\norm{\nabla\Pi_Kg}_{\infty,K}
			+\norm{\vect W\cdot\bn_K}_{L^1(\partial K)}
			\norm{\Pi_Kg}_{L^\infty(\partial K)}.
		\end{aligned}
	\]
	Estimate \eqref{eq:projection-bound} bounds the two test factors by $C/h_K$ and $C$. Taking the supremum over $g$ proves \eqref{eq:discrete-div-L1-detail}. The inverse bound \eqref{eq:local-inverse-trace} and $h_K\abs{\partial K}/\abs{K}\leq C$ then yield
	\[
		\begin{aligned}
			\norm{\mathscr R_K[\vect W]}_{\infty,K}
			&\leq\frac{C}{\abs{K}}\norm{\mathscr R_K[\vect W]}_{1,K}\leq\frac{C}{h_K}
			\left(1+\frac{h_K\abs{\partial K}}{\abs{K}}\right)
			\norm{\vect W}_{C(\overline K)},
		\end{aligned}
	\]
	which gives \eqref{eq:discrete-div-infinity-detail} after enlarging $C_{R,K}$.

	For \Cref{part:GLip-flux}, put $p=v-w$. Physical-flux Lipschitz continuity and the polynomial trace inequality imply
	\[
		\norm{\bfv(v)-\bfv(w)}_{L^1(K)}
		+h_K\norm{(\bfv(v)-\bfv(w))\cdot\bn_K}_{L^1(\partial K)}
		\leq C\alpha\norm{p}_{1,K}.
	\]
	Combining this inequality with \eqref{eq:discrete-div-L1-detail} proves the first estimate in \eqref{eq:GLip}. The second follows from \eqref{eq:discrete-div-infinity-detail} and $\norm{\bfv(v)-\bfv(w)}_{C(\overline K)}\leq\alpha\norm{v-w}_{\infty,K}$. Polynomial inverse and trace estimates are applied to $p$ and to the discrete divergence \eqref{eq:linear-spatial-discrete-div}, not to the general nonpolynomial flux difference.
\end{proof}

\subsubsection{cRKDG}\label{subsub:crk-estimates}

For the cRKDG method in \Cref{subsec:crk-predictor}, $\loc_K(v)-\loc_K(w)=\mathscr R_K[\bfv(v)-\bfv(w)]$. \Cref{lem:GLip} therefore bounds the stage operator. Define the stage majorants
\begin{equation}\label{eq:Pell}
	P_1(z)=1,
	\quad
	P_\ell(z)=1+z\sum_{j<\ell}\abs{a_{\ell j}}P_j(z),
	\quad
	B_{\Ark,\brk}(z)=\sum_{\ell=1}^s\abs{b_\ell}P_\ell(z).
\end{equation}

\begin{proposition}[Patchwise weak admissibility for cRKDG]\label{prop:rkLip}
	For the meshes and spaces of \Cref{subsec:mesh-space}, suppose \Cref{ass:scalar-data} holds. Use \eqref{eq:rkstage}--\eqref{eq:rkflux} with a strictly lower triangular $\Ark$ and $\sum_\ell b_\ell=1$. Let $v,w\in\ZK$ satisfy $m\leq v,w\leq M$ on $\overline K$, and set
	\[
		z_{K,1}=\frac{\alpha\tau C_{k_K}}{h_K},
		\quad
		z_{K,\infty}=\frac{\alpha\tau C_{k_K}^{(\infty)}}{h_K}.
	\]
	\begin{enumerate}
		\renewcommand{\labelenumi}{\normalfont(\roman{enumi})}
		\renewcommand{\theenumi}{(\roman{enumi})}
		\item\label[claim]{part:rk-local} \emph{Local estimates.} If the time step restriction holds 
		\begin{equation}\label{eq:rk-tubeCFL}
			R\max_\ell\{P_\ell(z_{K,\infty})-1\}\leq\varepsilon_I,
		\end{equation}
		then every comparison input in $[m,M]$ generates stages in $[\widetilde m,\widetilde M]$, and
		\begin{equation}\label{eq:rkLip}
			\norm{\FF_{K,\tau}(v)-\FF_{K,\tau}(w)}_{1,K}
			\leq \alpha B_{\Ark,\brk}(z_{K,1})\norm{v-w}_{1,K}.
		\end{equation}
		\item\label[claim]{part:rk-admissibility} \emph{Completed-step admissibility.} Also require \Cref{ass:pou}, $m\leq u_h^n\leq M$ cellwise, $\eta_{aK}>0$ for every active pair, and a finite conservative face flux \eqref{eq:corrector-face}. Then the raw cRKDG update is patchwise weakly admissible if, for every active pair, the CFL condition holds
		\begin{equation}\label{eq:rkCFL}
			\frac{\alpha\tau}{h_K}
			\left(h_K\norm{\nabla\psi_a}_{\infty,K}\right)
			B_{\Ark,\brk}\!\left(\frac{\alpha\tau C_{k_K}}{h_K}\right)
			\leq\eta_{aK},
		\end{equation}
		together with \eqref{eq:rk-tubeCFL}. For a globally Lipschitz physical flux, the interval-containment condition is unnecessary.
	\end{enumerate}
\end{proposition}

\begin{proof}
	\emph{Predictor containment.} Fix $v$ with $m\leq v\leq M$, and set $c_0=(m+M)/2$. For the first stage, $U_{1,K}(v)=v$ and $\norm{U_{1,K}(v)-c_0}_{\infty,K}\leq R$. Suppose the earlier stages are in the enlarged interval and satisfy $\norm{U_{j,K}(v)-c_0}_{\infty,K}\leq RP_j(z_{K,\infty})$. Since $\loc_K(c_0)=0$, \Cref{lem:GLip} gives
	\[
		\begin{aligned}
			\norm{U_{\ell,K}(v)-v}_{\infty,K}
			&\leq z_{K,\infty}\sum_{j<\ell}\abs{a_{\ell j}}
			\norm{U_{j,K}(v)-c_0}_{\infty,K}\\
			&\leq R\left(P_\ell(z_{K,\infty})-1\right)\leq\varepsilon_I.
		\end{aligned}
	\]
	The new stage therefore lies in $[\widetilde m,\widetilde M]$. Adding $\norm{v-c_0}_{\infty,K}\leq R$ proves the inductive bound for that stage as well. This argument applies to every comparison input.

	\emph{Flux estimate.} Both stage sequences for $v$ and $w$ lie in the flux domain. Subtracting their equations gives
	\[
		\begin{aligned}
			\norm{U_{\ell,K}(v)-U_{\ell,K}(w)}_{1,K}
			&\leq\norm{v-w}_{1,K}
			+z_{K,1}\sum_{j<\ell}\abs{a_{\ell j}}
			\norm{U_{j,K}(v)-U_{j,K}(w)}_{1,K}\\
			&\leq P_\ell(z_{K,1})\norm{v-w}_{1,K},
		\end{aligned}
	\]
	where the last step is induction from $P_1=1$. Consequently,
	\[
		\norm{\FF_{K,\tau}(v)-\FF_{K,\tau}(w)}_{1,K}
		\leq\alpha\sum_\ell\abs{b_\ell}P_\ell(z_{K,1})\norm{v-w}_{1,K}.
	\]
	This proves \eqref{eq:rkLip}.

	\emph{Completed-step admissibility.} All stages preserve a constant input, and $\sum_\ell b_\ell=1$ gives constant consistency of the volume flux. Equation \eqref{eq:rkLip} and \eqref{eq:volume-lipschitz} verify \Cref{ass:weighted-lipschitz} with $L_{aK}(\tau)=\alpha\norm{\nabla\psi_a}_{\infty,K} B_{\Ark,\brk}(z_{K,1})$. \Cref{thm:main} then gives \Cref{part:rk-admissibility}.
\end{proof}

\begin{remark}[Relation with stability polynomial]
	For entrywise nonnegative tableaux, $B_{\Ark,\brk}(z)=(R_{\mathrm{RK}}(z)-1)/z$, where $R_{\mathrm{RK}}$ is the usual stability polynomial. For signed tableaux the absolute-coefficient recursion \eqref{eq:Pell}, rather than the stability polynomial alone, is the relevant majorant. The patch argument does not require an SSP convex decomposition.
\end{remark}

\subsubsection{ADER-DG with local space-time predictors}\label{subsub:ader-estimates}

The weak equations and the two ADER predictor maps were defined in \Cref{sec:ader}. We first justify the linear time solve, then prove separate estimates for the converged and fixed-Picard predictors.

\begin{lemma}[The local time solve]\label{lem:time-solve}
	Let $\ZK$ be a finite-dimensional local space with the exact $L^2(K)$ inner product. Under \Cref{ass:ader-time}, the operator $\mathcal K_r$ in \eqref{eq:timeoperator} is well defined. For $p=1,\infty$, there are constants $C_{r,1}=C_r$ and $C_{r,\infty}=C_r^{(\infty)}$, depending only on the temporal degree and rule, such that
	\begin{equation}\label{eq:time-solve-bound}
		\norm{\mathcal K_r z}_{\mathbb{X},p,K}
		\leq C_{r,p}\max_i\norm{z(s_i)}_{p,K}.
	\end{equation}
\end{lemma}

\begin{proof}
	Let $y\in \mathbb{X}_K^r$ solve the homogeneous equation. Testing with $y$ gives
	\[
		0=\norm{y(1)}_{2,K}^2-\int_0^1(y,\partial_sy)_K\dd s
		=\frac{1}{2}\left(\norm{y(1)}_{2,K}^2+\norm{y(0)}_{2,K}^2\right).
	\]
	Both endpoint values vanish. Integration by parts then shows that $\int_0^1(\partial_sy,\phi)_K\dd s=0$, $\forall\phi\in \mathbb{X}_K^r$. For $r\geq 1$, choose $\phi=\partial_sy$ to obtain $\partial_sy=0$. For $r=0$, $y$ is already constant. In either case $y=0$, so the square time matrix is nonsingular.

	In a temporal basis, its inverse yields scalar polynomials $\kappa_i(s)$ with $\mathcal K_rz(s)=\sum_i\kappa_i(s)z(s_i)$. The triangle inequality gives \eqref{eq:time-solve-bound} with $\sup_{s\in[0,1]}\sum_i\abs{\kappa_i(s)}$ as a possible constant. This bound does not depend on the spatial mesh and works for any of the spatial norms used below.
\end{proof}

Set
\[
	\kappa_{K,1}=C_r C_{k_K}\frac{\alpha\tau}{h_K},
	\quad
	\kappa_{K,\infty}=C_r^{(\infty)}C_{k_K}^{(\infty)}\frac{\alpha\tau}{h_K},
	\quad z_{K,1}=C_{k_K}\frac{\alpha\tau}{h_K}.
\]

\begin{proposition}[Converged ADER predictor]\label{prop:ader-converged}
	For the meshes and spaces of \Cref{subsec:mesh-space}, suppose \Cref{ass:scalar-data,ass:ader-time} hold. Define the predictor and corrector flux by \eqref{eq:ader-fixed} and \eqref{eq:ader-flux}, and let $v,w\in\ZK$ satisfy $m\leq v,w\leq M$ on $\overline K$.

	\begin{enumerate}
		\renewcommand{\labelenumi}{\normalfont(\roman{enumi})}
		\renewcommand{\theenumi}{(\roman{enumi})}
		\item\label[claim]{part:ader-converged-local} \emph{Local estimates.} Suppose the time step restriction
		\begin{equation}\label{eq:ader-converged-range}
			\kappa_{K,1}<1,
			\quad
			\kappa_{K,\infty}<1,
			\quad
			R\frac{\kappa_{K,\infty}}{1-\kappa_{K,\infty}}\leq\varepsilon_I.
		\end{equation}
		Then \eqref{eq:ader-fixed} has a unique solution in the enlarged flux interval for every comparison input $v\in[m,M]$, and
		\begin{equation}\label{eq:ader-converged-lip}
			\norm{q_K(v)-q_K(w)}_{\mathbb{X},1,K}
			\leq\frac{1}{1-\kappa_{K,1}}\norm{v-w}_{1,K}.
		\end{equation}
		Consequently
		\begin{equation}\label{eq:ader-converged-flux}
			\norm{\FF_{K,\tau}^{\mathrm{ADER}}(v)-\FF_{K,\tau}^{\mathrm{ADER}}(w)}_{1,K}
			\leq\frac{\alpha}{1-\kappa_{K,1}}\norm{v-w}_{1,K}.
		\end{equation}
		\item\label[claim]{part:ader-converged-admissibility} \emph{Completed-step admissibility.} For \eqref{eq:corrector}, also require \Cref{ass:pou} and $m\leq u_h^n\leq M$ cellwise, $\eta_{aK}>0$, and a conservative face flux \eqref{eq:corrector-face}. The raw patch averages then lie in $[m,M]$ if, in addition to \eqref{eq:ader-converged-range} on every cell, the CFL condition holds
		\begin{equation}\label{eq:ader-converged-patch-cfl}
			\tau\norm{\nabla\psi_a}_{\infty,K}
			\frac{\alpha}{1-\kappa_{K,1}}\leq\eta_{aK}
			\quad\text{for every active pair.}
		\end{equation}
	\end{enumerate}
\end{proposition}

\begin{proof}
	\emph{Existence and uniqueness.} Extend $v$ as a constant in normalized time and set $c_0=(m+M)/2$. Put $D=R\kappa_{K,\infty}/(1-\kappa_{K,\infty})$ and consider the closed ball
	\[
		\mathcal B_v=\{q\in \mathbb{X}_K^r:\norm{q-v}_{\mathbb{X},\infty,K}\leq D\}.
	\]
	Because $D\leq\varepsilon_I$, every member takes values in the enlarged flux interval. In particular, the nonlinear map $\Phi_v(q)=v+\tau\mathcal K_r\mathscr R_K[\bfv(q)]$ is defined on this ball. Using $\mathscr R_K[\bfv(c_0)]=0$, \Cref{lem:GLip,lem:time-solve} give
	\[
		\norm{\Phi_v(q)-v}_{\mathbb{X},\infty,K}
		=  \tau\norm{
			\mathcal K_r\mathscr R_K
			[\mathbf f(q)-\mathbf f(c_0)] 
		}_{\mathbb{X},\infty,K}\leq\kappa_{K,\infty}\norm{q-c_0}_{\mathbb{X},\infty,K}
		\leq\kappa_{K,\infty}(D+R)=D.
	\]
	Thus $\Phi_v$ maps $\mathcal B_v$ into itself. For $q,\widetilde q\in\mathcal B_v$, \Cref{lem:GLip} gives
	\[
		\norm{\Phi_v(q)-\Phi_v(\widetilde q)}_{\mathbb{X},\infty,K} 
		= \tau\norm{
			\mathcal K_r\mathscr R_K
			[\mathbf f(q)-\mathbf f(\widetilde q)]
		}_{\mathbb{X},\infty,K} 
		\leq \kappa_{K,\infty} \norm{q-\widetilde q}_{\mathbb{X},\infty,K}.
	\]
	The factor is less than one, so the contraction theorem supplies a fixed point in the ball. Uniqueness also holds among all solutions taking values in the enlarged interval: their difference satisfies the same estimate with factor $\kappa_{K,\infty}<1$.

	\emph{Dependence on the input.} For $v,w$, subtract their fixed-point equations \eqref{eq:ader-fixed} and use the spatial $L^1$ bounds to obtain
	\[
		\norm{q_K(v)-q_K(w)}_{\mathbb{X},1,K}
		\leq\norm{v-w}_{1,K}
		+\kappa_{K,1}\norm{q_K(v)-q_K(w)}_{\mathbb{X},1,K}.
	\]
	Moving the last term to the left proves \eqref{eq:ader-converged-lip}. At the corrector nodes, positivity and normalization of the time weights give
	\[
		\norm{\FF_{K,\tau}^{\mathrm{ADER}}(v)-\FF_{K,\tau}^{\mathrm{ADER}}(w)}_{1,K}
		\leq\alpha\sum_i\vartheta_i\norm{q_K(v)(s_i)-q_K(w)(s_i)}_{1,K},
	\]
	which proves \eqref{eq:ader-converged-flux}.

	\emph{Completed-step admissibility.} A constant is a solution of the predictor equation. Uniqueness and $\sum_i\vartheta_i=1$ give \eqref{eq:constantflux}. The flux bound verifies \Cref{ass:weighted-lipschitz} with $L_{aK}(\tau)=\alpha\norm{\nabla\psi_a}_{\infty,K}/(1-\kappa_{K,1})$. \Cref{thm:main} proves \Cref{part:ader-converged-admissibility}.
\end{proof}

\begin{proposition}[Fixed-Picard ADER predictor]\label{prop:ader-fixed-picard}
	For the meshes and spaces of \Cref{subsec:mesh-space}, suppose \Cref{ass:scalar-data,ass:ader-time} hold. Fix an integer $J\ge0$ independently of the input, use \eqref{eq:ader-picard} and \eqref{eq:ader-flux}, and let $v,w\in\ZK$ satisfy $m\leq v,w\leq M$ on $\overline K$.

	\begin{enumerate}
		\renewcommand{\labelenumi}{\normalfont(\roman{enumi})}
		\renewcommand{\theenumi}{(\roman{enumi})}
		\item\label[claim]{part:ader-picard-local} \emph{Local estimates.} If the time step restriction holds
		\begin{equation}\label{eq:ader-containment}
			R\sum_{j=1}^J\kappa_{K,\infty}^j\leq\varepsilon_I,
		\end{equation}
		then every Picard iterate remains in the enlarged flux interval and
		\begin{equation}\label{eq:ader-growth}
			\norm{q_K^{[J]}(v)-q_K^{[J]}(w)}_{\mathbb{X},1,K}
			\leq\left(\sum_{j=0}^J\kappa_{K,1}^j\right)\norm{v-w}_{1,K}.
		\end{equation}
		The corrector flux consequently satisfies \eqref{eq:volume-lipschitz} with $L_K^{\mathrm{vol}}(\tau)=\alpha B_{r,J}^{\mathrm{ADER}}(z_{K,1})$, where
		\begin{equation}\label{eq:ader-B}
			B_{r,J}^{\mathrm{ADER}}(z)=\sum_{j=0}^J(C_rz)^j.
		\end{equation}
		\item\label[claim]{part:ader-picard-admissibility} \emph{Completed-step admissibility.} For \eqref{eq:corrector}, also require \Cref{ass:pou}, $m\leq u_h^n\leq M$ cellwise, $\eta_{aK}>0$, and a conservative face flux \eqref{eq:corrector-face}. The raw patch averages lie in $[m,M]$ if \eqref{eq:ader-containment} holds on every cell and the CFL condition holds 
		\begin{equation}\label{eq:ader-picard-patch-cfl}
			\tau\alpha\norm{\nabla\psi_a}_{\infty,K}
			B_{r,J}^{\mathrm{ADER}}(z_{K,1})\leq\eta_{aK}
			\quad\text{for every active pair.}
		\end{equation}
		No contraction assumption is required for a prescribed finite number of Picard updates.
	\end{enumerate}
\end{proposition}

\begin{proof}
	\emph{Predictor containment.} Set $c_0=(m+M)/2$. Starting with $q_K^{[0]}=v$, we prove simultaneously that the iterates are in the flux interval and that
	\[
		\norm{q_K^{[j]}-c_0}_{\mathbb{X},\infty,K}
		\leq R\sum_{i=0}^j\kappa_{K,\infty}^i,
		\quad
		\norm{q_K^{[j]}-v}_{\mathbb{X},\infty,K}
		\leq R\sum_{i=1}^j\kappa_{K,\infty}^i.
	\]
	The assertions hold at $j=0$. If they hold through $j$, then the local operator can be evaluated at $q_K^{[j]}$, and
	\[
		\norm{q_K^{[j+1]}-v}_{\mathbb{X},\infty,K}
		\leq\kappa_{K,\infty}\norm{q_K^{[j]}-c_0}_{\mathbb{X},\infty,K}
		\leq R\sum_{i=1}^{j+1}\kappa_{K,\infty}^i\leq\varepsilon_I.
	\]
	This gives containment of the new iterate. Adding $\norm{v-c_0}_\infty\leq R$ closes the induction.

	\emph{Flux estimate.} For two input polynomials define $d_j=\norm{q_K^{[j]}(v)-q_K^{[j]}(w)}_{\mathbb{X},1,K}$. The same operator estimates give $d_0=\norm{v-w}_{1,K}$ and $d_{j+1}\leq\norm{v-w}_{1,K}+\kappa_{K,1}d_j$. Induction yields \eqref{eq:ader-growth}. Applying the physical-flux Lipschitz bound at the time nodes gives
	\[
		\norm{\FF_{K,\tau}^{\mathrm{ADER}}(v)-\FF_{K,\tau}^{\mathrm{ADER}}(w)}_{1,K}
		\leq\alpha\left(\sum_{j=0}^J\kappa_{K,1}^j\right)\norm{v-w}_{1,K}.
	\]
	\emph{Completed-step admissibility.} Each Picard update fixes constants. The flux estimate verifies \Cref{ass:weighted-lipschitz} with $L_{aK}(\tau)=\alpha\norm{\nabla\psi_a}_{\infty,K} B_{r,J}^{\mathrm{ADER}}(z_{K,1})$, so \Cref{thm:main} applies. Neither induction requires $\kappa_{K,1}<1$ or $\kappa_{K,\infty}<1$.
\end{proof}

The converged and fixed-Picard maps therefore differ only in the temporal growth factor seen by the patch theorem: a resolvent for the nonlinear solve and a finite geometric sum for a prescribed number of corrections. The latter is the predictor used in the numerical experiments in \Cref{sec:numerical}. The iteration count for formal accuracy is described in \Cref{rem:ader-iteration-order}; efficient variants that increase the approximation degree during the iterations are developed in \cite{10.1016/j.amc.2023.128426,10.1007/s42967-023-00290-w}.

\subsubsection{LWDG}\label{subsub:lw-estimates}

For the method in \Cref{subsec:lwdg}, the estimate is obtained directly from the CK coefficients of the old polynomial. No predictor range induction is needed, but higher derivatives of the physical flux are required.

\begin{lemma}[Lipschitz estimates for the Taylor coefficients]\label{lem:ck-lipschitz}
	Fix integers $p\geq 1$ and $k_K$ in a bounded range. Let $\ZK$ consist of affine images of fixed reference polynomial spaces on a shape-regular family of affine cells, and suppose $\bfv\in C^p$ on an open neighborhood of $[m,M]$. Define $\mathtt A_j$ by \eqref{eq:lw-jets} using exact spatial derivatives. For $0\leq j\leq p-1$ there are constants $D_{j,K}$ such that
	\begin{equation}\label{eq:lw-coefficient-lip}
		\norm{\mathtt A_j(v)-\mathtt A_j(w)}_{1,K}
		\leq D_{j,K}h_K^{-j}\norm{v-w}_{1,K}
	\end{equation}
	for $v,w\in\ZK$ taking values in $[m,M]$. The constants are independent of $h_K$ on the stated reference-element family.
\end{lemma}

\begin{proof}
	Write the affine map as $T_K(\xi)=x_K + B_K\xi$. A spatial derivative pulls back as $B_K^{-\top}\widehat\nabla$. Induction in the CK recursion shows that each term in $\mathtt A_j$ has total spatial derivative order $j$: taking one more time derivative and substituting the conservation law adds one spatial differentiation. Consequently
	\[
		\widehat{\mathtt A}_{j,K}(\widehat v)
		:=h_K^j\mathtt A_j(v)\circ T_K
	\]
	depends on the geometry only through the bounded matrices $h_K B_K^{-1}$. It involves derivatives of $\widehat v$ of order at most $j$ and physical-flux derivatives of order at most $j$.

	The reference polynomials with values in $[m,M]$ form a compact convex set in coefficient space. Since $\bfv\in C^{j+1}$, the map $\widehat v\mapsto\widehat{\mathtt A}_{j,K}(\widehat v)$ is continuously differentiable on a neighborhood of this set. Its derivative is uniformly bounded there, including over the normalized shape-regular geometries. The fundamental theorem of calculus along $\widehat w+t(\widehat v-\widehat w)$, and equivalence of coefficient and reference $L^1$ norms, therefore give
	\[
		\norm{\widehat{\mathtt A}_{j,K}(\widehat v)
		-\widehat{\mathtt A}_{j,K}(\widehat w)}_{1,\widehat K}
		\leq D_{j,K}\norm{\widehat v-\widehat w}_{1,\widehat K}.
	\]
	Changing variables to $K$ cancels the constant Jacobian on both sides and restores the factor $h_K^{-j}$. This proves \eqref{eq:lw-coefficient-lip}.
\end{proof}

\begin{proposition}[Patchwise weak admissibility for LWDG]\label{prop:lw-scalar}
	Under the affine-cell, polynomial-space, and $C^p$ flux hypotheses of \Cref{lem:ck-lipschitz}, use the Taylor flux \eqref{eq:lw-flux} and its constants $D_{j,K}$. For $v,w\in\ZK$ with $m\leq v,w\leq M$ on $\overline K$, define
	\begin{equation}\label{eq:lw-Lvol}
		L_K^{\mathrm LW}(\tau)
		=\sum_{j=0}^{p-1}\frac{D_{j,K}}{(j+1)!}
		\left(\frac\tau{h_K}\right)^j.
	\end{equation}
	\begin{enumerate}
		\renewcommand{\labelenumi}{\normalfont(\roman{enumi})}
		\renewcommand{\theenumi}{(\roman{enumi})}
		\item\label[claim]{part:lw-local} \emph{Local estimate.} The averaged flux satisfies
		\[
			\norm{\FF_{K,\tau}^{\mathrm LW}(v)-\FF_{K,\tau}^{\mathrm LW}(w)}_{1,K}
			\leq L_K^{\mathrm LW}(\tau)\norm{v-w}_{1,K}.
		\]
		\item\label[claim]{part:lw-admissibility} \emph{Completed-step admissibility.} For \eqref{eq:corrector}, also use the mesh setting of \Cref{subsec:mesh-space} and require \Cref{ass:pou}, $m\leq u_h^n\leq M$ cellwise, $\eta_{aK}>0$, and a conservative face flux \eqref{eq:corrector-face}. The completed LWDG update is patchwise weakly admissible if the CFL condition holds
		\begin{equation}\label{eq:lw-CFL}
			\tau\norm{\nabla\psi_a}_{\infty,K}L_K^{\mathrm LW}(\tau)
			\leq\eta_{aK}
		\end{equation}
		for every active pair. No additional predictor-containment restriction is needed because there are no intermediate stages.
	\end{enumerate}
\end{proposition}

\begin{proof}
	Sum \eqref{eq:lw-coefficient-lip} with the Taylor weights in \eqref{eq:lw-flux}. This verifies the volume-flux Lipschitz hypothesis \eqref{eq:volume-lipschitz}. \Cref{thm:main} then gives the conclusion.
\end{proof}

The Taylor estimate uses higher derivatives of the physical flux but requires no intermediate predictor-range condition. \Cref{thm:main} therefore applies through different local estimates for the three classes of time discretizations. The Euler analysis below concerns only cRKDG and fixed-Picard ADER, for which positive predictor density can be established directly.

\begin{remark}[Temporal order]
	For any fixed temporal order, suppose the local flux maps are defined on the comparison inputs for sufficiently small $\tau$ and their weighted Lipschitz constants $L_{aK}(\tau)$ remain bounded as $\tau\downarrow0$. Since $\eta_{aK}>0$, the condition $\tau L_{aK}(\tau)\leq\eta_{aK}$ then holds for sufficiently small time steps. The temporal order enters through the local flux estimates and the resulting CFL constants. The estimates above verify this structure for the stated cRKDG, ADER, and LWDG constructions.
\end{remark}

\section{Patchwise weak admissibility for compressible Euler equations}\label{sec:euler}

For compressible Euler equations, PACL combines GQL linear admissibility tests with a normalization based on the old cell average. The GQL linear tests transfer the weighted positivity inequality to density and internal energy. The normalization expresses the old admissibility budget in dimensionless variables, allowing the completed temporal increment to be controlled for cRKDG and fixed-Picard ADER.

\subsection{GQL characterization of the admissible set}

For the ideal-gas Euler equations with $\gamma>1$,
\begin{equation}\label{eul:eq:euler-d}
	\partial_t\vect U+\sum_{r=1}^d\partial_{x_r}\vect f_r(\vect U)=0,
	\quad
	\vect U=(\rho,\vect m^\top,E)^\top,
\end{equation}
with
\[
	\vect f_r(\vect U)=
	\begin{pmatrix}
		m_r\\m_r\vect m/\rho+p\vect e_r\\(E+p)m_r/\rho
	\end{pmatrix}
	,
	\quad
	q(\vect U)=E-\frac{\abs{\vect m}^2}{2\rho},
	\quad
	p=(\gamma-1)q.
\]
Here, \(\mathbf{e}_r\) is the vector with the \(r\)-th component \(1\) and all others \(0\). In this section $\abs{\vect U}_\infty$ is the componentwise maximum norm, and $\norm{\vect U}_{\infty,K}=\max_{\overline K}\abs{\vect U}_\infty$. The unsubscripted vector norm  $\abs{\vect \ \cdot\ }$ is Euclidean. The admissible set is the convex cone
\begin{equation}\label{eul:eq:admissible-d}
	\GE=\{\vect U:\rho>0,\ q(\vect U)>0\}.
\end{equation}
For any fixed $\vect w\in\R^d$, define the linear functional
\begin{equation}\label{eul:eq:linear-test}
	\ell_{\vect w}(\vect U)
	=E-\vect m\cdot\vect w+\frac\rho2\abs{\vect w}^2.
\end{equation}
Completing the square gives
\begin{equation}\label{eul:eq:complete-square}
	\ell_{\vect w}(\vect U)
	=q(\vect U)+\frac\rho2\abs{\vect w-\frac{\vect m}{\rho}}^2,
\end{equation}
so the GQL characterization \cite{10.1137/21M1458247} is
\begin{equation}\label{eul:eq:scalar-admissibility}
	\vect U\in\GE
	\quad\Longleftrightarrow\quad
	\rho>0\ \text{and}\ \ell_{\vect w}(\vect U)>0
	\quad\forall\vect w\in\R^d.
\end{equation}
The advantage is that $\ell_{\vect w}$ is linear in the conservative variables. Thus, whenever the old polynomial is pointwise admissible, \eqref{eq:eta} implies the following scalar weighted inequalities for both $\rho_K^n$ and $\ell_{\vect w}(\vect U_K^n)$:
\begin{equation}\label{eul:eq:physical-old-margin}
	\begin{aligned}
		\frac{1}{\abs{K}}\int_K\psi_a\rho_K^n\dd\bx
		&\geq\eta_{aK}\bar\rho_K^n,\\
		\ell_{\vect w}\left(
		\frac{1}{\abs{K}}\int_K\psi_a\vect U_K^n\dd\bx\right)
		&\geq\eta_{aK}\ell_{\vect w}(\bar{\vect U}_K^n).
	\end{aligned}
\end{equation}
These inequalities quantify the admissibility budget available before the local temporal increment is added.

\subsection{The completed patch contribution}

Apply \Cref{lem:cancellation} componentwise. For either cRKDG or fixed-Picard ADER, let $\{\vect Y_{K,\nu}\}_\nu$ denote the local predictor values entering the conservative corrector, with temporal weights $w_\nu$ satisfying $\sum_\nu w_\nu=1$. Define
\[
	\vect F_{K,\tau,r}=\sum_\nu w_\nu\vect f_r(\vect Y_{K,\nu}),
	\quad
	\omega_{\mathrm{abs}}=\sum_\nu\abs{w_\nu}.
\]
The completed patch mass can be written as
\begin{equation}\label{eul:eq:moment-identity}
	\calM_a(\vect U_h^\star)
	=\sum_{K\in\mathcal T_a}\vect C_{aK},
\end{equation}
where
\begin{equation}\label{eul:eq:cell-contribution}
	\vect C_{aK}
	=\int_K\psi_a\vect U_K^n\dd\bx
	+\tau\int_K\sum_{r=1}^d\vect F_{K,\tau,r}\partial_r\psi_a\dd\bx.
\end{equation}
It suffices to prove that every $\vect C_{aK}$ is admissible, because sums and positive multiples of admissible states remain admissible by \eqref{eul:eq:scalar-admissibility}.

\subsection{Old-cell-average normalization}\label{eul:sec:normalization}

The budget in \eqref{eul:eq:physical-old-margin} scales with the old average density and internal energy. Near vacuum, an absolute temporal increment estimate in the original variables is therefore ineffective. We remove these physical scales before estimating the predictor.

\begin{assumption}[Old relative density]\label{eul:ass:relative}
	For some $0<\delta\leq 1$,
	\begin{equation}\label{eul:eq:relative-density}
		\vect U_K^n(\bx)\in\GE,
		\quad
		\rho_K^n(\bx)\geq\delta\bar\rho_K^n,
		\quad \forall\bx\in\overline K.
	\end{equation}
\end{assumption}
This permits arbitrarily small absolute density and pressure. It only prevents density inside a cell from being arbitrarily smaller than its own average.

Suppress the cell subscript temporarily and define
\begin{equation}\label{eul:eq:cell-scales}
	\bar{\vect v}=\bar{\vect m}/\bar\rho,
	\quad
	c=\sqrt{2q(\bar{\vect U})/\bar\rho}.
\end{equation}
Combining a Galilean shift by $\bar{\vect v}$ \cite[Section~3.4]{10.1111/j.1365-2966.2009.15715.x} with cellwise nondimensionalization with $\bar{\rho}$ and $c$ \cite[Section~3.11.1]{I_do_like_CFD}, define
\begin{equation}\label{eul:eq:state-normalization}
	\vect Z=\stateT\vect U=
	\begin{pmatrix}
		\rho/\bar\rho\\[2pt]
		(\vect m-\bar{\vect v}\rho)/(\bar\rho c)\\[2pt]
		(E-\bar{\vect v}\cdot\vect m+\abs{\bar{\vect v}}^2\rho/2)/(\bar\rho c^2).
	\end{pmatrix}
\end{equation}
This \emph{linear} transformation maps the old cell average into the canonical state $(1,\vect0^\top,1/2)^\top$.

\begin{lemma}[Normalization identities]\label{eul:lem:normalization}
	For the ideal-gas flux \eqref{eul:eq:euler-d} with $\gamma>1$, fix an admissible reference state $\bar{\vect U}\in\GE$. Define $\bar{\vect v}$ and $c$ by \eqref{eul:eq:cell-scales} and the invertible linear map $\stateT$ by \eqref{eul:eq:state-normalization}. Then, for every state $\vect U$ with $\rho(\vect U)>0$,
	\begin{equation}\label{eul:eq:normalized-mean}
		\stateT\bar{\vect U}=(1,\vect0^\top,1/2)^\top,
		\quad
		q(\stateT\vect U)=\frac{q(\vect U)}{\bar\rho c^2}.
	\end{equation}
	Hence $\stateT$ and $\stateT^{-1}$ preserve admissibility. For every $\vect Z$ with positive density, the transformed flux satisfies \begin{equation}\label{eul:eq:covariance}
		\vect g_r(\vect Z)
		:=\stateT\vect f_r(\stateT^{-1}\vect Z)
		=\bar v_r\vect Z+c\vect f_r(\vect Z),
		\quad 1\leq r\leq d.
	\end{equation}
	This identity also holds when the pressure is nonpositive.
\end{lemma}

\begin{proof}
	Write $\vect Z=(r,\vect j^\top,H)^\top$. The inverse transformation is
	\[
		\rho=\bar\rho r,\quad
		\vect m=\bar\rho(\bar{\vect v}r+c\vect j),\quad
		E=\bar\rho\left(\frac{1}{2}\abs{\bar{\vect v}}^2r
		+c\bar{\vect v}\cdot\vect j+c^2H\right).
	\]
	Substituting in $E-\abs{\vect m}^2/(2\rho)$ cancels the terms containing $\bar{\vect v}$ and leaves
	\[
		q(\vect U)=\bar\rho c^2
		\left(H-\frac{\abs{\vect j}^2}{2r}\right)
		=\bar\rho c^2q(\vect Z).
	\]
	Since $\rho=\bar\rho r$ and $q(\vect U)=\bar\rho c^2q(\vect Z)$, with $\bar\rho>0$ and $c>0$, positivity of density and internal energy is equivalent in the original and normalized variables.
	Hence $\vect U\in\GE$ if and only if $\vect Z\in\GE$. Applying the transformation to the old cell average gives normalized density one and momentum zero. The definition of $c$ gives normalized energy $H(\bar{\vect U}) = q(\bar{\vect U})/(\bar\rho c^2)=1/2$.

	For the flux identity, write $p_Z=(\gamma-1)q(\vect Z)$. The preceding calculation gives $p(\vect U)=\bar\rho c^2p_Z$. Substitution of the inverse transformation into the three blocks of $\stateT\vect f_r(\vect U)$ yields
	\[
		\begin{pmatrix}
			\bar v_r r+cj_r\\
			\bar v_r\vect j+c(j_r\vect j/r+p_Z\vect e_r)\\
			\bar v_rH+c(H+p_Z)j_r/r
		\end{pmatrix}
		=\bar v_r\vect Z+c\vect f_r(\vect Z).
	\]
	This proves \eqref{eul:eq:covariance} for every state of positive density.
\end{proof}

Set $\vect Z_1=\stateT_K\vect U_K^n$ and
\[
	\vect W_{aK}=\frac{1}{\abs{K}}\int_K\psi_a\vect Z_1\dd\bx.
\]
Because $\stateT_K$ is fixed on the cell, it commutes with averaging. Consequently,
\[
	\bar{\vect Z}_1=(1,\vect0^\top,1/2)^\top,
	\quad
	\frac{1}{\abs{K}}\int_K\ell_{\vect w}(\vect Z_1)\dd\bx
	=\ell_{\vect w}(\bar{\vect Z}_1)=\frac{1+\abs{\vect w}^2}{2}.
\]
For each fixed $\vect w$, both $(\vect Z_1)_\rho$ and $\ell_{\vect w}(\vect Z_1)$ belong to the local approximation space and are positive. Applying \eqref{eq:eta} to these scalar functions gives the explicit old budget
\begin{equation}\label{eul:eq:old-weighted-margin}
	(\vect W_{aK})_\rho\geq\eta_{aK},
	\quad
	\ell_{\vect w}(\vect W_{aK})
	\geq\frac{\eta_{aK}}2(1+\abs{\vect w}^2).
\end{equation}
The lower bounds in \eqref{eul:eq:old-weighted-margin} are independent of the absolute old cell average density and internal energy.

\subsection{Admissibility under normalized temporal increments}

Write the normalized cell contribution as
\begin{equation}\label{eul:eq:old-plus-error}
	\frac{\stateT_K\vect C_{aK}}{\abs{K}}
	=\vect W_{aK}+\vect r_{aK}.
\end{equation}
The old contribution $\vect W_{aK}$ has the lower bounds \eqref{eul:eq:old-weighted-margin}. The predictor determines the temporal increment $\vect r_{aK}$. The following lemma gives a sufficient componentwise bound on this temporal increment for all linear admissibility tests.

Define the dimension-dependent constant
\begin{equation}\label{eul:eq:test-constant}
	C_d=
	\sup_{\vect w\in\R^d}
	\frac{1+\abs{\vect w}_1+\abs{\vect w}^2/2}{1+\abs{\vect w}^2}
	=\frac{3+\sqrt{1+4d}}4.
\end{equation}

\begin{lemma}[Admissibility under normalized temporal increments]\label{eul:lem:budget}
	Fix $\eta_{aK}>0$. Assume that $\vect W_{aK}$ satisfies both lower bounds in \eqref{eul:eq:old-weighted-margin}, and let $\vect r_{aK}$ be any temporal increment vector. With $C_d$ defined by \eqref{eul:eq:test-constant}, if $\varepsilon\ge0$ and $C_d\abs{\vect r_{aK}}_\infty\leq\varepsilon$, then
	\begin{equation}\label{eul:eq:cell-contribution-bound}
		\begin{aligned}
			(\vect W_{aK}+\vect r_{aK})_\rho
			&\geq\eta_{aK}-\varepsilon/C_d,\\
			\ell_{\vect w}(\vect W_{aK}+\vect r_{aK})
			&\geq(\eta_{aK}/2-\varepsilon)(1+\abs{\vect w}^2)
			\quad\forall\vect w\in\R^d.
		\end{aligned}
	\end{equation}
	In particular, $\varepsilon\leq\eta_{aK}/4$ retains a fixed positive fraction of the old admissibility budget.
\end{lemma}

\begin{proof}
	For $t=\abs{\vect w}$, the inequality $\abs{\vect w}_1\leq\sqrt d\,t$ reduces the ratio defining $C_d$ to $(1+\sqrt d\,t+t^2/2)/(1+t^2)$. Its maximum over $t\ge0$ is $(3+\sqrt{1+4d})/4$. Equality in the first inequality is possible when all components of $\vect w$ have equal magnitude.

	The density component changes by at most $\abs{\vect r_{aK}}_\infty\leq\varepsilon/C_d$. For each linear GQL test, the triangle inequality gives
	\[
		\abs{\ell_{\vect w}(\vect r_{aK})}
		\leq\abs{\vect r_{aK}}_\infty
		\left(1+\abs{\vect w}_1+\frac{1}{2}\abs{\vect w}^2\right)
		\leq\varepsilon(1+\abs{\vect w}^2).
	\]
	Subtracting these bounds from \eqref{eul:eq:old-weighted-margin} proves \eqref{eul:eq:cell-contribution-bound}. In particular, if $\varepsilon\leq\eta_{aK}/4$, then every linear GQL test is at least $\eta_{aK}(1+\abs{\vect w}^2)/4$ and the density is at least $\eta_{aK}(1-1/(4C_d))>0$.
\end{proof}
It remains to bound the actual temporal increment $\vect r_{aK}$.

\subsection{Predictor density and correction size}\label{eul:sec:stage}

Choose a nonnegative-polynomial inverse constant $\iota_K\geq 1$ such that
\begin{equation}\label{eul:eq:positive-inverse}
	\norm{p}_{\infty,K}
	\leq\iota_K\frac{1}{\abs{K}}\int_Kp\dd\bx,
	\quad \forall p\in\ZK,
	\quad p\ge0.
\end{equation}
This is the $r=1$ case of \eqref{eq:local-inverse-trace}, restricted to nonnegative functions. Apply it to $(\vect Z_1)_\rho$ and $2(\vect Z_1)_E$, whose means are both one. Admissibility of $q(\vect Z_1)$ also gives $\abs{(\vect Z_1)_m}^2<2(\vect Z_1)_\rho(\vect Z_1)_E$. Together with the relative-density assumption, these facts yield
\begin{equation}\label{eul:eq:old-normalized-bound}
	\delta\leq (\vect Z_1)_\rho\leq\iota_K,
	\quad
	0<2(\vect Z_1)_E\leq\iota_K,
	\quad
	\abs{(\vect Z_1)_m}\leq\iota_K.
\end{equation}
For each $\bx\in\overline K$, consider the normalized states within
half the old normalized density of $\vect Z_1(\bx)$ in the
componentwise maximum norm. Their union is
\begin{equation}\label{eul:eq:normalized-domain}
	\mathcal D_K
	=\bigcup_{\bx\in\overline K}
	\left\{
		\vect Z:
		\abs{\vect Z-\vect Z_1(\bx)}_\infty
		\leq \frac12(\vect Z_1(\bx))_\rho
	\right\}.
\end{equation}
Define the old-state wave speed and local CFL number by
\begin{equation}\label{eul:eq:old-wave-speed}
	\alpha_K^n
	=\max_{\bx\in\overline K}
	\left(
		\abs{\vect v_K^n(\bx)}
		+\sqrt{\frac{\gamma p_K^n(\bx)}{\rho_K^n(\bx)}}
	\right),
	\quad
	\nu_K^n=\frac{\tau\alpha_K^n}{h_K}.
\end{equation}

\begin{lemma}[Compressible Euler flux bound on predictor stages]\label{eul:lem:predictor-flux}
	Suppose $\vect U_K^n(\bx)\in\GE$ for every $\bx\in\overline K$, and let $\iota_K$ satisfy \eqref{eul:eq:positive-inverse}.
	For $\mathcal D_K$ and $\alpha_K^n$ defined in \eqref{eul:eq:normalized-domain} and \eqref{eul:eq:old-wave-speed}, the transformed Euler flux satisfies
	\begin{equation}\label{eul:eq:Phi-bound}
		\max_{\vect Z\in\mathcal D_K}\max_{1\leq s\leq d}
		\abs{\vect g_s(\vect Z)}_\infty
		\leq C_{\gamma,d}\,\iota_K\alpha_K^n, 
	\end{equation}
	where $C_{\gamma,d}$ depends only on $\gamma$ and $d$.
\end{lemma}

\begin{proof}
Fix $\vect Z=(r,\vect j^\top,H)^\top\in\mathcal D_K$ and choose $\bx\in\overline K$ as in \eqref{eul:eq:normalized-domain}.
Write $\vect Z_1(\bx)=(r_1,\vect j_1^\top,H_1)^\top$.
Then $r\geq r_1/2$ and $\abs{\vect j-\vect j_1}\leq\sqrt d\,r_1/2$.
Using the upper bounds in \eqref{eul:eq:old-normalized-bound} and $\abs{\vect j_1}^2<2r_1H_1$, we obtain
\[
	\abs{\vect Z}_\infty\leq\frac32\iota_K,
	\quad
	\abs{H}\leq\iota_K,
	\quad
	\frac{\abs{\vect j}^2}{r}
	\leq 4\frac{\abs{\vect j_1}^2}{r_1}+d r_1
	\leq (4+d)\iota_K .
\]
Hence $\abs{p(\vect Z)}\leq P_{\gamma,d}\iota_K$, where $P_{\gamma,d}=(\gamma-1)(3+d/2)$.

Put $b_\gamma=(1+2/(\gamma(\gamma-1)))^{1/2}$.
The old mean velocity satisfies $\abs{\bar{\vect v}_K}\leq\alpha_K^n$.
The ideal-gas energy identity gives
\[
	2E_K^n
	=\rho_K^n\abs{\vect v_K^n}^2+\frac{2p_K^n}{\gamma-1}
	\leq b_\gamma^2(\alpha_K^n)^2\rho_K^n .
\]
Taking cell averages and using \eqref{eul:eq:cell-scales} yields $c_K^2\leq 2\bar E_K^n/\bar\rho_K^n \leq b_\gamma^2(\alpha_K^n)^2$.
For each direction $s$, the definition of $\mathcal D_K$ gives
\[
	c_K\frac{\abs{j_s}}{r}
	\leq 2c_K\frac{\abs{j_{1,s}}}{r_1}+c_K
	=2\abs{v_{K,s}^n(\bx)-\bar v_{K,s}}+c_K
	\leq (4+b_\gamma)\alpha_K^n .
\]
With $v_s=\bar v_{K,s}+c_Kj_s/r$, \eqref{eul:eq:covariance} gives
\[
	\vect g_s(\vect Z)
	=v_s\vect Z+p(\vect Z)
	\begin{pmatrix}
		0\\ c_K\vect e_s\\ v_s-\bar v_{K,s}
	\end{pmatrix}.
\]
These estimates prove \eqref{eul:eq:Phi-bound} with
\[
	C_{\gamma,d}
	=\frac32(5+b_\gamma)+(4+b_\gamma)P_{\gamma,d}.
\]
\end{proof}

Let the componentwise discrete divergence \eqref{eq:linear-spatial-discrete-div} satisfy
\begin{equation}\label{eul:eq:discrete-div-bound}
	\norm{\mathscr R_K[\mat W]}_{\infty,K}
	\leq\frac{C_{G,K}}{h_K}
	\max_r\norm{\vect w_r}_{\infty,K}.
\end{equation}
Here the discrete divergence \eqref{eq:linear-spatial-discrete-div} acts componentwise on the flux matrix $\mat W=(\vect w_1,\ldots,\vect w_d)$. For the mesh family of \Cref{subsec:mesh-space}, \eqref{eul:eq:discrete-div-bound} follows componentwise from \eqref{eq:discrete-div-infinity-detail} in \Cref{part:GLip-discrete-div} of \Cref{lem:GLip}, with the finite sum over directions absorbed in $C_{G,K}$. 
Define the geometric factor
\begin{equation}\label{eul:eq:beta}
	\beta_{aK}
	=\frac{h_K}{\abs{K}}\int_K\abs{\nabla\psi_a}_1\dd\bx .
\end{equation}
The temporal constants are 
\begin{equation}\label{eul:eq:temporal-constants}
	(C_{\mathrm{pred}},\omega_{\mathrm{abs}})=
	\begin{cases}
		(\max_\ell\sum_{j<\ell}\abs{a_{\ell j}},\ \sum_j\abs{b_j}),&\text{cRKDG},\\
		(C_r^{(\infty)},\ 1),&\text{fixed-Picard ADER}.
	\end{cases}
\end{equation}

\begin{lemma}[Predictor-legality and completed-correction estimates]\label{eul:lem:predictor-correction}
	Assume \Cref{eul:ass:relative} and the bounds \eqref{eul:eq:positive-inverse} and \eqref{eul:eq:discrete-div-bound}.
	Use either the componentwise cRKDG equations \eqref{eq:rkstage}--\eqref{eq:rkflux} with an explicit tableau and $\sum_j b_j=1$, or the fixed-$J$ Picard predictor \eqref{eq:ader-picard}--\eqref{eq:ader-flux} under \Cref{ass:ader-time}. Let $\nu_K^n$ be given by \eqref{eul:eq:old-wave-speed}, and let $C_{\mathrm{pred}}, \omega_{\mathrm{abs}}$ be given by \eqref{eul:eq:temporal-constants}.

	If all preceding predictor values lie in $\mathcal D_K$, the next predictor state satisfies
	\begin{equation}\label{eul:eq:simple-stage-bound}
		\norm{\vect Z-\vect Z_1}_{\infty,K}
		\leq C_{\mathrm{pred}} C_{G,K}
		C_{\gamma,d}\iota_K\nu_K^n ,
	\end{equation}
	uniformly in predictor time for ADER. Thus the condition
	\begin{equation}\label{eul:eq:density-CFL}
		C_{\mathrm{pred}} C_{G,K}
		C_{\gamma,d}\iota_K\nu_K^n
		\leq \delta/2
	\end{equation}
	keeps every required predictor value in $\mathcal D_K$ and gives positive predictor density.

	If all predictor values entering the corrector lie in $\mathcal D_K$, the normalized temporal increment satisfies
	\begin{equation}\label{eul:eq:simple-correction-bound}
		C_d\abs{\vect r_{aK}}_\infty
		\leq C_d\,\omega_{\mathrm{abs}}\,\beta_{aK}
		C_{\gamma,d}\iota_K\nu_K^n .
	\end{equation}
	For each active pair with $\eta_{aK}>0$, the old weighted bound \eqref{eul:eq:old-weighted-margin} and the condition
	\begin{equation}\label{eul:eq:weighted-CFL}
		C_d\,\omega_{\mathrm{abs}}\,\beta_{aK}
		C_{\gamma,d}\iota_K\nu_K^n
		\leq \eta_{aK}/4
	\end{equation}
	ensure that the completed cell contribution $\vect C_{aK}$ belongs to $\GE$.
\end{lemma}

\begin{proof}
	Fix a cell $K$ and write $\mat g(\vect Z)=(\vect g_1(\vect Z),\ldots,\vect g_d(\vect Z))$. 
	The map $\stateT_K$ is linear and fixed during the time step, so it commutes with the componentwise discrete divergence and the local time solve. By \eqref{eul:eq:covariance}, the cRKDG stages satisfy
	\[
		\vect Z_\ell
		=\vect Z_1+\tau\sum_{j<\ell}a_{\ell j}
		\mathscr R_K[\mat g(\vect Z_j)] .
	\]
	If the preceding stage values lie in $\mathcal D_K$, applying \Cref{eul:lem:predictor-flux} to each stage and using \eqref{eul:eq:discrete-div-bound} give
	\[
		\norm{\vect Z_\ell-\vect Z_1}_{\infty,K}
		\leq
		\frac{\tau}{h_K} C_{G,K}
		C_{\gamma,d}\iota_K\alpha_K^n
		\sum_{j<\ell}\abs{a_{\ell j}}
		\leq C_{\mathrm{pred}} C_{G,K}
		C_{\gamma,d}\iota_K\nu_K^n .
	\]
	For fixed-Picard ADER, set $\vect Z^{[j]}=\stateT_K q_K^{[j]}$.
	The normalized iteration is
	\[
		\vect Z^{[j+1]}
		=\vect Z_1+\tau\mathcal K_r
		\mathscr R_K[\mat g(\vect Z^{[j]})] .
	\]
	If the preceding iterate takes values in $\mathcal D_K$, \Cref{lem:time-solve} gives the same estimate uniformly in predictor time, with $C_{\mathrm{pred}}=C_r^{(\infty)}$. 
	This proves \eqref{eul:eq:simple-stage-bound} for both predictors.

	The initial predictor is $\vect Z_1$. 
	Under \eqref{eul:eq:density-CFL}, each new predictor value at $\bx$ differs from $\vect Z_1(\bx)$ by at most $\delta/2\leq(\vect Z_1(\bx))_\rho/2$ in the componentwise maximum norm. 
	It therefore belongs to the neighborhood centered at $\vect Z_1(\bx)$ in \eqref{eul:eq:normalized-domain}, and its density is at least $(\vect Z_1(\bx))_\rho/2>0$. 
	Induction gives the claim for every RK stage or Picard iterate.

	For the completed step, \eqref{eul:eq:cell-contribution} and \eqref{eul:eq:covariance} give
	\[
		\vect r_{aK}
		=\frac{\tau}{\abs{K}}
		\int_K\sum_{r=1}^d\sum_\nu
		w_\nu\vect g_r(\vect Z_{K,\nu})\,
		\partial_r\psi_a\dd\bx .
	\]
	When the values entering this integral lie in $\mathcal D_K$, \eqref{eul:eq:Phi-bound} yields
	\[
		\abs{\vect r_{aK}}_\infty
		\leq
		\frac{\tau}{\abs{K}}\,
		\omega_{\mathrm{abs}} C_{\gamma,d}\iota_K\alpha_K^n
		\int_K\abs{\nabla\psi_a}_1\dd\bx
		=
		\omega_{\mathrm{abs}}\,\beta_{aK}
		C_{\gamma,d}\iota_K\nu_K^n .
	\]
	This proves \eqref{eul:eq:simple-correction-bound}. 
	Under \eqref{eul:eq:weighted-CFL}, \Cref{eul:lem:budget} gives $\vect W_{aK}+\vect r_{aK}\in\GE$. 
	Applying $\stateT_K^{-1}$ and restoring the positive cell-volume factor gives $\vect C_{aK}\in\GE$.
\end{proof}

The two inequalities have distinct roles: \eqref{eul:eq:density-CFL} is a predictor-legality condition, while \eqref{eul:eq:weighted-CFL} is a completed-step admissibility condition.

\begin{theorem}[Euler patchwise weak admissibility]\label{eul:thm:exact}
	For ideal-gas Euler with $\gamma>1$, use the meshes and spaces of \Cref{subsec:mesh-space}. Suppose \Cref{ass:pou,eul:ass:relative} hold, and $\eta_{aK}>0$ for every active pair. Use the componentwise completed update \eqref{eq:corrector} with a face flux satisfying \eqref{eq:corrector-face} and finite at every evaluated predictor trace. Use one of the two predictors specified in \Cref{eul:lem:predictor-correction}, with the inverse bound \eqref{eul:eq:positive-inverse} and the discrete divergence bound \eqref{eul:eq:discrete-div-bound}. If \eqref{eul:eq:density-CFL} holds for every cell and \eqref{eul:eq:weighted-CFL} for every active pair, then all required predictor values have positive density and
	\begin{equation}\label{eul:eq:weak-exact-conclusion}
		\frac{\calM_a(\vect U_h^\star)}{\mathfrak m_a}\in\GE
		\quad\text{for every patch }a.
	\end{equation}
	No admissibility assumption is imposed on predictor pressure or on predictor cell averages.
\end{theorem}

\begin{proof}
	\Cref{eul:lem:predictor-correction} first gives positive density at every required predictor evaluation. The local physical fluxes are therefore defined, and the same lemma bounds each temporal increment. By \Cref{eul:lem:budget},
	\[
		\vect X_{aK}:=\frac{\stateT_K\vect C_{aK}}{\abs{K}}\in\GE.
	\]
	Undoing the normalization gives $\vect C_{aK}=\abs{K}\stateT_K^{-1}\vect X_{aK}\in\GE$. Although different cells use different normalization maps, admissibility has now been established in the same physical variables. Hence
	\[
		\rho\left(\sum_{K\in\mathcal T_a}\vect C_{aK}\right)>0,
		\quad
		\ell_{\vect w}\left(\sum_{K\in\mathcal T_a}\vect C_{aK}\right)
		=\sum_{K\in\mathcal T_a}\ell_{\vect w}(\vect C_{aK})>0
		\quad\forall\vect w\in\R^d.
	\]
	The characterization \eqref{eul:eq:scalar-admissibility} proves admissibility of the sum. The patch identity \eqref{eul:eq:moment-identity} and division by $\mathfrak m_a>0$ give \eqref{eul:eq:weak-exact-conclusion}.
\end{proof}

\begin{remark}[Intermediate pressure]\label{rem:predictor-pressure}
	For cRKDG and fixed-Picard ADER, the theorem permits nonpositive
	pressure during the local evolution, while ensuring positive density.
	The physical Euler flux therefore remains well defined.
	However, numerical interface fluxes may require an extension to
	such states: for example, a wave-speed estimate involving
	$\sqrt{\gamma p/\rho}$ must be modified when $p<0$.
\end{remark}

\begin{remark}[Time-step conditions and forward Euler]\label{rem:euler-forward-euler}
	Both \eqref{eul:eq:density-CFL} and \eqref{eul:eq:weighted-CFL} use the old-state CFL number $\nu_K^n=\tau\alpha_K^n/h_K$ from \eqref{eul:eq:old-wave-speed}. 
	The relative-density parameter $\delta$ enters only the predictor condition \eqref{eul:eq:density-CFL}. We emphasize again that \(\delta\) is only used to control the intermediate predictor stages and will \emph{not} appear if LWDG or predictors that automatically enforce admissibility\footnote{Local predictors are not restricted to use conservative variables; see \Cref{rem:predictor-in-primitive-var}. We only cover the conservative variable case in this paper.} are used. 

	For Forward Euler, $C_{\mathrm{pred}}=0$ and $\omega_{\mathrm{abs}}=1$. The only flux evaluation uses the old state, so $\vect Z_1(\bx)$ belongs to its own neighborhood in \eqref{eul:eq:normalized-domain} at every $\bx\in\overline K$. 
	Thus pointwise admissibility of $\vect U_K^n$ suffices for the flux estimate \eqref{eul:eq:Phi-bound}. 
	Under the mesh and flux hypotheses of \Cref{eul:thm:exact}, the completed patch averages are admissible if
	\[
		\frac{\tau\alpha_K^n}{h_K}\,\beta_{aK}
		\leq \frac{\eta_{aK}}{4C_d C_{\gamma,d}\iota_K}
		\quad\text{for every active pair }(a,K).
	\]
	Thus the forward Euler specialization recovers the characteristic-speed CFL scaling used in the positivity analysis of Zhang and Shu \cite{10.1016/j.jcp.2010.08.016}. Here the sufficient CFL coefficient comes from the weighted patch estimate, and pointwise admissibility of the old state suffices without a prescribed relative-density bound.
\end{remark}

\begin{remark}[Extension to converged ADER and LWDG]
	The Euler analysis also extends naturally to converged ADER and LWDG. For converged ADER, one additionally establishes contraction of the predictor map in the normalized positive-density region. For LWDG, the predictor estimates are replaced by bounds on the normalized Taylor flux coefficients, with no intermediate-stage admissibility requirement. In both cases, the remaining argument uses the same weighted admissibility budget and conservative patch limiting. Since the corresponding proof strategies have already been illustrated in the scalar analysis, we omit these extensions for brevity.
\end{remark}

\section{Conservative patch limiting and accuracy}\label{sec:patch-limiting}

PACL converts the admissible patch averages established in \Cref{sec:framework,sec:euler} into a bound-preserving polynomial in the original DG space. A common contraction within each patch makes its weighted cell states admissible while preserving their assigned mass. Convex recombination then gives admissible physical cell averages, which serve as reference states for the final Zhang--Shu scaling. This construction realizes the second part of \Cref{fig:pacl-vs-zhang-shu}(b). \Cref{thm:pacl-accuracy} establishes approximation-order preservation for a single application of the scalar limiting procedure, under the stated local mesh and scaling assumptions.

\subsection{Patch-limiting algorithm}\label{subsec:limiter-algorithm}

The formulas below use scalar notation. Their vector counterparts
act on all conservative components with the same coefficient.
Let $\mathcal E_K\subseteq\overline K$ be a nonempty compact
admissibility enforcement set for the final Zhang--Shu scaling.
In the quadrature formulation of \Cref{sec:quad}, this set contains all
volume quadrature nodes and all face flux-evaluation nodes.
The choice $\mathcal E_K=\overline K$ enforces whole-cell
pointwise admissibility, while a finite choice enforces
admissibility at the prescribed points.

\begin{algorithm}[Conservative patch limiting]\label{alg:patch-limiting}
	\textit{Input:} raw DG polynomials $u_h^\star$, scalar bounds $[m,M]$ or Euler threshold parameters, admissibility enforcement sets $\mathcal E_K$, and the patch masses $m_{aK},\mathfrak m_a$ from \eqref{eq:patch-masses}. The masses and patch connectivity are precomputed on a fixed mesh.
	\begin{enumerate}
		\item\label[algstep]{step:raw-moments} \emph{Compute the raw moments and check patch admissibility.} For every active pair and every patch, form
		\begin{equation}\label{eq:patch-V-M}
			\avg{u_h^\star}_{aK} = \frac{1}{m_{aK}}\int_K\psi_a u_h^\star \dd\bx,
			\quad
			\avg{u_h^\star}_a =\frac{\calM_a(u_h^\star)}{\mathfrak m_a} = \frac{1}{\mathfrak m_a} \sum_{K\in\mathcal T_a} m_{aK} \avg{u_h^\star}_{aK}.
		\end{equation}
		Recover the raw cell average as $\bar u_h^\star|_K = \abs{K}^{-1}\sum_a m_{aK} \avg{u_h^\star}_{aK}$. Require the weighted patch average $\avg{u_h^\star}_a \in [m,M]$ or $\avg{\vect U_h^\star}_a\in\GE$, respectively. For Euler, choose the positive safeguard thresholds $\epsilon_\rho,\epsilon_q$ sufficiently small so that all admissible patch averages belong to $\GE^\epsilon$. If a weighted patch average fails this GQL test, return failure to the time-step procedure and restart with half the time step. This is possibly due to an insufficient CFL condition.

		\item\label[algstep]{step:patch-coefficient} \emph{Compute one blending coefficient for each patch.} Use \eqref{eq:scalar-patch-theta} for scalar data or \eqref{eq:euler-patch-limiter} for Euler to compute $\theta_a$. The resulting $\theta_a$ is shared by all cells in $\mathcal T_a$.

		\item\label[algstep]{step:patch-recombination} \emph{Correct the weighted states and form the new cell averages.}
		For each active pair $(a,K)$, set
		\begin{equation}\label{eq:patch-limiter}
			\widetilde{\avg{u_h^\star}_{aK}} = \avg{u_h^\star}_a + \theta_a(\avg{u_h^\star}_{aK} - \avg{u_h^\star}_a).
		\end{equation}
		Then, for each cell $K$, collect the corrected weighted mass
		contributions from all patches containing $K$ to obtain its new
		cell average:
		\begin{equation}\label{eq:corrected-cell-average-prop}
			c_K=\frac{1}{\abs{K}}
			\sum_{a:K\in\mathcal T_a}m_{aK}\widetilde{\avg{u_h^\star}_{aK}}.
		\end{equation}
		All patch corrections use the raw moments from
		\Cref{step:raw-moments}, so overlapping patches are processed
		independently rather than by sequentially modifying $u_h^\star|_K$.

		\item\label[algstep]{step:mean-shift} \emph{Shift the original polynomial to the corrected cell average.} Set
		\begin{equation}\label{eq:mean-shift}
			v_K = u_h^\star|_K - \bar u_h^\star|_K + c_K.
		\end{equation}
		For Euler, write the shifted cell average as $\vect C_K$ and the shifted polynomial as $\vect V_K$. An optional oscillation-control procedure may be applied to \(v_K\) or $\vect V_K$ before the final scaling in the next step, provided that it preserves the corrected physical cell average and the local DG space. Examples include suitable TVB \cite{10.1090/S0025-5718-1989-0983311-4,10.48550/arXiv.2606.26692}, WENO \cite{10.1016/j.jcp.2012.08.028,10.1016/j.jcp.2013.04.012,10.1016/j.jcp.2022.111424}, and OEDG \cite{10.1090/mcom/3998} limiters.

		\item\label[algstep]{step:cell-scaling} \emph{Compute the cellwise blending coefficient and apply Zhang--Shu scaling.} Use \eqref{eq:scalar-cell-theta} or \eqref{eq:euler-cell-theta} to compute $\vartheta_K$, respectively, and set
		\begin{align}
			u_h^{n+1}|_{K} := u_h^{\mathrm{lim}}|_K &= c_K + \vartheta_K (v_K-c_K), \label{eq:scalar-final-scaling} \\
			\vect U_h^{n+1}|_K := \vect U_h^{\mathrm{lim}}|_K &= \vect C_K + \vartheta_K (\vect V_K-\vect C_K), \label{eq:euler-final-scaling}
		\end{align}
		for scalar and system cases, respectively. 
		The reference state $c_K$ or $\vect C_K$ is unchanged during oscillation control and scaling.
	\end{enumerate}
\end{algorithm}

\begin{figure}[!ht]
	\centering
	\definecolor{leftblue}{HTML}{175A89}
\definecolor{rightred}{HTML}{B5412C}
\tikzset{
	every node/.style={font=\fontsize{10}{12}\selectfont},
	title/.style={font=\fontsize{10}{12}\selectfont\bfseries},
	curve/.style={line width=1.15pt,line cap=round,line join=round},
	mean/.style={black,line width=1.35pt},
	guide/.style={black!30,densely dotted,line width=.5pt},
	old/.style={black!50,dashed,line width=.7pt},
	flow/.style={-{Latex[length=1.8mm]},black!65,line width=.75pt},
	change/.style={-{Latex[length=1.3mm]},rightred,line width=.75pt},
	patchL/.style={leftblue,dash pattern=on 3pt off 2pt,line width=.85pt},
	patchR/.style={rightred,dash pattern=on 3pt off 2pt,line width=.85pt},
	box/.style={draw=black!50,rounded corners=2pt,inner sep=7pt,fill=white},
}

\newcommand{\Axes}{%
  \foreach \x in {1,2,3}
	\draw[guide] (\x,-.25)--(\x,.85);
  \draw[black!55,line width=.5pt] (0,0)--(3.02,0);
  \node[anchor=east,inner sep=2pt] at (0,0) {$0$};
  \foreach \x/\name in {.5/{K_-},1.5/{K},2.5/{K_+}}
	\node[anchor=north,inner sep=3pt] at (\x,-.265) {$\name$};
}

\newcommand{\Polynomials}[4]{%
  \draw[#4,domain=0:1,samples=70,variable=\s]
	plot ({\s},{.60+.45*\s-.69*\s*\s});
  \draw[#4,domain=0:1,samples=70,variable=\s]
	plot ({1+\s},{-.04+(#1)+(#3)*(1.50*\s*\s-1.98*\s+.49)});
  \draw[#4,domain=0:1,samples=70,variable=\s]
	plot ({2+\s},{.08+.60*\s+.12*\s*\s+(#2)});
}
\newcommand{\Means}[2]{%
  \draw[mean] (.25,.595)--(.75,.595);
  \draw[mean] (1.25,{-.04+(#1)})--(1.75,{-.04+(#1)});
  \draw[mean] (2.25,{.42+(#2)})--(2.75,{.42+(#2)});
}

\begin{tikzpicture}


\node[title] at (2.53,3.85) {(a) Raw DG function};
\begin{scope}[shift={(.50,0)},x=1.35cm,y=4cm]
  \Axes
  \Polynomials{0}{0}{1}{curve,leftblue}
  \Means{0}{0}

  \draw[patchL] (0,.2975)--(2,.2975);
  \draw[patchR] (1,.0900)--(3,.0900);
  \foreach \x in {0,2}
    \draw[leftblue,fill=white,line width=.75pt]
      (\x,.2975) circle[radius=1.4pt];
  \foreach \x in {1,3}
    \draw[rightred,fill=white,line width=.75pt]
      (\x,.0900) circle[radius=1.4pt];
  \node[text=leftblue,anchor=south,inner sep=2pt]
    at (.40,.3025) {$M_L$};
  \node[text=rightred,anchor=south,inner sep=2pt]
    at (2.65,.0950) {$M_R$};
\end{scope}

\begin{scope}[shift={(.50,-2.40)},x=1.35cm,y=.45cm]
  \draw[black!30,line width=.4pt] (0,0)--(3,0);
  \draw[leftblue,line width=.8pt] (0,0)--(1,1)--(2,0);
  \draw[rightred,line width=.8pt] (1,0)--(2,1)--(3,0);
\end{scope}
\node[text=leftblue] at (1.30,-1.72) {Left patch};
\node[text=rightred] at (3.70,-1.72) {Right patch};
\node at (2.53,-2.88) {$M_L=0.2975,\quad M_R=0.09$};

\node[title] at (7.83,3.85) {(b) Corrected cell averages};
\begin{scope}[shift={(5.80,0)},x=1.35cm,y=4cm]
  \Axes
  \Polynomials{0}{0}{1}{old}
  \draw[old] (1.25,-.040)--(1.75,-.040);
  \draw[old] (2.25,.420)--(2.75,.420);
  \Polynomials{.065}{-.065}{1}{curve,leftblue}
  \Means{.065}{-.065}
  \draw[change] (1.82,-.040)--(1.82,.025);
  \draw[change] (2.82,.420)--(2.82,.355);
\end{scope}
\node at (7.83,-1.82) {$c_K=0.025>0$};
\node at (7.83,-2.40) {Pointwise violations remain};

\node[title] at (13.13,3.85) {(c) Cellwise positivity scaling};
\begin{scope}[shift={(11.10,0)},x=1.35cm,y=4cm]
  \Axes
  \draw[leftblue,dashed,line width=.75pt,
        domain=0:1,samples=70,variable=\s]
    plot ({1+\s},{1.50*\s*\s-1.98*\s+.515});
  \Polynomials{.065}{-.065}{75/817}{curve,rightred}
  \Means{.065}{-.065}
  \draw[change,shorten <=1.5pt,shorten >=1.5pt]
    (1.08,.3662)--(1.08,{.025+(75/817)*.3412});
  \draw[change,shorten <=1pt,shorten >=1pt]
    (1.66,-.1384)--(1.66,.0100);
\end{scope}
\node at (13.13,-1.82) {$\min u_K^{n+1}=\epsilon=0.01$};
\node at (13.13,-2.40) {$\bar u_K^{n+1}=c_K=0.025$};

\draw[flow] (4.78,1.45)--(5.50,1.45);
\draw[flow] (10.08,1.45)--(10.80,1.45);

\draw[black!35,line width=.5pt] (.10,-3.60)--(5.30,-3.60);
\draw[black!35,line width=.5pt] (10.35,-3.60)--(15.50,-3.60);
\node at (7.83,-3.60) {Patch correction in (b)};


\node[title] at (2.53,-4.40) {(b.1) Split the cell mass};
\node[text=leftblue] at (1.15,-4.96) {$\psi_L=1-s$};
\node[text=rightred] at (3.72,-4.96) {$\psi_R=s$};
\begin{scope}[shift={(.50,-5.75)},x=4.05cm,y=.55cm]
  \draw[black!35,line width=.4pt] (0,0)--(1,0);
  \draw[leftblue,line width=1pt] (0,1)--(1,0);
  \draw[rightred,dashed,line width=1pt] (0,0)--(1,1);
\end{scope}
\node at (2.53,-6.05) {$s\in[0,1]\quad\text{on }K$};

\node[text=leftblue] at (1.30,-6.52) {$\psi_LW_K$};
\node[text=rightred] at (3.85,-6.52) {$\psi_RW_K$};

\begin{scope}[shift={(.40,-8.16)},x=1.80cm,y=3.20cm]
  \path[fill=leftblue!12]
    (0,0)
    plot[domain=0:1,samples=90,variable=\s]
      ({\s},{(1-\s)*(.45-1.98*\s+1.50*\s*\s)})
    --(1,0)--cycle;
  \draw[black!50,line width=.5pt] (0,0)--(1,0);
  \node[anchor=east,inner sep=2pt] at (0,0) {$0$};
  \draw[leftblue,curve,domain=0:1,samples=90,variable=\s]
    plot ({\s},{(1-\s)*(.45-1.98*\s+1.50*\s*\s)});
  \node[anchor=north,inner sep=3pt] at (.5,-.19) {on $K$};
\end{scope}

\begin{scope}[shift={(2.95,-8.16)},x=1.80cm,y=3.20cm]
  \path[fill=rightred!12]
    (0,0)
    plot[domain=0:1,samples=90,variable=\s]
      ({\s},{\s*(.45-1.98*\s+1.50*\s*\s)})
    --(1,0)--cycle;
  \draw[black!50,line width=.5pt] (0,0)--(1,0);
  \node[anchor=east,inner sep=2pt] at (0,0) {$0$};
  \draw[rightred,curve,domain=0:1,samples=90,variable=\s]
    plot ({\s},{\s*(.45-1.98*\s+1.50*\s*\s)});
  \node[anchor=north,inner sep=3pt] at (.5,-.19) {on $K$};
\end{scope}
\node[text=leftblue] at (2.53,-9.66)
  {$\int_K\psi_LW_K\,dx=0.020$};
\node[text=rightred] at (2.53,-10.21)
  {$\int_K\psi_RW_K\,dx=-0.060$};
\node at (2.53,-10.85) {$0.020-0.060=-0.040$};

\node[title] at (7.83,-4.40) {(b.2) Scale weighted states};
\node at (7.83,-5.02)
  {$V_{aJ}=\frac{1}{m_{aJ}}\int_J\psi_aW_J\,dx$};
\node at (7.83,-5.53) {$m_{aJ}=\tfrac12$};

\node[text=leftblue] at (7.83,-6.12)
  {$M_L=0.2975,\quad\theta_L=1$};
\draw[leftblue,line width=.6pt] (5.90,-6.48) rectangle (9.76,-7.32);
\draw[leftblue!60,line width=.5pt] (7.83,-6.48)--(7.83,-7.32);
\node at (6.87,-6.70) {$K_-$};
\node at (8.80,-6.70) {$K$};
\node[text=leftblue] at (6.87,-7.10) {$0.555$};
\node[text=leftblue] at (8.80,-7.10) {$0.040$};
\node at (7.83,-7.73) {$\tfrac12(0.555+0.040)=M_L$};

\node[text=rightred] at (7.83,-8.30)
  {$M_R=0.09,\quad\theta_R=8/21$};
\begin{scope}[shift={(6.10,-10.30)},x=1.15cm,y=4cm]
  \draw[black!35,line width=.5pt] (0,0)--(3,0);
  \draw[guide] (0,.09)--(3,.09);
  \node[anchor=south,inner sep=2pt] at (1.50,.095) {$M_R$};

  \foreach \x/\before/\after in {.65/-.12/.01,2.35/.30/.17}{
    \draw[old] ({\x-.25},\before)--({\x+.25},\before);
    \draw[rightred,line width=1.15pt]
      ({\x-.25},\after)--({\x+.25},\after);
    \draw[change,shorten <=1pt,shorten >=1pt]
      (\x,\before)--(\x,\after);
  }
  \node[anchor=east,text=black!60] at (.30,-.12) {$-0.120$};
  \node[anchor=east,text=rightred] at (.30,.01) {$0.010$};
  \node[anchor=west,text=black!60] at (2.70,.30) {$0.300$};
  \node[anchor=west,text=rightred] at (2.70,.17) {$0.170$};
  \node[anchor=north,inner sep=3pt] at (.65,-.20) {$K$};
  \node[anchor=north,inner sep=3pt] at (2.35,-.20) {$K_+$};
\end{scope}

\node[title] at (13.13,-4.40) {(b.3) Reweight and collect};
\node[text=leftblue] at (12.00,-5.12) {Left patch};
\node[text=rightred] at (14.25,-5.12) {Right patch};
\node[text=leftblue] at (12.00,-5.78) {$\widetilde V_{LK}$};
\node[text=rightred] at (14.25,-5.78) {$\widetilde V_{RK}$};
\node[text=leftblue] at (12.00,-6.28) {$0.040$};
\node[text=rightred] at (14.25,-6.28) {$0.010$};

\draw[flow,leftblue] (12.00,-6.65)--(12.00,-7.60);
\draw[flow,rightred] (14.25,-6.65)--(14.25,-7.60);
\node[anchor=west] at (12.10,-7.12) {$\times\,\tfrac12$};
\node[anchor=west] at (14.35,-7.12) {$\times\,\tfrac12$};
\node[text=leftblue] at (12.00,-7.96) {$0.020$};
\node[text=rightred] at (14.25,-7.96) {$0.005$};

\node[box,align=center] (newmean) at (13.13,-9.50)
  {$c_K=0.020+0.005$\\[5pt]$=0.025$};
\draw[flow,leftblue] (12.00,-8.30)--($(newmean.north)+(-.55,0)$);
\draw[flow,rightred] (14.25,-8.30)--($(newmean.north)+(.55,0)$);
\node at (13.13,-10.52) {$|K|=1$};

\draw[black!35,line width=.5pt] (.10,-11.75)--(15.50,-11.75);
\node at (7.83,-12.10)
  {Unit cells; positivity threshold $\epsilon=0.01$;
   total mass $0.975$ is preserved.};

\end{tikzpicture}
	\caption{Graphical illustration of the conservative patch limiting procedure: an example.}
	\label{fig:CPL.jpg}
\end{figure}
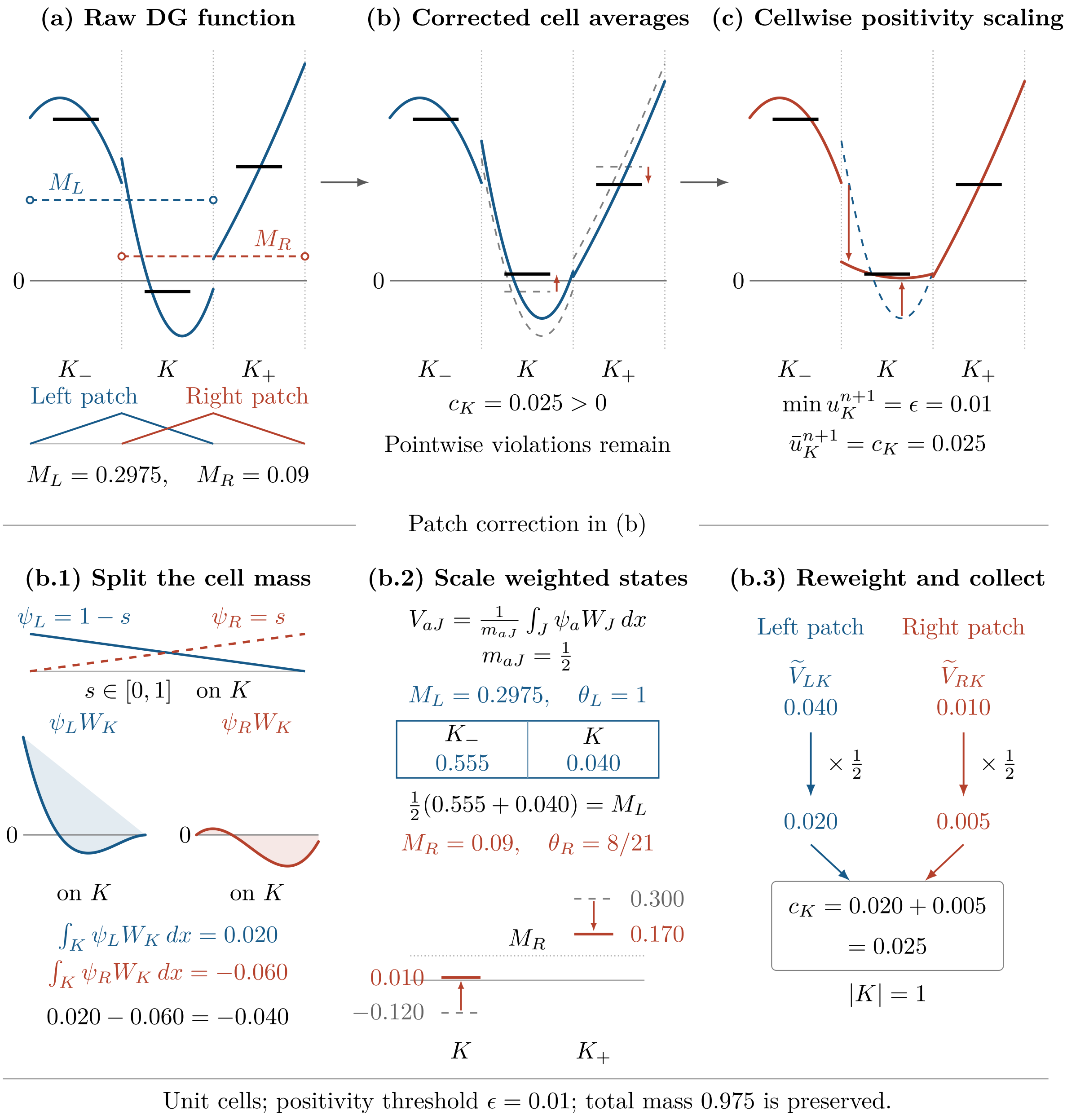

Here $u_h^{\mathrm{lim}}\in\Vh$ denotes the limited scalar output for the raw input $u_h^\star$ defined in \eqref{eq:scalar-final-scaling}. This notation describes one application of the limiter independently of the time-step index. \Cref{alg:patch-limiting} is illustrated by \Cref{fig:CPL.jpg}. 

With the compatible quadrature of \Cref{sec:quad}, all moments in \Cref{step:raw-moments} are evaluated using the same positive volume rule, including the physical Jacobian on mapped cells. Exact state--hat pairing makes these the physical moments used in the analysis. Only the coefficient formulas differ between the scalar and Euler cases.

\begin{remark}[Local nature of the limiting]\label{rem:locality}
	The patch correction in a cell $K$ depends only on the raw data in the patches containing $K$. These patches are processed independently, so the correction does not propagate beyond their union. For the vertex patches used here, this union consists only of cells sharing a vertex with $K$. The subsequent cell average shift and Zhang--Shu scaling are cell-local. 
\end{remark}

\begin{remark}[Cells of degree zero]\label{rem:degree-zero-cells}
On cells with $k_K\geq1$, use the weighted states in \eqref{eq:patch-V-M}. On a cell with $k_K=0$, by direct calculation, we have the auxiliary vertex states 
\[
  m_{aK}\avg{u_h^\star}_{aK}
  = m_{aK} \bar u_K^n
	+ \tau\int_K \bfv(\bar u_K^n)\cdot\nabla\psi_a\,\dd\bx
	- \tau\sum_{e\subset\partial K}
	  \int_e \psi_a\,\widehat{\mathscr F}_{K,e}\,\dd s .
\]
This construction is equivalent to embedding $\bar u_K^n$ as a constant in the first-degree space containing the vertex hats, applying the DG update with volume flux $\bfv(\bar u_K^n)$ and the original face fluxes, and taking the vertex moments of the resulting polynomial. 

The patch averages $\avg{u_h^\star}_a$ are assembled from these auxiliary moments on degree-zero cells and from the raw DG moments on cells of positive degree. Whenever these averages are admissible, \Cref{alg:patch-limiting} applies the same conservative correction on the entire patch. On a degree-zero cell, the corrected average is the final solution. In the scalar case, under the input bounds and boundary assumptions of \Cref{thm:main}, a sufficient CFL condition is \eqref{eq:scalar-patch-cfl} on cells of positive degree and
\[
  \tau L\left|\int_K\nabla\psi_a\,\dd\bx\right|\leq m_{aK}
\]
on cells of degree zero, where $L$ is a Lipschitz constant of $\bfv$ on $[m,M]$.
\end{remark}

\subsection{Limiting coefficients}\label{subsec:limiter-coefficients}

\subsubsection{Scalar equations}

For an admissible reference state $c\in[m,M]$ and a value $z\in\R$, define
\begin{equation}\label{eq:scalar-theta}
	\Theta(c,z)=
	\begin{cases}
		(M-c)/(z-c),&z>M,\\
		(c-m)/(c-z),&z<m,\\
		1,&m\leq z\leq M.
	\end{cases}
\end{equation}
This is the largest $t\in[0,1]$ for which $c+t(z-c)\in[m,M]$. The denominators in the first two cases are strictly positive. The patch and cellwise coefficients in \Cref{alg:patch-limiting} are, respectively,
\begin{align}
	\theta_a&=\min_{K\in\mathcal T_a}\Theta(\avg{u_h^\star}_a,\avg{u_h^\star}_{aK}),
	\label{eq:scalar-patch-theta}\\
	\vartheta_K&=\inf_{\bx\in \mathcal E_K}\Theta(c_K,v_K(\bx)).
	\label{eq:scalar-cell-theta}
\end{align}
For a finite admissibility enforcement set, the infimum is a minimum over its nodes.

\subsubsection{Compressible Euler equations}

Choose artificial safeguards $0<\epsilon_\rho,\epsilon_q\ll1$
as in \Cref{alg:patch-limiting} for the sake of implementation.
Let $\vect C\in\GE$ be the patch reference state satisfying \(\rho(\mathbf C)\ge\epsilon_\rho\) and \(q(\mathbf C)\ge\epsilon_q\), and let $\mathcal S$ be the
bounded, nonempty set of states in conservative variables to be limited. 

As in the scalar case, for a lower bound $\epsilon$ and a reference
value $c\geq\epsilon$, define
\[
	\Theta_\epsilon(c,z)=
	\begin{cases}
		1, & z\geq\epsilon,\\[1mm]
		(c-\epsilon)/(c-z), & z<\epsilon.
	\end{cases}
\]
This is the largest $t\in[0,1]$ for which $c+t(z-c)\geq\epsilon$.
We apply this formula twice, following the concavity argument
in \cite{10.1016/j.jcp.2011.10.002}.

\emph{First, enforce positive density.}
Compute the smallest density and scale all states toward $\vect C$:
\begin{equation}\label{eq:euler-density-cutoff}
	\rho_{\min}=\inf_{\vect U\in\mathcal S}\rho(\vect U),
	\quad
	t_\rho=\Theta_{\epsilon_\rho}\left(\rho(\vect C),\rho_{\min}\right),
	\quad
	\vect U^\rho=\vect C+t_\rho(\vect U-\vect C).
\end{equation}
Since density is linear, every intermediate state satisfies
$\rho(\vect U^\rho)\geq\epsilon_\rho$, so $q(\vect U^\rho)$ is well defined.

\emph{Next, enforce positive internal energy.}
Using the intermediate states, compute
\begin{equation}\label{eq:euler-energy-cutoff}
	q_{\min}=\inf_{\vect U\in\mathcal S}q(\vect U^\rho),
	\quad
	t_q=\Theta_{\epsilon_q}\left(q(\vect C),q_{\min}\right).
\end{equation}
The final state $\widetilde{\vect U}=(1-t_q)\vect C+t_q\vect U^\rho$ retains the density lower bound $\epsilon_\rho$. By concavity of $q$,
\[
	q(\widetilde{\vect U})
	\geq (1-t_q)q(\vect C)+t_q q(\vect U^\rho)
	\geq (1-t_q)q(\vect C)+t_q q_{\min}
	\geq \epsilon_q.
\]
The two scalings therefore give a single coefficient:
\begin{equation}\label{eq:euler-line-coefficient}
	\Theta_{\mathrm E}(\vect C,\mathcal S;\epsilon_\rho,\epsilon_q)
	=t_\rho t_q,
	\quad
	\widetilde{\vect U}=\vect C+t_\rho t_q(\vect U-\vect C).
\end{equation}

The patch and cellwise coefficients in \Cref{alg:patch-limiting}
are, respectively,
\begin{align}
	\theta_a
	&=\Theta_{\mathrm E}\left(
		\avg{\vect U_h^\star}_a,\,
		\{\avg{\vect U_h^\star}_{aK}:K\in\mathcal T_a\};
		\epsilon_\rho,\epsilon_q\right),
	\label{eq:euler-patch-limiter}\\
	\vartheta_K
	&=\Theta_{\mathrm E}\left(
		\vect C_K,\,
		\{\vect V_K(\bx):\bx\in\mathcal E_K\};
		\epsilon_\rho,\epsilon_q\right).
	\label{eq:euler-cell-theta}
\end{align}
For a finite set of states, both infima are minima. For an ideal polytropic gas, a pressure bound $\epsilon_p$ corresponds to
$\epsilon_q = \epsilon_p/(\gamma-1)$.

Consequently, every corrected patch state belongs to
\begin{equation}\label{eq:euler-patch-target}
	\GE^\epsilon
	=\{\vect U:\rho(\vect U)\geq\epsilon_\rho,\ 
				 q(\vect U)\geq\epsilon_q\}.
\end{equation}

\subsection{Conservation and bound preservation}

\begin{proposition}[Conservation and bound preservation]\label{prop:patch-correction}
	Let the weights satisfy \Cref{ass:pou}. Use exact moments in \eqref{eq:patch-V-M} and suppose that every patch reference state is in $[m,M]$ for scalar data or in $\GE^\epsilon$ for Euler. 

	\begin{enumerate}
		\renewcommand{\labelenumi}{\normalfont(\roman{enumi})}
		\renewcommand{\theenumi}{(\roman{enumi})}
		\item\label[claim]{part:patch-conservation} With the coefficients \eqref{eq:scalar-patch-theta} or \eqref{eq:euler-patch-limiter}, \Cref{step:patch-coefficient,step:patch-recombination} of \Cref{alg:patch-limiting} give admissible weighted states and cell averages. Each patch preserves its assigned mass:
		\begin{equation}\label{eq:patch-mass-preservation}
			\sum_{K\in\mathcal T_a}m_{aK}\widetilde{\avg{u_h^\star}_{aK}}
			=\mathfrak m_a \avg{u_h^\star}_a
			=\sum_{K\in\mathcal T_a} m_{aK} \avg{u_h^\star}_{aK}.
		\end{equation}
		The identity holds componentwise for Euler, and recombination preserves every conservative component globally. 

		\item\label[claim]{part:patch-scaling} Use the local spaces of \Cref{subsec:mesh-space} and a nonempty admissibility enforcement set $\mathcal E_K\subseteq\overline K$. Choose positive Euler cell thresholds (typically of machine precision size) that are not greater than the corresponding density and internal energy of $\vect C_K$. \Cref{step:mean-shift,step:cell-scaling} give a polynomial in the original DG space satisfying the prescribed bounds on $\mathcal E_K$ and preserving the shifted cell average. This remains true with an additional limiter in \Cref{step:mean-shift} that preserves the cell average and the local space.
	\end{enumerate}
\end{proposition}

\begin{proof}
	Each line coefficient is feasible, and the minimum over a patch keeps all corrected values in the scalar interval or in $\GE^\epsilon$. Because the same $\theta_a$ is used throughout the patch,
	\[
		\sum_Km_{aK}\widetilde{\avg{u_h^\star}_{aK}}
		=\mathfrak m_a \avg{u_h^\star}_a+\theta_a\sum_K m_{aK} (\avg{u_h^\star}_{aK} - \avg{u_h^\star}_a)
		=\mathfrak m_a \avg{u_h^\star}_a.
	\]
	Recombination weights $m_{aK}/\abs{K}$ are nonnegative and sum to one, so each cell average is admissible by convexity. Summing the patch identities and using the partition of unity proves global conservation.

	For \Cref{part:patch-scaling}, the cell average shift only changes the constant mode. The final scaling preserves the resulting cell average and enforces bounds on admissibility enforcement sets by \eqref{eq:scalar-cell-theta} or \eqref{eq:euler-cell-theta}, as in \cite{10.1016/j.jcp.2009.12.030,10.1016/j.jcp.2010.08.016}. An intervening cell average-preserving limiter does not change this conclusion. The choice \(\mathcal E_K=\overline K\) guarantees admissibility throughout the cell, while a discrete \(\mathcal E_K\) guarantees admissibility only at the prescribed nodes in \(\mathcal E_K\).
\end{proof}

\begin{remark}[Local conservation]\label{rem:local-conservation}
	Let
	\[
		\mathcal A_{\rm act}=\{a:\theta_a<1\},
		\quad
		\mathcal T_{\rm act}
		=\bigcup_{a\in\mathcal A_{\rm act}}\mathcal T_a.
	\]
	Writing
	\[
		\delta_{aK}
		=m_{aK}\left(\widetilde{\avg{u_h^\star}_{aK}} - \avg{u_h^\star}_{aK}\right),
	\]
	the patch conservation identity gives
	\[
		\sum_{K\in\mathcal T_a}\delta_{aK}=0,
		\quad
		\abs{K}\left(c_K - \bar u_h^\star|_K\right)
		=\sum_{a:K\in\mathcal T_a}\delta_{aK}.
	\]
	Consequently,
	\[
		\sum_{K\in\mathcal T_{\rm act}}\abs{K}c_K
		=\sum_{K\in\mathcal T_{\rm act}}\abs{K}\bar u_h^\star|_K,
		\quad
		c_K=\bar u_h^\star|_K
		\quad\text{for }K\notin\mathcal T_{\rm act}.
	\]
	These identities hold componentwise for Euler and remain valid after the final mean-preserving scaling. Thus the redistribution preserves every conserved quantity over the union of the active patches. In the interpretation of \Cref{rem:dual-cell}, \eqref{eq:patch-mass-preservation} states that the patch correction preserves the total conserved quantity assigned to each vertex.
\end{remark}

\subsection{Accuracy of the scalar patch limiter}\label{subsec:patch-accuracy}

The following estimate controls the change in the physical cell average by the error of the raw approximation. The proof follows the correction-estimate strategy used for Zhang--Shu scaling~\cite[Theorem~4 and Appendix~C]{10.1016/j.jcp.2016.10.002}, applied here to conservative transfers between cells. Combined with the final mean-preserving scaling, \Cref{thm:patch-accuracy} yields the full scalar limiter estimate in \Cref{thm:pacl-accuracy}.

Let $u$ satisfy $m\leq u\leq M$, and let $u_h^\star$ be the raw scalar approximation. Define the patch neighborhood by
\begin{equation}
	\omega_K = \bigcup_{a:K\in\mathcal T_a}\bigcup_{J\in\mathcal T_a}J,
	\label{eq:patch-neighborhood}
\end{equation}

\begin{lemma}[Accuracy of the corrected cell average]\label{thm:patch-accuracy}
	Under \Cref{ass:pou}, fix $m<M$. Let $u\in L^\infty(\Omega)$ satisfy $m\leq u\leq M$ and let $u_h^\star\in\Vh$ have $\avg{u_h^\star}_a\in[m,M]$ on every patch. Use exact weighted moments, the maximal common scalar coefficient \eqref{eq:scalar-patch-theta}, and the recombination \eqref{eq:corrected-cell-average-prop}. Let $\omega_K$ be defined by \eqref{eq:patch-neighborhood}. Then
	\begin{equation}\label{eq:meanerror}
		\begin{aligned}
			\abs{c_K-\bar u_h^\star|_K}
			&\leq \frac{\abs{\omega_K}}{\abs{K}}
			\norm{u_h^\star-u}_{L^\infty(\omega_K)},\\
			\abs{c_K-\bar u_K}
			&\leq \left(1+\frac{\abs{\omega_K}}{\abs{K}}\right)
			\norm{u_h^\star-u}_{L^\infty(\omega_K)}.
		\end{aligned}
	\end{equation}
\end{lemma}

\begin{proof}
	Let $\varepsilon_K=\norm{u_h^\star-u}_{L^\infty(\omega_K)}$.
	Fix a patch $a$ incident to $K$. Positivity of the weighted averages and $m\leq u\leq M$ imply
	\[
		m-\varepsilon_K\leq \avg{u_h^\star}_{a J}\leq M+\varepsilon_K
		\quad \text{for } J\in\mathcal T_a.
	\]
	Suppose the lower bound determines the patch coefficient. Let 
	\[
		s=\min_{J\in\mathcal T_a}\avg{u_h^\star}_{a J}<m,
		\quad
		\alpha_{a K}=\frac{m_{aK}}{\mathfrak m_a}.
	\]
	By maximality of $\theta_a$,
	\[
		1-\theta_a=\frac{m-s}{\avg{u_h^\star}_a-s},
		\quad
		m-s\leq\varepsilon_K.
	\]
	Moreover,
	\[
		\avg{u_h^\star}_a-s
		=\sum_{J\in\mathcal T_a}
		\frac{m_{aJ}}{\mathfrak m_a}(\avg{u_h^\star}_{a J}-s),
	\]
	and every term on the right is nonnegative. Hence
	\[
		0\leq \avg{u_h^\star}_{aK}-s\leq\frac{\avg{u_h^\star}_a-s}{\alpha_{a K}}.
	\]
	Subtracting $\avg{u_h^\star}_a-s$ from this interval and multiplying by $\alpha_{a K}$ give
	\[
		\begin{aligned}
			-(\avg{u_h^\star}_a-s)&\leq \avg{u_h^\star}_{aK}-\avg{u_h^\star}_a\leq(\alpha_{a K}^{-1}-1)(\avg{u_h^\star}_a-s),\\
			\alpha_{a K}\abs{\avg{u_h^\star}_{aK}-\avg{u_h^\star}_a}&\leq\max\{\alpha_{a K},1-\alpha_{a K}\}(\avg{u_h^\star}_a-s).
		\end{aligned}
	\]
	Multiplying by $\mathfrak m_a(1-\theta_a)$ yields
	\begin{equation}\label{eq:transfer-bound}
		m_{aK}\abs{\widetilde{\avg{u_h^\star}_{aK}}-\avg{u_h^\star}_{aK}}
		\leq\max\{m_{aK},\mathfrak m_a-m_{aK}\}\varepsilon_K.
	\end{equation}
	If an upper constraint determines the coefficient, apply the same calculation to $-u_h^\star$ and the interval $[-M,-m]$. When $\theta_a=1$, the patch transfer is zero. These cases exhaust the maximal common coefficient, including a tie between active constraints. Finally,
	\[
		\bar u_h^\star|_K=\frac{1}{\abs{K}}\sum_a m_{aK}\avg{u_h^\star}_{aK}.
	\]
	Summing \eqref{eq:transfer-bound} and using the recombination gives
	\[
		\abs{c_K-\bar u_h^\star|_K}
		\leq \frac{\varepsilon_K}{\abs{K}}
		\sum_{a:K\in\mathcal T_a}
		\max\{m_{aK},\mathfrak m_a-m_{aK}\}
		\leq \frac{\abs{\omega_K}}{\abs{K}}\varepsilon_K.
	\]
	Here, nonnegativity and the partition of unity imply
	\[
		\sum_{a:K\in\mathcal T_a}\mathfrak m_a
		=\int_{\omega_K}\sum_{a:K\in\mathcal T_a}\psi_a\dd\bx
		\leq\abs{\omega_K}.
	\]
	The second inequality in \eqref{eq:meanerror} follows from
	$\abs{\bar u_h^\star|_K-\bar u_K}\leq\varepsilon_K$.
\end{proof}

For \Cref{alg:patch-limiting} without the optional limiter, \eqref{eq:mean-shift} and \Cref{thm:patch-accuracy} give
\begin{equation}\label{eq:shiftaccuracy}
	\norm{v_K-u}_{L^\infty(K)}
	\leq \left(1+\frac{\abs{\omega_K}}{\abs{K}}\right)
	\norm{u_h^\star-u}_{L^\infty(\omega_K)}.
\end{equation}
The shifted polynomial $v_K$ has the admissible physical cell average $c_K$, so the final Zhang--Shu scaling preserves its approximation accuracy.

\begin{theorem}[Accuracy preservation of the scalar PACL procedure]\label{thm:pacl-accuracy}
	Under the assumptions of \Cref{thm:patch-accuracy}, let $u_h^{\mathrm{lim}}$ be defined by \eqref{eq:scalar-final-scaling}, with the optional oscillation-control limiter disabled. Assume that the polynomial degrees and the cellwise Zhang--Shu scaling constants are uniformly bounded on the mesh family. Then there is a constant $C$, independent of $h$ and of the distance of $u$ or $c_K$ from the bounds, such that
	\begin{equation}\label{eq:pacl-accuracy}
		\norm{u_h^{\mathrm{lim}}-u}_{L^\infty(K)}
		\leq C\left(1+\frac{\abs{\omega_K}}{\abs{K}}\right)
		\norm{u_h^\star-u}_{L^\infty(\omega_K)}.
	\end{equation}
	If $\abs{\omega_K}/\abs{K}$ is uniformly bounded, one application of the scalar PACL procedure preserves the approximation order of $u_h^\star$, including at active bounds and flat extrema.
\end{theorem}

\begin{proof}
	Apply the Zhang--Shu scaling estimate \cite[Theorem~4, Remark~4, and Appendix~C]{10.1016/j.jcp.2016.10.002} to $v_K$ with its physical cell average $c_K$, and combine it with \eqref{eq:shiftaccuracy} and the triangle inequality.
\end{proof}

\begin{remark}[Local mesh condition]\label{rem:pacl-local-mesh}
	A uniform bound on $\abs{\omega_K}/\abs{K}$ controls only local patch neighborhoods and permits global mesh grading without global quasi-uniformity. In one dimension, a uniform bound on adjacent cell-length ratios provides this control.
\end{remark}

\section{Quadrature version}\label{sec:quad}

Compatible quadrature gives a discrete counterpart of the exact-integration PACL analysis. Exact state--hat pairings \eqref{eq:physical-compatibility} retain the physical cell and patch moments. Shared face rules preserve conservative flux cancellation and lead to the discrete patch identity \eqref{eq:quadrature-patch-identity}, while the discrete Green identity \eqref{eq:discrete-green} preserves constants in the local discrete divergence. The nodal weighted positivity constant \eqref{eq:etaQ}, under the positivity conditions of \Cref{lem:etaQ}, then supplies the admissibility budget for the completed update.

\subsection{Compatible volume and face quadrature}

Let
\begin{equation}\label{eq:QK}
	Q_K(g)=\abs{K}\sum_qw_{Kq}g(\bx_{Kq}),
	\quad w_{Kq}>0,
	\quad \sum_qw_{Kq}=1,
\end{equation}
and let $Q_e^f$ denote the face rules. Fix a nonempty compact
admissibility enforcement set $\mathcal E_K\subseteq\overline K$
containing all volume quadrature nodes and all face nodes at which
a physical flux is evaluated. This set may be a finite collection
of points or the whole cell $\overline K$.
We use an exact cell mass matrix and quadrature only for the
nonlinear volume and face terms.

\begin{assumption}[Compatible quadrature]\label{ass:quadrature}
	Use the local spaces of \Cref{subsec:mesh-space}, the positive normalized volume rule \eqref{eq:QK}, and exact integration for the cell mass matrix. For every cell and active patch pair, require:
	\begin{enumerate}
		\item\label[condition]{quad:unisolvence} the volume nodes determine every function in $\ZK$;
		\item\label[condition]{quad:moments} the discrete state--hat pairings are exact,
		\begin{equation}\label{eq:physical-compatibility}
			Q_K(\psi_av)=\int_K\psi_av\dd\bx,
			\quad \forall v\in\ZK;
		\end{equation}
		\item\label[condition]{quad:faces} the two cells adjacent to a shared face use the same physical face nodes and weights, and a conservative numerical flux;
		\item\label[condition]{quad:green} the discrete Green identity holds,
		\begin{equation}\label{eq:discrete-green}
			Q_K(\partial_r\varphi)
			=\sum_{e\subset\partial K}Q_e^f(n_{K,r}\varphi),
			\quad \forall\varphi\in\ZK,
			\quad r=1,\ldots,d.
		\end{equation}
	\end{enumerate}
\end{assumption}
\Cref{quad:unisolvence} makes $Q_K(\abs{\cdot})$ a norm on the local space. \Cref{quad:moments} identifies the discrete state--hat pairings with the physical weighted moments. \Cref{quad:faces,quad:green} preserve face cancellation and the constant-state identity, respectively. Positivity of the nodal weighted constant is addressed in \Cref{subsec:nodal-estimates}. Face-flux cancellation uses the shared physical nodes and weights, continuity of the patch hat function, and flux conservation. It requires neither positive face weights nor polynomial exactness of the face rule. The discrete Green identity \eqref{eq:discrete-green} serves the separate role of preserving constants.

\begin{remark}[Discrete Green identity on mapped cells]
	\label{rem:green-quadrature}
	In two dimensions, let $T_K:\widehat K\to K$ be the cell map, with $J_K=\det DT_K>0$ and $C_K=J_KDT_K^{-T}$. For $\widehat\varphi=\varphi\circ T_K$, the transformed volume and boundary integrands in \eqref{eq:discrete-green} are
	\[
		(C_K\widehat\nabla\widehat\varphi)_r
		\quad\text{and}\quad
		(C_K\widehat n)_r\widehat\varphi,
	\]
	respectively. Exact reference quadrature for these integrands therefore implies \eqref{eq:discrete-green}, provided that volume and boundary geometry are evaluated consistently.

	For affine triangular $\Pk^k$ cells, volume exactness through degree $k-1$ and edge exactness through degree $k$ suffice. For bilinearly mapped $\Qk^k$ quadrilaterals, volume exactness on $\Qk^k$ and edge exactness through degree $k$ suffice. The mass matrix and state--hat pairings require separate exactness checks.
\end{remark}

\begin{proposition}[Discrete patch identity]\label{prop:quadrature-patch-identity}
	For the meshes and spaces of \Cref{subsec:mesh-space}, suppose \Cref{ass:pou,ass:quadrature} hold. Define the quadrature completed update by replacing the volume and face integrals in \eqref{eq:corrector} by $Q_K$ and $Q_e^f$. Then
	\begin{equation}\label{eq:quadrature-patch-identity}
		\calM_a(u_h^\star)
		=\calM_a(u_h^n)
		+\tau\sum_{K\in\mathcal T_a}
		Q_K\!\left(\FF_{K,\tau}^Q(u_K^n)\cdot\nabla\psi_a\right),
	\end{equation}
	where the left-hand side is the same physical patch mass as in \eqref{eq:mass}. If $\FF_{K,\tau}^Q(c)=\bfv(c)$, constant states make the right-hand correction vanish.
\end{proposition}

\begin{proof}
	The quadrature corrector is defined by replacing the nonlinear volume and face integrals in \eqref{eq:corrector} by $Q_K$ and $Q_e^f$, while retaining the exact mass pairing. Test with $\psi_a|_K$ and sum over the patch. On a shared face the two terms are
	\[
		-\tau Q_e^f\!\left(\widehat{\mathscr F}_{K,e}\psi_a^-
		+\widehat{\mathscr F}_{K',e}\psi_a^+\right)=0.
	\]
	The cancellation uses the same physical nodes and weights on both sides as well as continuity of $\psi_a$. The same cancellation holds on periodically paired faces. The mass terms are physical integrals, giving \eqref{eq:quadrature-patch-identity}. \Cref{quad:moments} allows these moments to be computed with $Q_K$ as well.

	For constant data, use the discrete Green identity \eqref{eq:discrete-green} with $\varphi=\psi_a|_K$. After summing over cells, all shared-face traces cancel and $\sum_KQ_K(\partial_r\psi_a)=0$ for each $r$. Thus $\sum_KQ_K(\bfv(c)\cdot\nabla\psi_a)=0$, proving the final assertion.
\end{proof}

\subsection{Nodal positivity and operator estimates}\label{subsec:nodal-estimates}

The next-step constraints are imposed on the admissibility enforcement set $\mathcal E_K$ defined above. Introduce the discrete patch constant
\begin{equation}\label{eq:etaQ}
	\eta_{aK}^Q
	:=\inf_{\substack{0\ne p\in\ZK\\ p\ge0\text{ on }\mathcal E_K,\, Q_K(|p|)> 0}}
	\frac{Q_K(\psi_ap)}{Q_K(\abs{p})}.
\end{equation}
The volume-node unisolvence in \Cref{ass:quadrature} makes the denominator a norm on $\ZK$.

\begin{lemma}[Positivity of the nodal weighted constant]\label{lem:etaQ}
	Require \Cref{ass:pou}, positive volume weights \eqref{eq:QK}, and volume-node unisolvence in \Cref{quad:unisolvence} of \Cref{ass:quadrature}. Let $\mathcal E_K\subseteq\overline K$ be compact and contain
the volume nodes. For an active pair $(a,K)$, define
	\[
		X_{aK}^{+}
		=\{\bx_{Kq}:\psi_a(\bx_{Kq})>0\}.
	\]
	Then $\eta_{aK}^Q>0$ if and only if
	\[
		p\in\ZK,\quad
		p\geq0\text{ on }\mathcal E_K,\quad
		p=0\text{ on }X_{aK}^{+}
		\quad\Longrightarrow\quad p=0.
	\]
	In particular, either of the following conditions ensures $\eta_{aK}^Q>0$:
	\begin{enumerate}
		\item The hat $\psi_a$ is strictly positive at every
		volume node. In this case,
		\begin{equation}\label{eq:etaQ-lower}
			\eta_{aK}^Q\geq\min_q\psi_a(\bx_{Kq})>0.
		\end{equation}
		\item The nodes in $X_{aK}^{+}$ determine every
		function in $\ZK$.
	\end{enumerate}
	For a fixed volume rule and patch hat, enlarging the admissibility enforcement set cannot decrease the constant:
	\[
		\mathcal E_K\subseteq\widetilde{\mathcal E}_K
		\quad\Longrightarrow\quad
		\eta_{aK}^Q(\mathcal E_K)
		\leq\eta_{aK}^Q(\widetilde{\mathcal E}_K),
	\]
	where the argument indicates the admissibility enforcement set used in \eqref{eq:etaQ}.
\end{lemma}

\begin{proof}
	For $p\geq0$ on $\mathcal E_K$, positivity of the volume weights and nonnegativity of $\psi_a$ imply
	\[
		Q_K(\psi_ap)=0
		\quad\Longleftrightarrow\quad
		p=0\text{ on }X_{aK}^{+}.
	\]
	By volume-node unisolvence, $Q_K(\abs{\cdot})$ is a norm on
the finite-dimensional space $\mathbb Z^{k_K}(K)$.
The condition $p\geq0$ on $\mathcal E_K$ defines a closed subset
of this space. Hence the normalized set
\[
    \bigl\{
        p\in \mathbb Z^{k_K}(K):
        p\geq0\text{ on }\mathcal E_K,\;
        Q_K(\abs{p})=1
    \bigr\}
\]
is compact. The continuous functional $Q_K(\psi_ap)$ therefore has a strictly positive minimum precisely when the stated implication holds.

	If $\psi_a$ is positive at every volume node, then
	\[
		Q_K(\psi_ap)
		\geq\left(\min_q\psi_a(\bx_{Kq})\right)Q_K(p)
		=\left(\min_q\psi_a(\bx_{Kq})\right)Q_K(\abs{p}),
	\]
	which proves \eqref{eq:etaQ-lower}. If $X_{aK}^{+}$ determines $\ZK$, vanishing on $X_{aK}^{+}$ implies $p=0$, so the second sufficient condition follows from the general criterion.

	Finally, enlarging $\mathcal E_K$ restricts the set of polynomials in the infimum \eqref{eq:etaQ}, while leaving its numerator and denominator unchanged. The infimum therefore cannot decrease.
\end{proof}

For the inverse and discrete divergence estimates, write $\norm{v}_{\mathcal E,K}
    =\sup_{x\in\mathcal E_K}\abs{v(x)}$ and 
$\norm{v}_{1,Q,K}:=Q_K(\abs{v})$.
For continuous functions, the supremum is attained because
$\mathcal E_K$ is compact. Use componentwise suprema for vectors.

\begin{lemma}[Nodal inverse and local-operator bounds]\label{lem:quadrature-operator}
	Under \Cref{ass:quadrature}, let $\mathcal E_K\subseteq\overline K$ be compact and contain all nodes entering the volume and face evaluations. There are finite constants $\iota_K^Q\geq 1$ and $C_{G,K}^Q$ such that
	\begin{align}
		\norm{v}_{\mathcal E,K}
		&\leq\frac{\iota_K^Q}{\abs{K}}Q_K(\abs{v}),\quad \forall v\in\ZK,
		\label{eq:nodal-inverse}\\
		\norm{\mathscr R_K^Q[\mat W]}_{\mathcal E,K}
		&\leq\frac{C_{G,K}^Q}{h_K}\max_r\norm{\vect w_r}_{\mathcal E,K}.
		\label{eq:nodal-discrete-div}
	\end{align}
	Here $\mathscr R_K^Q$ is the componentwise discrete divergence defined by
	\[
		(\mathscr R_K^Q[\mat W],\varphi)_K
		=Q_K\!\left(\sum_r\vect w_r\partial_r\varphi\right)
		-\sum_e Q_e^f\!\left(\sum_r n_{K,r}\vect w_r\varphi\right),\quad\forall\varphi\in\ZK.
	\]
	For a scalar flux Lipschitz with constant $\alpha$ on an interval $I$, and $v,w\in\ZK$ taking values in $I$ at every node of $\mathcal E_K$, the discrete divergence of the flux difference satisfies
	\begin{align}
		\norm{\mathscr R_K^Q[\bfv(v)-\bfv(w)]}_{\mathcal E,K} & \leq\frac{\alpha C_{G,K}^Q}{h_K}\norm{v-w}_{\mathcal E,K}, \label{eq:nodal-G-infinity}\\
		\norm{\mathscr R_K^Q[\bfv(v)-\bfv(w)]}_{1,Q,K} & \leq \frac{\alpha C_{G,K}^Q\iota_K^Q}{h_K}\norm{v-w}_{1,Q,K}. \label{eq:nodal-G-L1}
	\end{align}
\end{lemma}

\begin{proof}
	Positive weights and volume-node unisolvence make $Q_K(\abs{\cdot})$ a norm on the local space. The evaluation norm $\norm{\cdot}_{\mathcal E,K}$ is finite on this space because local functions are continuous on $\overline K$, and it is a norm because $\mathcal E_K$ contains the volume nodes. Norm equivalence gives the inverse bound. The discrete divergence depends on finitely many volume and face flux samples and takes values in the local space. Its composition with evaluation on $\mathcal E_K$ is bounded, giving the local-operator bound. On fixed reference rules and uniformly regular elements, the constants are uniform after the displayed factors $\abs{K}^{-1}$ and $h_K^{-1}$ are extracted.

	Apply \eqref{eq:nodal-discrete-div} to $\bfv(v)-\bfv(w)$ and use the pointwise flux bound to prove \eqref{eq:nodal-G-infinity}. Also $Q_K(\abs{z})\leq\abs{K}\norm{z}_{\mathcal E,K}$ for any local polynomial $z$. Applying this to the discrete divergence difference, followed by \eqref{eq:nodal-inverse} for $v-w$, proves \eqref{eq:nodal-G-L1}. The discrete Green identity \eqref{eq:discrete-green} additionally gives $\mathscr R_K^Q[\bfv(c)]=0$.
\end{proof}
These constants involve all evaluation nodes, including face nodes. An estimate using only the volume-node maximum would not, by itself, control flux evaluations at extrapolated face values.

\subsection{Discrete admissibility results}\label{subsec:discrete-admissibility}

For the discrete local volume flux define
\[
	\mathscr F_{aK,\tau}^Q(v)
	=Q_K\!\left(\FF_{K,\tau}^Q(v)\cdot\nabla\psi_a\right).
\]

\begin{theorem}[Scalar patchwise weak admissibility with quadrature]\label{thm:quadrature-scalar}
	Use the quadrature completed update of \Cref{prop:quadrature-patch-identity} on the meshes and spaces of \Cref{subsec:mesh-space}, under \Cref{ass:pou,ass:quadrature}. Require $\eta_{aK}^Q>0$ for every active pair and $\FF_{K,\tau}^Q(c)=\bfv(c)$. Suppose that the discrete local flux is well defined for every input in $\mathcal V_K^Q=\{v\in\ZK:m\leq v\leq M\text{ on }\mathcal E_K\}$ and satisfies
	\begin{equation}\label{eq:quadrature-lipschitz}
		\abs{\mathscr F_{aK,\tau}^Q(v)-\mathscr F_{aK,\tau}^Q(w)}
		\leq L_{aK}^Q(\tau)Q_K(\abs{v-w})
	\end{equation}
	for all $v,w\in\mathcal V_K^Q$. If
	\begin{equation}\label{eq:quadrature-cfl}
		\tau L_{aK}^Q(\tau)\leq\eta_{aK}^Q
		\quad\text{for every active pair},
	\end{equation}
	and the old state satisfies $m\leq u_h^n\leq M$ on $\mathcal E_K$, then the \emph{physical} patch averages of the raw completed update satisfy
	\[
		m\leq\avg{u_h^\star}_a\leq M
		\quad\text{for every patch }a.
	\]
\end{theorem}

\begin{proof}
	For a patch input collection define
	\[
		H_a^Q(\{v_K\})=\frac{1}{\mathfrak m_a}
		\sum_{K\in\mathcal T_a}
		\left[Q_K(\psi_av_K)+\tau\mathscr F_{aK,\tau}^Q(v_K)\right].
	\]
	Suppose one input is increased from $v_K$ to $w_K$ in the nodal order, and put $p=w_K-v_K\ge0$ on $\mathcal E_K$. Positivity of the volume weights gives $Q_K(\abs{p})=Q_K(p)$, even though $p$ need not be nonnegative between nodes. Therefore
	\[
		\begin{aligned}
			\mathfrak m_a\left[H_a^Q(\ldots,w_K,\ldots)
			-H_a^Q(\ldots,v_K,\ldots)\right]
			&\geq Q_K(\psi_ap)-\tau L_{aK}^Q(\tau)Q_K(p)\\
			&\geq\left(\eta_{aK}^Q-\tau L_{aK}^Q(\tau)\right)Q_K(p)\ge0.
		\end{aligned}
	\]
	The constant-state part of \Cref{prop:quadrature-patch-identity} gives $H_a^Q(\{c\})=c$. Replace the inputs one cell at a time to compare the old collection with $\{m\}$ and $\{M\}$. The same proposition and exact state--hat pairings identify $H_a^Q(\{u_K^n\})$ with the physical patch average, proving the claim.
\end{proof}

For cRKDG and both ADER predictors, \Cref{lem:quadrature-operator} verifies the spatial input to the proofs in \Cref{subsec:flux-estimates}: use $\norm{\cdot}_{1,Q,K}$ and $\norm{\cdot}_{\mathcal E,K}$, with $C_{k_K}$ replaced by $C_{G,K}^Q\iota_K^Q$ and $C_{k_K}^{(\infty)}$ by $C_{G,K}^Q$. The same stage and Picard inductions give the stated growth factors. The predictor-admissibility condition must hold at all evaluation nodes. For the Taylor flux, the coefficient estimate is applied to the actually evaluated CK maps: polynomial derivatives are evaluated at the nodes and the physical-flux derivatives are taken at the nodal old states. Compactness of the nodally bounded coefficient set and the argument of \Cref{lem:ck-lipschitz} give the corresponding $Q_K$-norm bound. A different flux-interpolation procedure requires its own coefficient-map bound.

\begin{corollary}[Euler positivity with quadrature]\label{cor:euler-quadrature}
	For ideal-gas Euler with $\gamma>1$, use the componentwise quadrature update of \Cref{prop:quadrature-patch-identity} on the meshes and spaces of \Cref{subsec:mesh-space}, under \Cref{ass:pou,ass:quadrature}. Require $\eta_{aK}^Q>0$ for every active pair and the nodal old-state bounds
	\[
		\vect U_K^n(x)\in\GE,\quad
		\rho_K^n(x)\geq\delta\bar\rho_K^n,
		\quad \forall x\in\mathcal E_K,
		\quad 0<\delta\leq 1.
	\]
	Use $\iota_K^Q,C_{G,K}^Q$ from \Cref{lem:quadrature-operator}. 
	Define the nodal old-state wave speed and local CFL number by
	\[
		\alpha_{K,Q}^n
		=\max_{\bx\in\mathcal E_K}
		\left(
			\abs{\vect v_K^n(\bx)}
			+\sqrt{\frac{\gamma p_K^n(\bx)}{\rho_K^n(\bx)}}
		\right),
		\quad
		\nu_{K,Q}^n=\frac{\tau\alpha_{K,Q}^n}{h_K}.
	\]
	With $\vect Z_1=\stateT_K\vect U_K^n$, set 
		\begin{align}
			\beta_{aK}^Q
			=\frac{h_K}{\abs{K}}Q_K(\abs{\nabla\psi_a}_1),\quad\text{and}\quad
			\mathcal D_K^Q
			=\bigcup_{\bx\in\mathcal E_K}
			\left\{
				\vect Z:
				\abs{\vect Z-\vect Z_1(\bx)}_\infty
				\leq \frac12(\vect Z_1(\bx))_\rho
			\right\}.
		\end{align}
	Use either the explicit cRKDG predictor \eqref{eq:rkstage}--\eqref{eq:rkflux} with $\sum_j b_j=1$, or the fixed-$J$ Picard predictor
	\eqref{eq:ader-picard}--\eqref{eq:ader-flux} under \Cref{ass:ader-time}, with spatial discrete divergence $\mathscr R_K^Q$. Set $C_{\mathrm{pred}}$ and $\omega_{\mathrm{abs}}$ by \eqref{eul:eq:temporal-constants}. 
	Suppose the numerical face flux is finite at every evaluated predictor trace and that
	\begin{equation}\label{eq:euler-Q-CFL}
		\begin{aligned}
			C_{\mathrm{pred}} C_{G,K}^Q
			C_{\gamma,d}\iota_K^Q\nu_{K,Q}^n
			&\leq \delta/2
			&&\text{for every cell }K,\\
			C_d\,\omega_{\mathrm{abs}}\,\beta_{aK}^Q
			C_{\gamma,d}\iota_K^Q\nu_{K,Q}^n
			&\leq \eta_{aK}^Q/4
			&&\text{for every active pair }(a,K).
		\end{aligned}
	\end{equation}
	Then every required predictor evaluation has positive density and
	\[
		\calM_a(\vect U_h^\star)/\mathfrak m_a\in\GE
		\quad\text{for every patch }a.
	\]
\end{corollary}

\begin{proof}
	The old state is required to be admissible only at the admissibility enforcement set $\mathcal E_K$. Nevertheless, positive volume weights and exact cell averages give $\bar{\vect U}_K^n=Q_K(\vect U_K^n)/\abs{K}\in\GE$, so its normalization is well defined. With $\vect Z_1=\stateT_K\vect U_K^n$,
	\[
		\frac{Q_K(\vect Z_1)}{\abs{K}}=(1,\vect0^\top,1/2)^\top.
	\]
	Apply \eqref{eq:nodal-inverse} to the nodally positive density and energy components. This gives \eqref{eul:eq:old-normalized-bound} on $\mathcal E_K$ with $\iota_K^Q$. The momentum bound again follows from nodal admissibility.
	At each volume node, the ideal-gas energy identity gives
	\[
		2E_K^n
		\leq b_\gamma^2(\alpha_{K,Q}^n)^2\rho_K^n,
		\quad
		b_\gamma^2=1+\frac{2}{\gamma(\gamma-1)}.
	\]
	Applying the positive rule $Q_K$ and using exactness of the cell means gives $c_K\leq b_\gamma\alpha_{K,Q}^n$. 
	The same volume weights give $\abs{\bar{\vect v}_K}\leq\alpha_{K,Q}^n$. 
	Each state in $\mathcal D_K^Q$ lies in a neighborhood centered at $\vect Z_1(\bx)$ for some $\bx\in\mathcal E_K$.
	Applying the pointwise argument in the proof of \Cref{eul:lem:predictor-flux} at this node, with $\iota_K^Q$ in place of $\iota_K$, gives
	\[
		\max_{\vect Z\in\mathcal D_K^Q}\max_{1\leq s\leq d}
		\abs{\vect g_s(\vect Z)}_\infty
		\leq C_{\gamma,d}\,\iota_K^Q\alpha_{K,Q}^n.
	\]

	The old weighted contribution $\vect W_{aK}^Q=Q_K(\psi_a\vect Z_1)/\abs{K}$ satisfies
	\[
		(\vect W_{aK}^Q)_\rho\geq\eta_{aK}^Q,
		\quad
		\ell_{\vect w}(\vect W_{aK}^Q)
		\geq\frac{\eta_{aK}^Q}{2}(1+\abs{\vect w}^2).
	\]
	Indeed, for each fixed $\vect w\in\R^d$, both
$p=(\vect Z_1)_\rho$ and $p=\ell_{\vect w}(\vect Z_1)$ belong to
$\ZK$ and are positive on $\mathcal E_K$. Hence
\eqref{eq:etaQ} gives $Q_K(\psi_a p)\geq\eta_{aK}^Q Q_K(p)$
for each choice of $p$. Their quadrature means are $1$ and
$(1+\abs{\vect w}^2)/2$, respectively, by the normalization of
$\vect Z_1$. The first restriction in \eqref{eq:euler-Q-CFL} and the nodal discrete divergence bound give a predictor excursion of at most $\delta/2$ on $\mathcal E_K$, by the induction in \Cref{eul:lem:predictor-correction}. 
At each evaluation node $\bx$, this excursion is at most $(\vect Z_1(\bx))_\rho/2$, so every predictor value lies in $\mathcal D_K^Q$. 
For the completed step, the nodal flux bound and the positive volume weights give
\[
	\abs{\vect r_{aK}^Q}_\infty
	\leq \omega_{\mathrm{abs}}\,\beta_{aK}^Q
	C_{\gamma,d}\iota_K^Q\nu_{K,Q}^n .
\]
The second restriction in \eqref{eq:euler-Q-CFL} and \Cref{eul:lem:budget} yield $\vect W_{aK}^Q+\vect r_{aK}^Q\in\GE$. 
Undoing the normalization and summing the physical cell contributions as in \Cref{eul:thm:exact} gives the result by \Cref{prop:quadrature-patch-identity}.
\end{proof}

The factor $\iota_K^Q$ in \eqref{eq:euler-Q-CFL} comes from the nodal inverse estimate \eqref{eq:nodal-inverse} on $\mathcal E_K$. Once the physical patch averages are admissible, the patch limiter and the physical-mean correction of \Cref{sec:patch-limiting} apply without a change of definition.

\subsection{Quadrature rules and admissibility enforcement sets}\label{subsec:quadrature-choices}

\Cref{tab:quadrature-choices} gives practical configurations satisfying \Cref{ass:quadrature} and ensuring $\eta_{aK}^Q>0$ for standard vertex hats.
Write $k=k_K\geq1$, and choose
\begin{equation}\label{eq:quadrature-check-set}
	\mathcal E_K
	=(X_K^{\mathrm{vol}} \cup X_K^{\mathrm{face}}) \sqcup X_K^{\mathrm{add}}.
\end{equation}
The final scaling enforces admissibility on $\mathcal E_K$.

Here $\mathrm{GL}_n$ and $\mathrm{LGL}_n$ denote the $n$-point Gauss--Legendre and Gauss--Lobatto node sets, with their standard quadrature weights.
Tensor products are formed on the reference cell or face and mapped to the physical element.
Tensor-product cells use $\Qk^k$ spaces, with mapped spaces on bilinear quadrilaterals; triangles use $\Pk^k$ spaces.
For $d=1$, the face nodes are the interval endpoints.

\begin{table}[htb]
	\centering
	\caption{Quadrature rules and additional admissibility
		enforcement points.}
	\label{tab:quadrature-choices}
	\begin{tabular}{lccc}
		\toprule
		Cell & Volume nodes & Face nodes & $X_K^{\mathrm{add}}$\\
		\midrule
		Affine tensor-product
		& $\mathrm{GL}_{k+1}^{\otimes d}$
		& $\mathrm{GL}_{k+1}^{\otimes(d-1)}$
		& $\varnothing$\\
		Bilinear quadrilateral
		& $\mathrm{GL}_{k+1}^{\otimes2}$
		& $\mathrm{GL}_{k+1}$
		& $\varnothing$\\
		Affine triangle
		& Interior cubature nodes
		& $\mathrm{GL}_{k+1}$
		& $\varnothing$\\
		\midrule
		Affine tensor-product
		& $\mathrm{LGL}_{k+2}^{\otimes d}$
		& $\mathrm{LGL}_{k+2}^{\otimes(d-1)}$
		& $\varnothing$\\
		Bilinear quadrilateral
		& $\mathrm{LGL}_{k+2}^{\otimes2}$
		& $\mathrm{LGL}_{k+2}$
		& $\varnothing$\\
		\midrule
		Affine tensor-product, $k\geq2$
		& $\mathrm{LGL}_{k+1}^{\otimes d}$
		& $\mathrm{LGL}_{k+1}^{\otimes(d-1)}$
		& chosen so that $\eta_{aK}^Q>0$ \\
		Bilinear quadrilateral, $k\geq3$
		& $\mathrm{LGL}_{k+1}^{\otimes2}$
		& $\mathrm{LGL}_{k+1}$
		& chosen so that $\eta_{aK}^Q>0$ \\
		\bottomrule
	\end{tabular}
\end{table}

The mass matrix is computed exactly; when the volume rule is exact for the mass products, the same rule may be used for its assembly. On mapped cells, the physical Jacobian is included in the volume weights.
Adjacent cells use the same physical face nodes and weights.
The triangular cubature has positive weights, interior nodes, and total-degree exactness at least $2k$; see \cite{10.2307/43693493,10.1016/j.cam.2023.115451}.
All listed configurations integrate the state--hat pairings exactly and satisfy \eqref{eq:discrete-green}.
For the Gauss and triangular configurations, the interior volume nodes determine the local polynomial space and all vertex hats are strictly positive there, so \Cref{lem:etaQ} gives $\eta_{aK}^Q>0$.

\begin{remark}[Gauss--Lobatto quadrature]\label{rem:lobatto-enforcement}
	For the $(k+1)$-point Lobatto configurations in \Cref{tab:quadrature-choices}, additional nodes, such as the corresponding Gauss nodes, are needed to ensure $\eta_{aK}^Q>0$.
	Using only the Lobatto nodes can give $\eta_{aK}^Q=0$. 
	The justification and a one-dimensional counterexample are given in \Cref{app:lobatto-enforcement}.
	Alternatively, $(k+2)$-point Lobatto quadrature in each volume and face coordinate works for all $k\geq1$ on affine tensor-product cells and bilinear quadrilaterals, with $X_K^{\mathrm{add}}=\varnothing$.
	In this case, the positive volume nodes of each vertex hat determine $\Qk^k$, so \Cref{lem:etaQ} applies.
\end{remark}

\begin{remark}[Relation to convex-decomposition point sets]
	The set $\mathcal E_K$ plays a role analogous to the point sets in the Zhang--Shu framework \cite{10.1016/j.jcp.2009.12.030} and its optimal cell-average decompositions \cite{10.1137/23M1549365}: admissibility at selected points yields admissibility of an updated average under a CFL restriction. 
	Here the updated average is a weighted patch average, and $\mathcal E_K$ enters the sufficient CFL condition through $\eta_{aK}^Q$ and the local operator bounds.
\end{remark}

\section{Numerical experiments}\label{sec:numerical}

The experiments examine accuracy at active scalar bounds, scalar shock formation, Euler positivity near vacuum and in strong-shock flows, and the number of rejected trial steps in the practical admissibility safeguard.

\subsection{Implementation}

Numerical experiments are implemented using the MFEM finite element library \cite{MFEM410}. We use degree two in space: $\Pk^2$ on intervals and triangles and mapped $\Qk^2$ on quadrilaterals. The two-dimensional hybrid meshes consist of affine triangular and bilinear quadrilateral elements. We use 6-point positive interior quadrature with algebraic degree of exactness 4 for triangles and a tensor product of 3-point Gauss--Legendre rules for bilinear quadrilaterals. Each edge uses 3-point Gauss--Legendre quadrature. These choices integrate the mass matrix and the state--hat pairings exactly and satisfy the discrete Green identity \eqref{eq:discrete-green}. For scalar equations, we choose the admissibility enforcement set $\mathcal E_K = \overline{K}$. For systems, we use \eqref{eq:quadrature-check-set} with $X_K^{\mathrm{add}} = \varnothing$.

Time integration uses the third-order Heun-cRKDG method (cRKDG3) \cite[Eqn.~(2.10)]{10.1137/23M158629X} and third-order fixed-Picard ADER (ADER3) with three-point Gauss--Legendre time quadrature and three predictor updates. All examples use patch limiting followed by cellwise Zhang--Shu scaling (PACL). OE/RIOE damping \cite{10.1090/mcom/3998, 10.1016/j.jcp.2025.113769} may be applied after the mean shift in \Cref{step:mean-shift} and before the cellwise Zhang--Shu scaling in \Cref{step:cell-scaling}. The corrected cell average from patch limiting remains fixed throughout this step.

Following \Cref{rem:predictor-pressure}, we evaluate the Euler physical flux using the actual conservative state, allowing negative predictor pressure, and use $\max(p,0)$ in the Einfeldt wave-speed estimates for the HLL numerical flux. The adiabatic constant is taken to be $\gamma=1.4$ unless otherwise stated. Absolute density and pressure safety guards in the final scaling are $10^{-13}$. The accepted and rejected time step attempts are reported in \Cref{subsec:num-rejections}. For Euler computations, we recommend choosing the trial time step using the maximum characteristic speed over all volume and face quadrature points in each cell. With $X_K^{\mathrm{add}}=\varnothing$ as in these experiments, this maximum is the old-state wave speed $\alpha_{K,Q}^n$ in \eqref{eq:euler-Q-CFL}. If a required predictor density or raw patch admissibility check fails, the attempt is rejected and repeated with half the step size.

\subsection{1D examples}\label{subsec:num-1d}

\begin{figure}[!ht]
	\centering
	\includegraphics[width=0.4\linewidth]{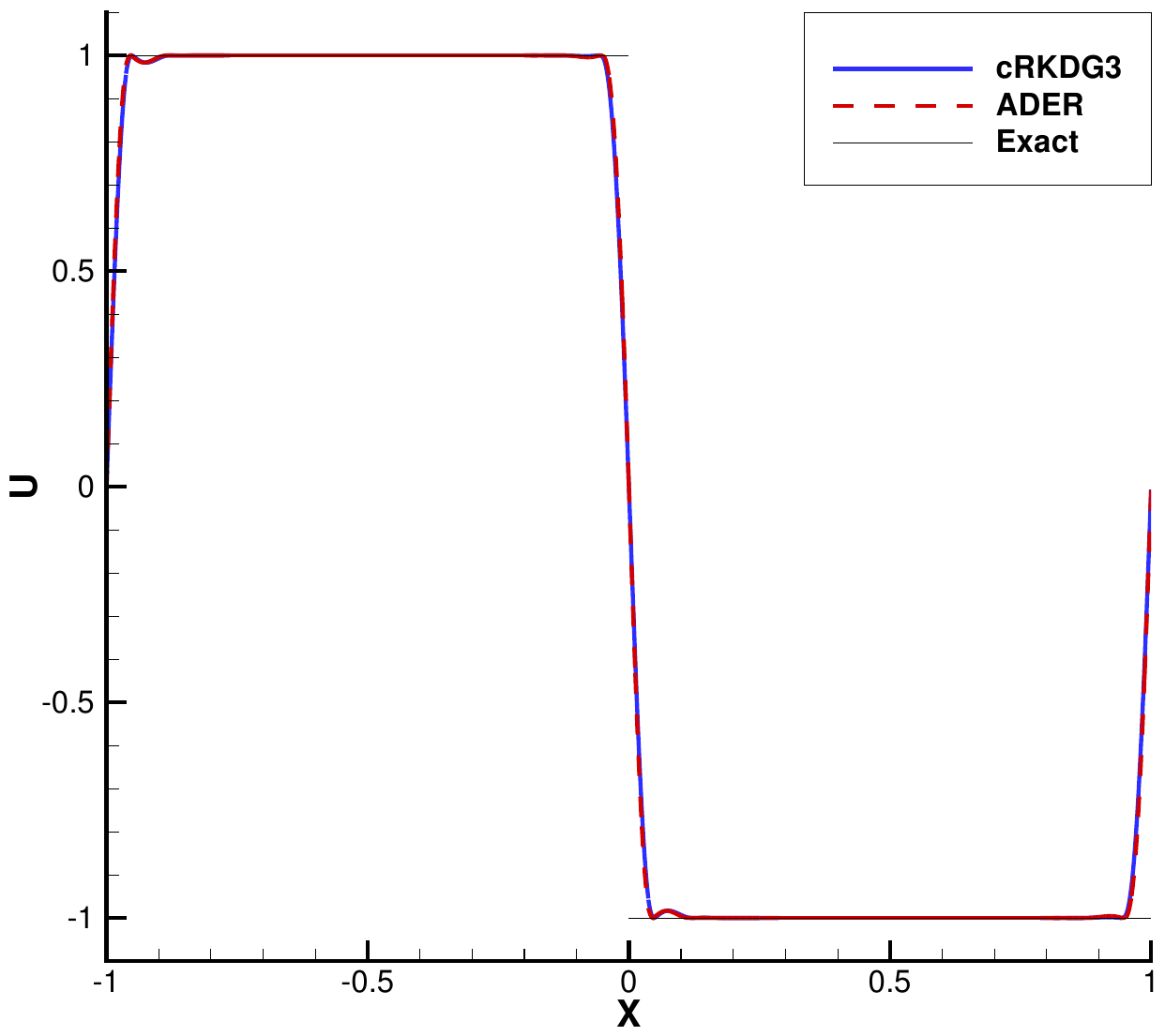}
	\caption{Linear square-wave transport after 50 periods, $N=160$, PACL without oscillation suppression.}
	\label{fig:num-square}
\end{figure}

\begin{example}[Long-time square-wave transport]\label{ex:num-square}
	We solve $u_t+u_x=0$ on the periodic interval $[-1,1]$ with square-wave initial data $u_0=1$ on $[-1,0]$ and $u_0=-1$ on $(0,1]$ \cite[Example~6.1]{10.1016/j.jcp.2009.12.030}. We use $N=160$ and $T=100$, corresponding to 50 periods. \Cref{fig:num-square} shows the numerical solutions and the exact solution. The profiles show the accumulated diffusion and oscillations after 50 periods. The numerical solutions stay well within the bound.
\end{example}

\begin{figure}[!ht]
	\centering
	\includegraphics[width=0.4\linewidth]{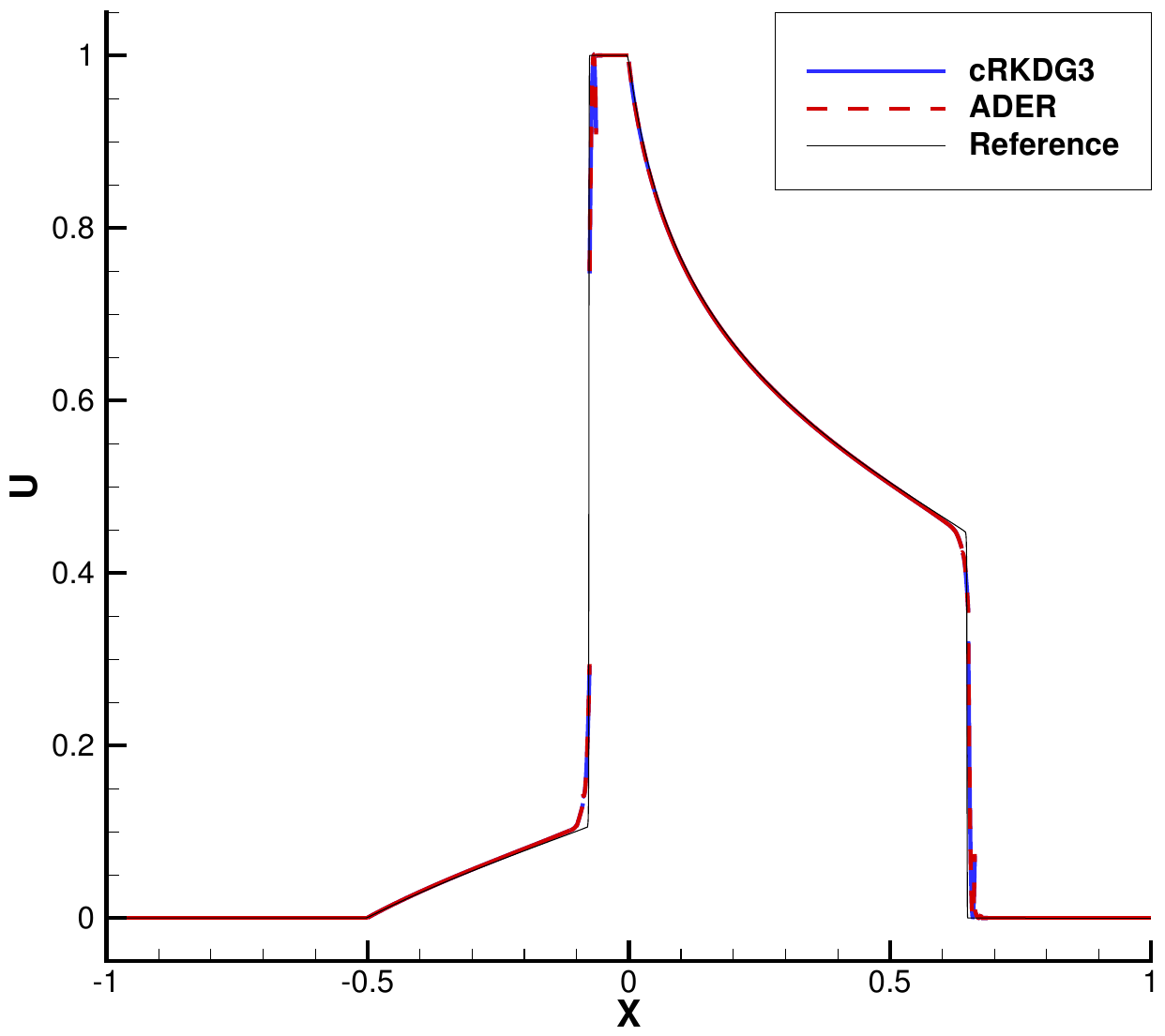}
	\caption{Buckley--Leverett, PACL and OE, $N=160$, $T=0.4$.}
	\label{fig:num-buckley}
\end{figure}

\begin{example}[Buckley--Leverett equation]\label{ex:num-buckley}
	We solve the Buckley--Leverett equation with flux and initial data
	\[
		f(u)=\frac{4u^2}{4u^2+(1-u)^2},\quad
		u_0(x)=\begin{cases}
			1,&-1/2\leq x\leq0,\\
			0,&\text{elsewhere}.
		\end{cases}
	\]
	The initial data form a compact pulse with bounds $[0,1]$ \cite[Example~5.3]{10.1016/j.jcp.2009.12.030}. We use $N=160$ and $T=0.4$. The computational domain is $[-1,1]$ with periodic boundaries. \Cref{fig:num-buckley} shows the numerical solutions. The reference solution is computed by SSP-RK3-FD-WENO5 with 4096 grid cells. The numerical solution stays within the bound \([0,1]\), which is ensured by our numerical routines.
\end{example}

\begin{figure}[!ht]
	\centering
	\begin{subfigure}{0.3\linewidth}
		\centering
		\includegraphics[width=\linewidth]{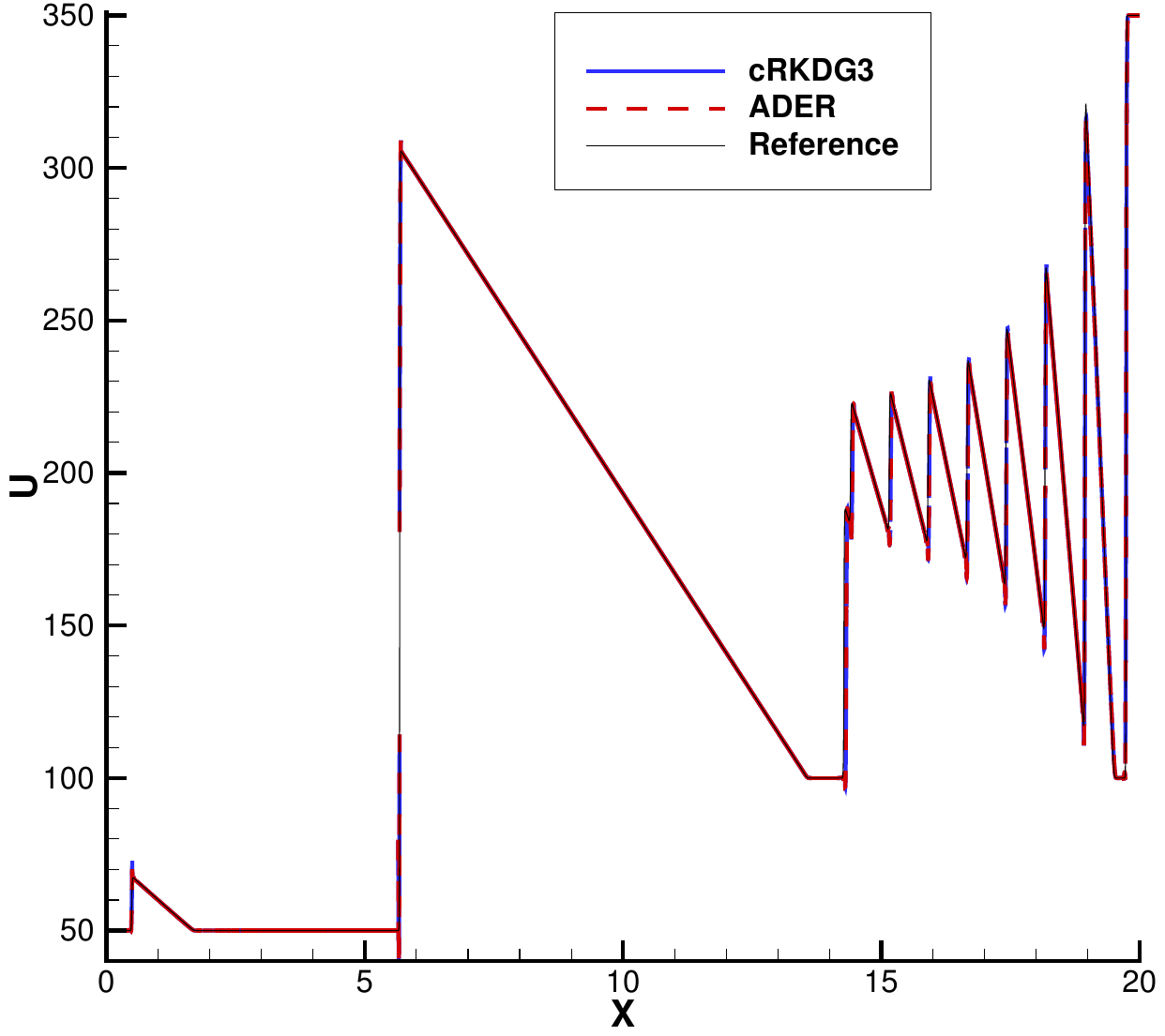}
		\caption{\(t = 0.5\)}
	\end{subfigure}
	\begin{subfigure}{0.3\linewidth}
		\centering
		\includegraphics[width=\linewidth]{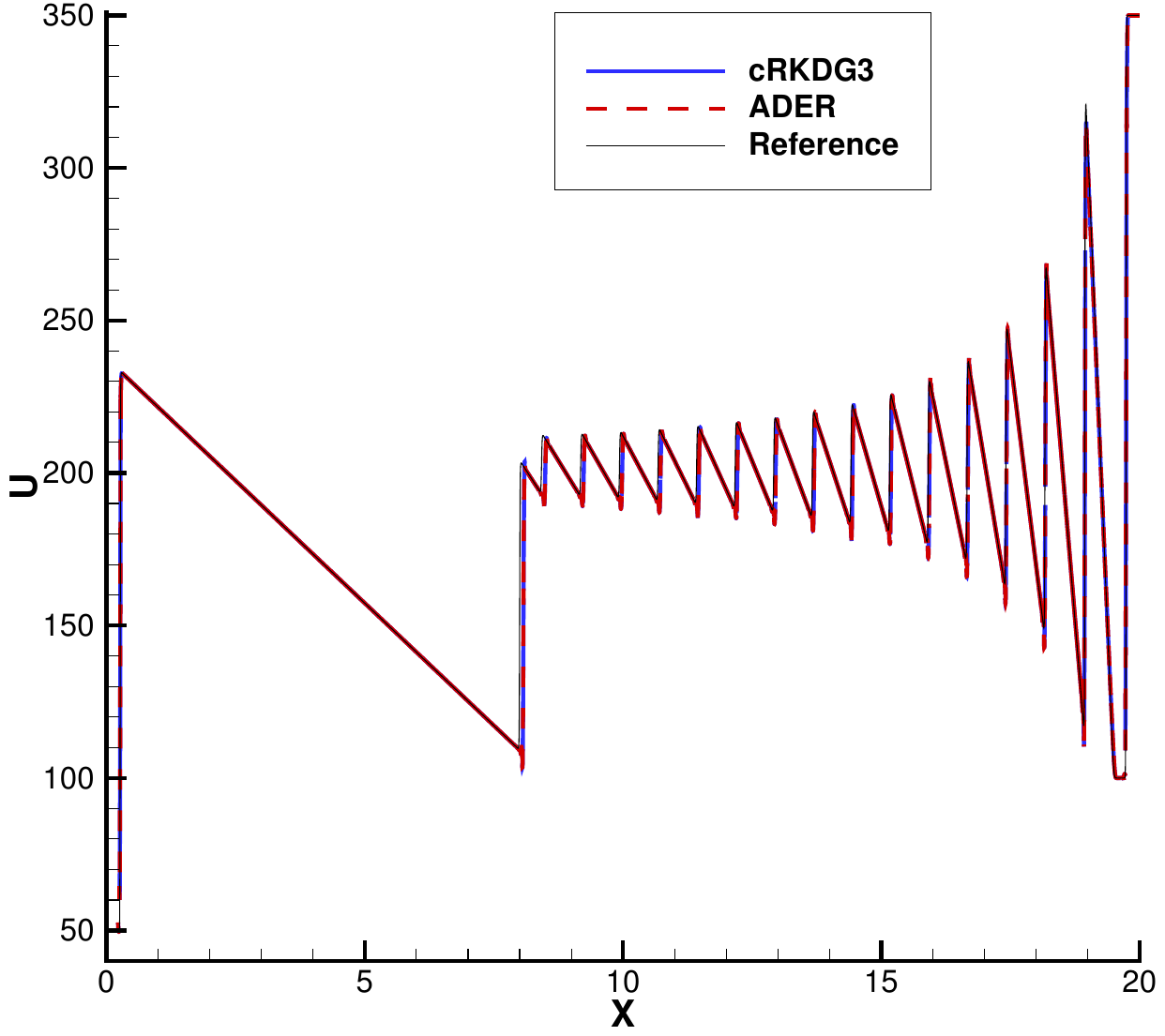}
		\caption{\(t = 1.0\)}
	\end{subfigure}
	\begin{subfigure}{0.3\linewidth}
		\centering
		\includegraphics[width=\linewidth]{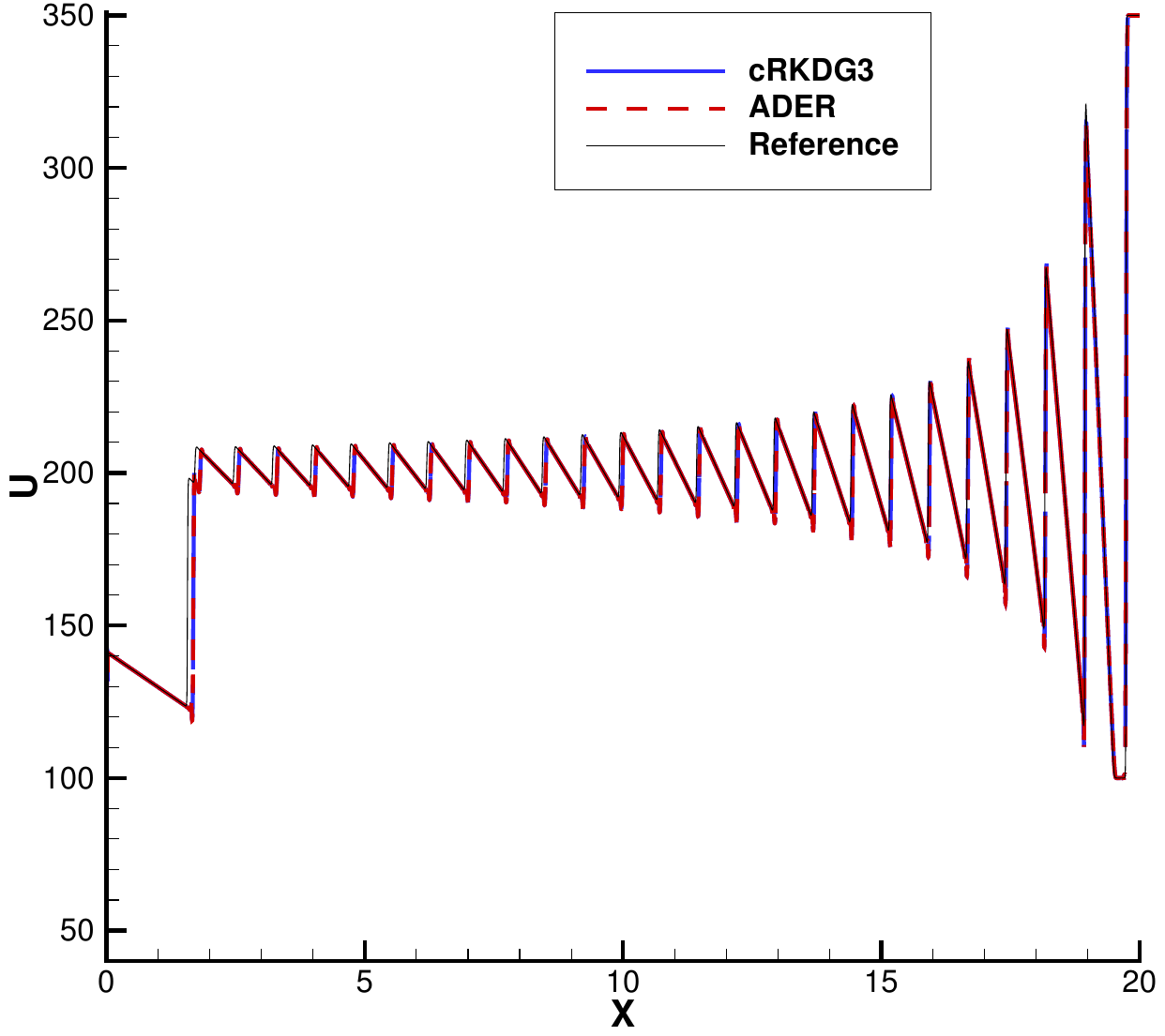}
		\caption{\(t = 1.5\)}
	\end{subfigure}
	\caption{Traffic flow with a piecewise-quadratic flux, $N=800$, PACL and OE.}
	\label{fig:num-traffic56}
\end{figure}

\begin{example}[Traffic flow with a piecewise-quadratic flux]\label{ex:num-traffic56}
	We model traffic flow on a road of length $20$ km with a piecewise-quadratic flux. Time is measured in hours, and the density $\rho$ is measured in vehicles per kilometer, with bounds $[0,350]$. The flux is
	\[
		f(\rho)=\begin{cases}
			-0.4\rho^2+100\rho,&\rho<50,\\
			-0.1\rho^2+15\rho+3500,&50\leq\rho<100,\\
			-0.024\rho^2-5.2\rho+4760,&\rho\geq100.
		\end{cases}
	\]
	The initial density is $50$ on $[0,10)$ and $350$ on $[10,15]$. On $[15,20]$, it decreases linearly between the points $(15,350)$, $(18.571,100)$, $(19.286,50)$ and $(20,0)$ \cite{10.1016/j.trb.2007.08.004}. At $x=0$, the prescribed density is $0$ for the first ten minutes, $75$ for the next ten minutes, and $50$ thereafter. At $x=20$, a traffic signal imposes a repeating three-minute cycle: the density is $0$ for two minutes and $350$ for one minute \cite[Example~5.6]{10.1016/j.jcp.2009.12.030}. We use $N=800$ cells and apply PACL and OE with both cRKDG3 and ADER3. The reference solution is computed by SSP-RK3-FD-WENO5 with 4096 grid cells. \Cref{fig:num-traffic56} shows the solutions at $T=0.5$, $1.0$ and $1.5$ hours. Our numerical solution aligns well with the literature.
\end{example}

\begin{figure}[!ht]
	\centering
	\begin{subfigure}{0.3\linewidth}
		\centering
		\includegraphics[width=\linewidth]{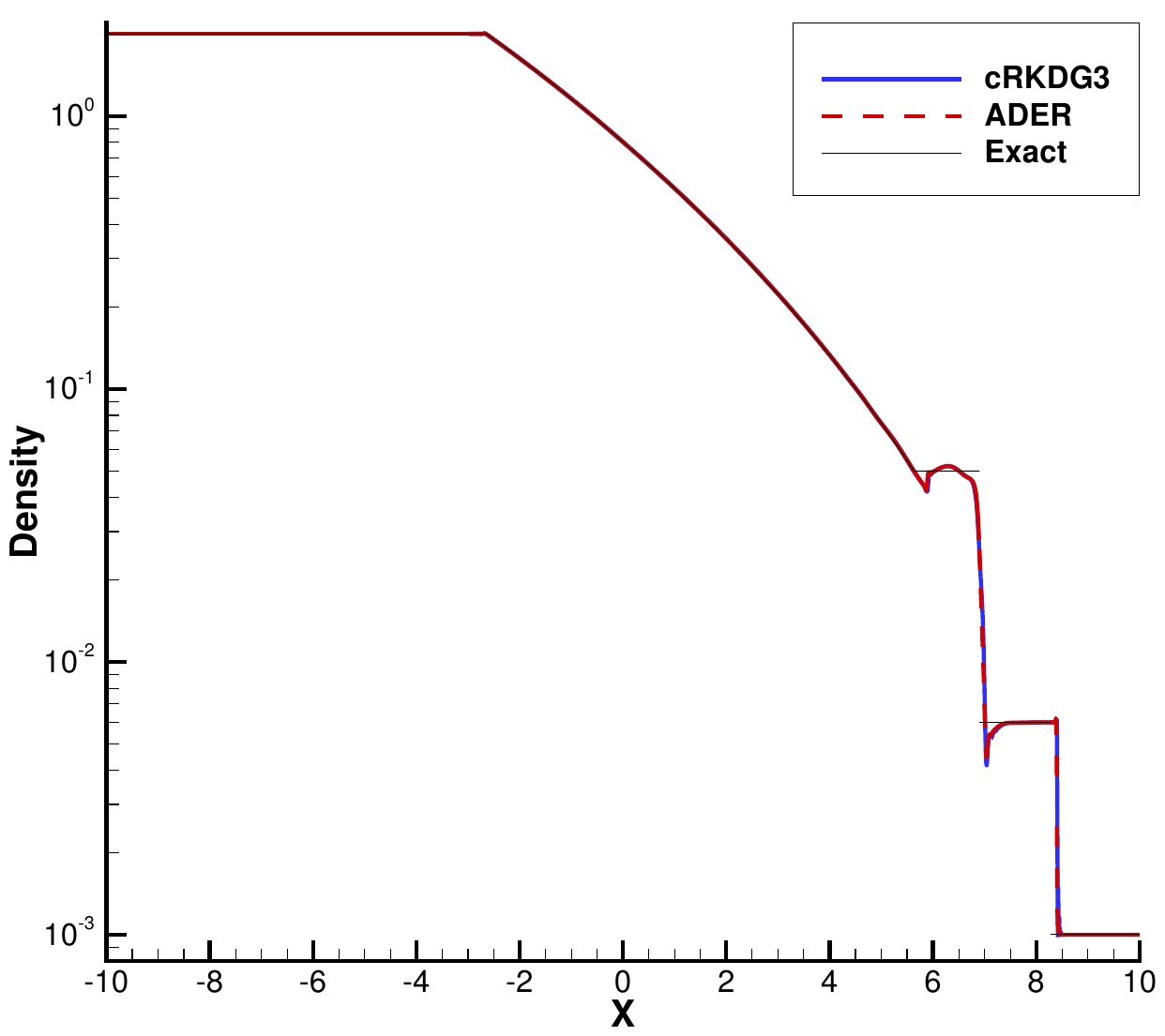}
		\caption{Density}
	\end{subfigure}
	\begin{subfigure}{0.3\linewidth}
		\centering
		\includegraphics[width=\linewidth]{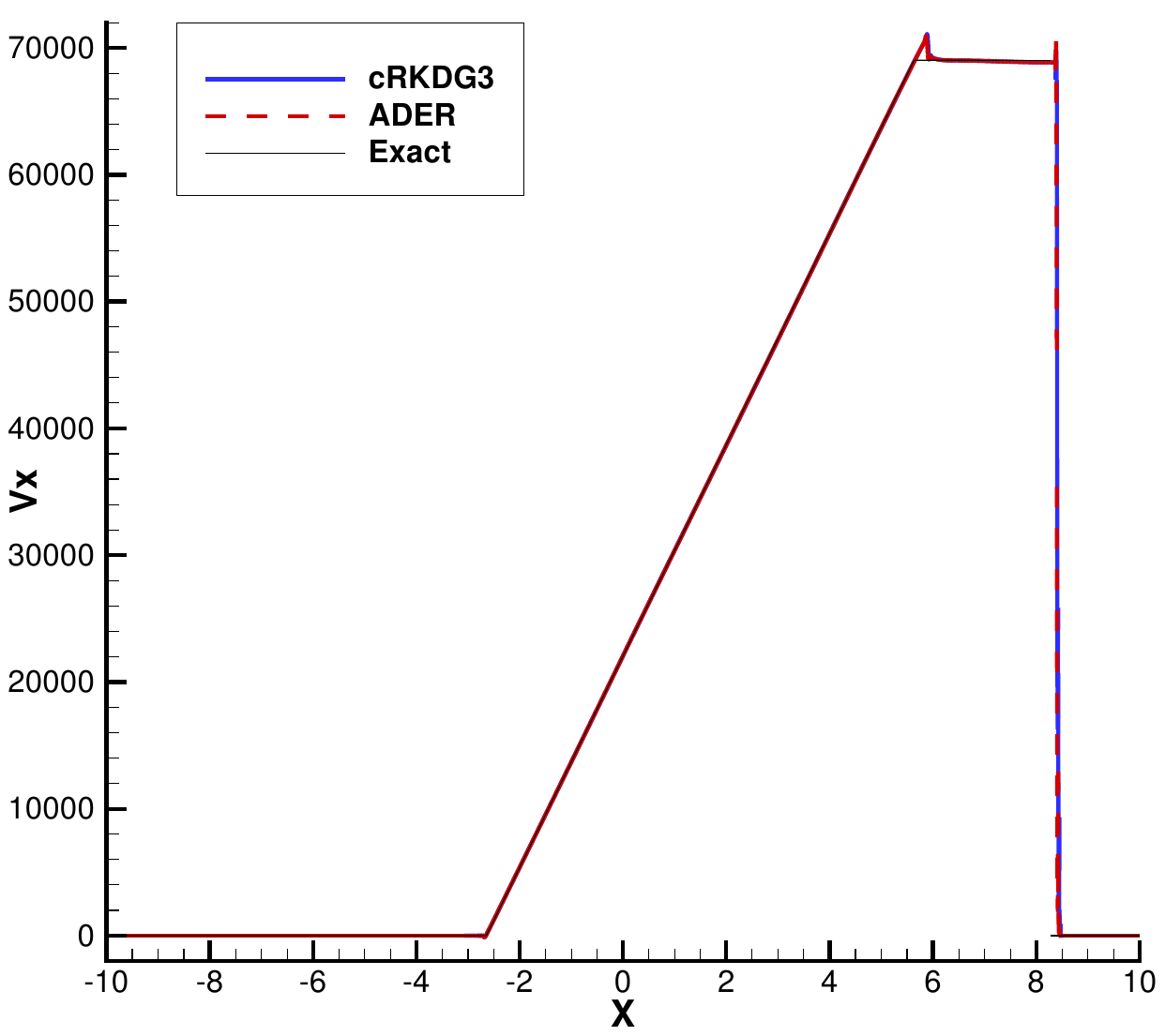}
		\caption{Velocity}
	\end{subfigure}
	\begin{subfigure}{0.3\linewidth}
		\centering
		\includegraphics[width=\linewidth]{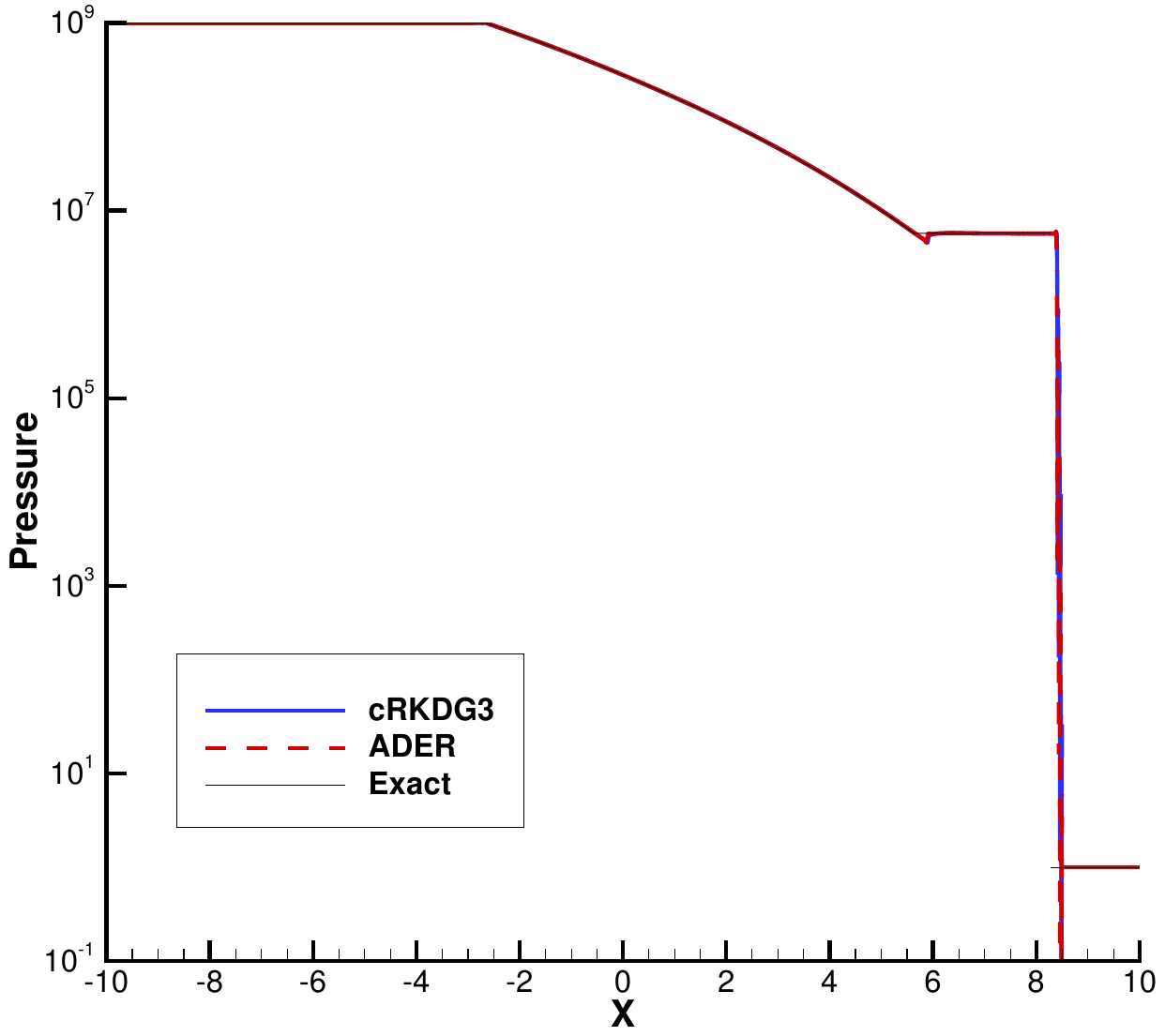}
		\caption{Pressure}
	\end{subfigure}
	\caption{LeBlanc shock tube, $N=800$, $T=10^{-4}$, PACL without oscillation suppression.}
	\label{fig:num-leblanc}
\end{figure}

\begin{example}[LeBlanc shock tube]\label{ex:num-leblanc}
	The LeBlanc shock tube is posed on $[-10,10]$ with initial states $(\rho,v,p)_L=(2,0,10^9)$ and $(\rho,v,p)_R=(10^{-3},0,1)$, separated at $x=0$ \cite[Section~4.3]{10.1016/j.jcp.2010.08.016}. No oscillation-suppressing limiter is used. \Cref{fig:num-leblanc} compares density, pressure and velocity with the exact Riemann solution. This verifies the robustness of our method in the vicinity of the extremely large pressure ratio.
\end{example}

\begin{figure}[!ht]
	\centering
	\begin{subfigure}{0.3\linewidth}
		\centering
		\includegraphics[width=\linewidth]{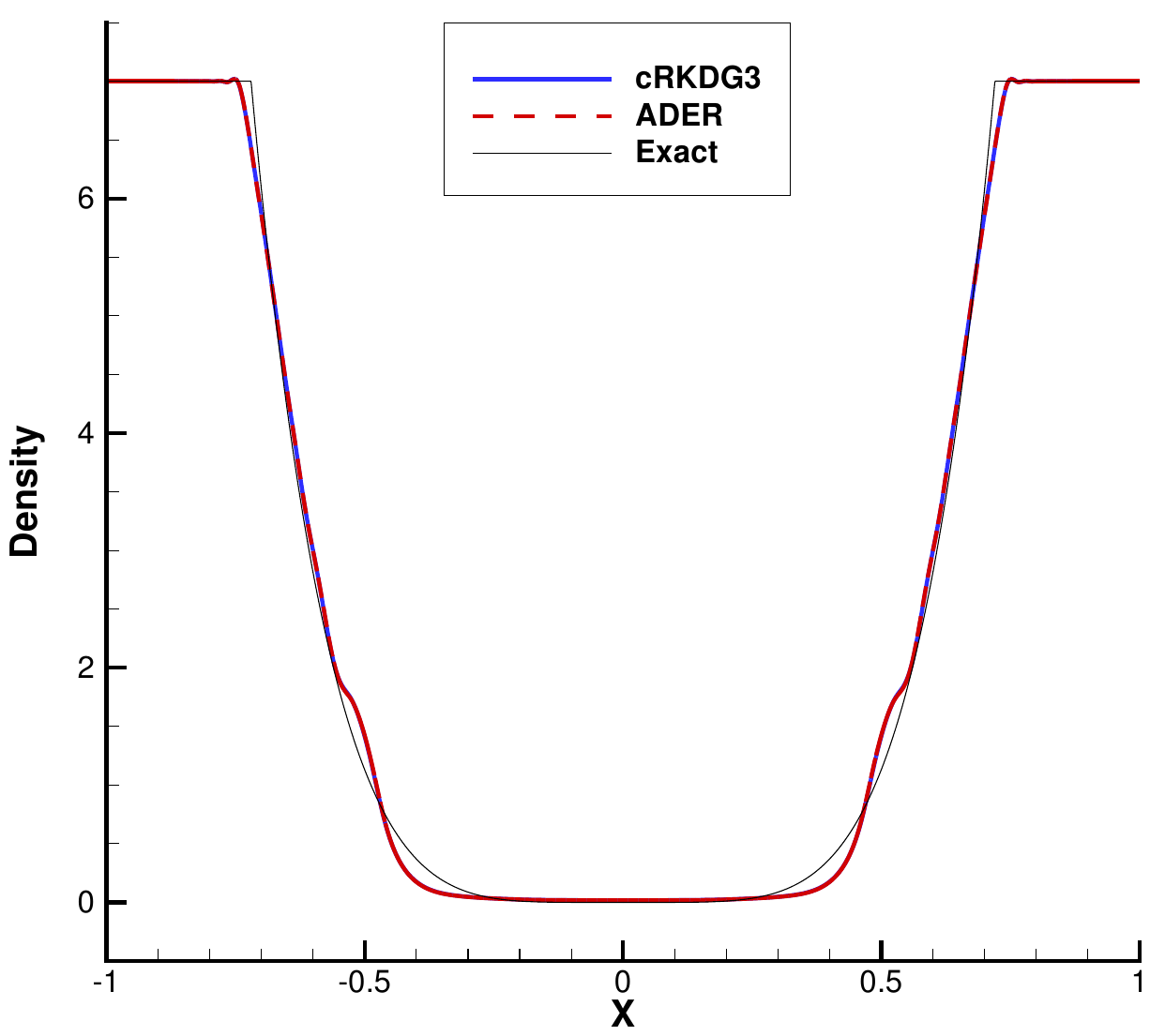}
		\caption{Density}
	\end{subfigure}
	\begin{subfigure}{0.3\linewidth}
		\centering
		\includegraphics[width=\linewidth]{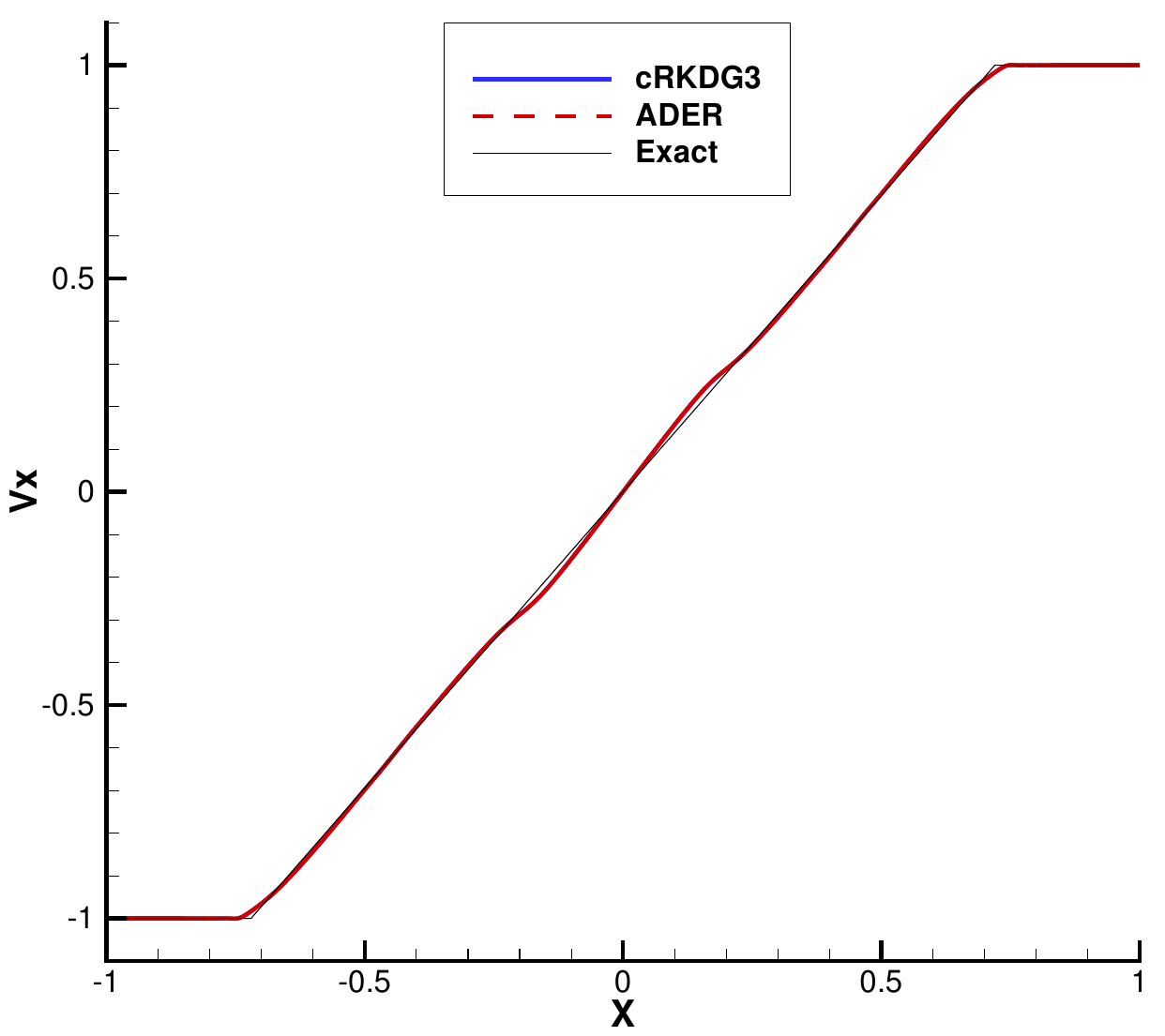}
		\caption{Velocity}
	\end{subfigure}
	\begin{subfigure}{0.3\linewidth}
		\centering
		\includegraphics[width=\linewidth]{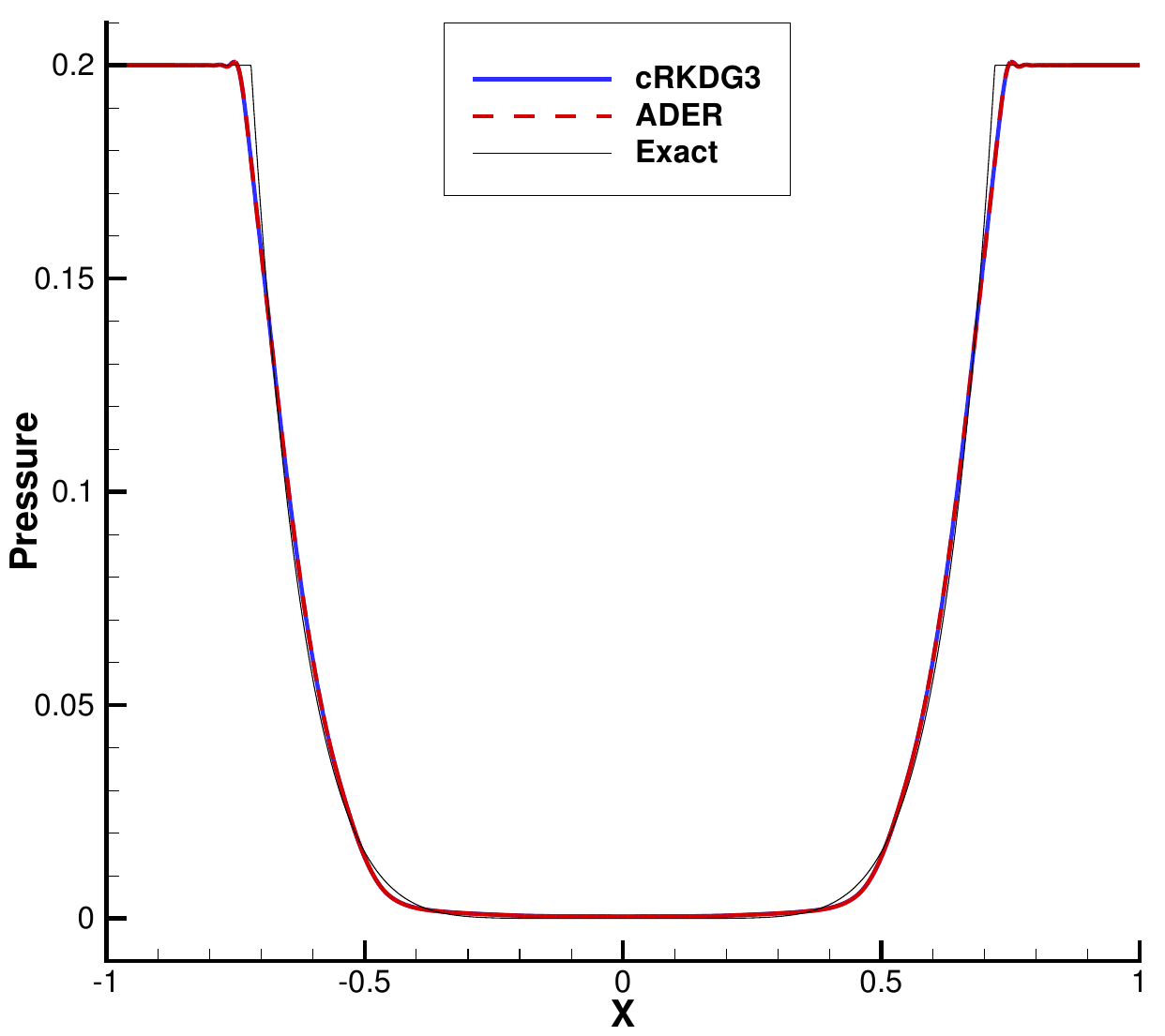}
		\caption{Pressure}
	\end{subfigure}
	\caption{One-dimensional double rarefaction, PACL and OE, $N=200$, $T=0.600$.}
	\label{fig:num-rarefaction}
\end{figure}

\begin{example}[Double rarefaction]\label{ex:num-rarefaction}
	We consider a Riemann problem consisting of two rarefaction waves and an intermediate vacuum state on $[-1,1]$, with initial density $\rho_0=7$, pressure $p_0=0.2$, and velocity $v_0=-1$ for $x<0$ and $v_0=1$ for $x>0$ \cite[Section~4.3]{10.1016/j.jcp.2010.08.016}. We take $N=200$ and $T=0.600$. The exact solution has a vacuum at the origin. \Cref{fig:num-rarefaction} compares the numerical solutions with the exact Riemann solution. Our numerical method remains robust. 
\end{example}

\subsection{2D examples}\label{subsec:num-2d}

The smooth tests use PACL. For the two smooth tests below, an $N\times N$ rectangular background grid is used, and the rectangles in half of the domain are split along a diagonal into right triangles. This gives $N^2/2$ quadrilaterals and $N^2$ triangles, with $\Qk^2$ and $\Pk^2$ approximation spaces, respectively. Both cRKDG3 and ADER3 use the same mesh, and each uniform refinement produces four children per cell.

\begin{example}[Two-dimensional smooth Burgers equation]\label{ex:num-burgers2d-acc}
	We consider the two-dimensional Burgers equation with smooth initial data
	\[
		u_t+(u^2/2)_x+(u^2/2)_y=0,\quad u_0(x,y)=0.5+\sin(\pi(x+y)),
	\]
	on the periodic square $[-1,1]^2$, with bounds $[-0.5,1.5]$ \cite[Example~5.5]{10.1016/j.jcp.2009.12.030}. At $T=0.05$, the smooth solution satisfies $u(x,y,t)=u_0(x-tu,y-tu)$. \Cref{tab:num-burgers} reports cRKDG3 and ADER3 errors on hybrid meshes. The $L^1$ orders remain close to three under refinement.
\end{example}

\begin{table}[htb]
	\centering
	\caption{Two-dimensional Burgers equation at $T=0.05$ on hybrid meshes.}
	\label{tab:num-burgers}
	\begin{tabular}{crrrrrrr}
		\toprule
		& & \multicolumn{2}{c}{$L^1$ error}
		& \multicolumn{2}{c}{$L^2$ error}
		& \multicolumn{2}{c}{$L^\infty$ error} \\
		\cmidrule(lr){3-4}
		\cmidrule(lr){5-6}
		\cmidrule(lr){7-8}
		Method & $N$ & Error & Order & Error & Order & Error & Order \\
		\midrule

		\multirow{4}{*}{cRKDG3}
		& 40  & 1.411e-04 & --    & 1.340e-04 & --    & 9.684e-04 & --    \\
		& 80  & 1.780e-05 & 2.987 & 1.797e-05 & 2.899 & 1.389e-04 & 2.802 \\
		& 160 & 2.246e-06 & 2.986 & 2.377e-06 & 2.918 & 1.961e-05 & 2.824 \\
		& 320 & 2.878e-07 & 2.964 & 3.172e-07 & 2.906 & 2.623e-06 & 2.902 \\
		\midrule

		\multirow{4}{*}{ADER3}
		& 40  & 1.412e-04 & --    & 1.341e-04 & --    & 9.693e-04 & --    \\
		& 80  & 1.781e-05 & 2.987 & 1.798e-05 & 2.899 & 1.390e-04 & 2.802 \\
		& 160 & 2.247e-06 & 2.987 & 2.378e-06 & 2.919 & 1.962e-05 & 2.824 \\
		& 320 & 2.878e-07 & 2.965 & 3.172e-07 & 2.906 & 2.625e-06 & 2.902 \\

		\bottomrule
	\end{tabular}
\end{table}

\begin{example}[Smooth low-density Euler transport]\label{ex:num-euler-acc}
	We advect a smooth density profile with constant velocity and pressure on the periodic square $[0,1]^2$:
	\[
		\rho_0(x,y)=1+0.99\sin(2\pi(x+y)),\quad (v_1,v_2)=(1,1),\quad p=1.
	\]
	The exact density is $\rho(x,y,t)=1+0.99\sin(2\pi(x+y-2t))$, with minimum $0.01$ \cite[Example~6.6]{10.1007/s10915-011-9472-8}. \Cref{tab:num-euler-acc} compares cRKDG3 and ADER3 density errors at $T=0.1$ on hybrid meshes. The finest $L^1$ and $L^2$ orders are approximately $2.966$.
\end{example}

\begin{table}[!ht]
	\centering
	\caption{Low-density Euler transport at $T=0.1$ on hybrid $\Pk^2/\Qk^2$ meshes: density errors.}
	\label{tab:num-euler-acc}
	\begin{tabular}{crrrrrrr}
		\toprule
		& & \multicolumn{2}{c}{$L^1$ error}
		& \multicolumn{2}{c}{$L^2$ error}
		& \multicolumn{2}{c}{$L^\infty$ error} \\
		\cmidrule(lr){3-4}
		\cmidrule(lr){5-6}
		\cmidrule(lr){7-8}
		Method & $N$ & Error & Order & Error & Order & Error & Order \\
		\midrule

		\multirow{4}{*}{cRKDG3}
		& 20  & 2.506e-04 & --    & 3.784e-04 & --    & 2.807e-03 & --    \\
		& 40  & 3.253e-05 & 2.945 & 4.903e-05 & 2.948 & 3.499e-04 & 3.004 \\
		& 80  & 4.208e-06 & 2.951 & 6.328e-06 & 2.954 & 4.313e-05 & 3.020 \\
		& 160 & 5.384e-07 & 2.966 & 8.096e-07 & 2.966 & 5.339e-06 & 3.014 \\
		\midrule

		\multirow{4}{*}{ADER3}
		& 20  & 2.506e-04 & --    & 3.785e-04 & --    & 2.808e-03 & --    \\
		& 40  & 3.253e-05 & 2.945 & 4.904e-05 & 2.948 & 3.499e-04 & 3.004 \\
		& 80  & 4.208e-06 & 2.951 & 6.329e-06 & 2.954 & 4.314e-05 & 3.020 \\
		& 160 & 5.384e-07 & 2.966 & 8.097e-07 & 2.966 & 5.340e-06 & 3.014 \\

		\bottomrule
	\end{tabular}
\end{table}

\begin{figure}[!ht]
	\centering
	\begin{subfigure}{0.3\linewidth}
		\centering
		\includegraphics[width=\linewidth]{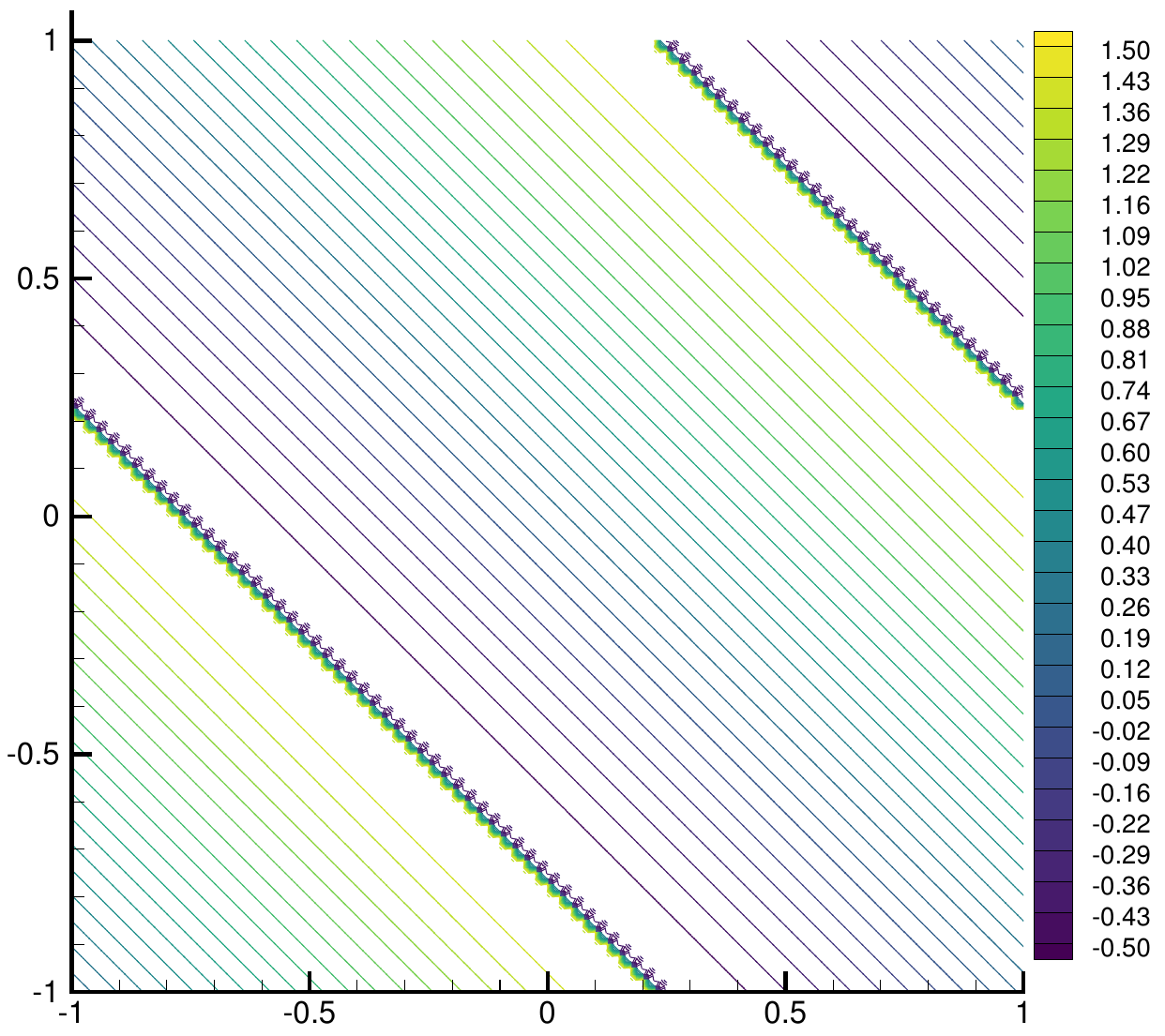}
		\caption{cRKDG3, \(T=0.23\)}
	\end{subfigure}
	\begin{subfigure}{0.3\linewidth}
		\centering
		\includegraphics[width=\linewidth]{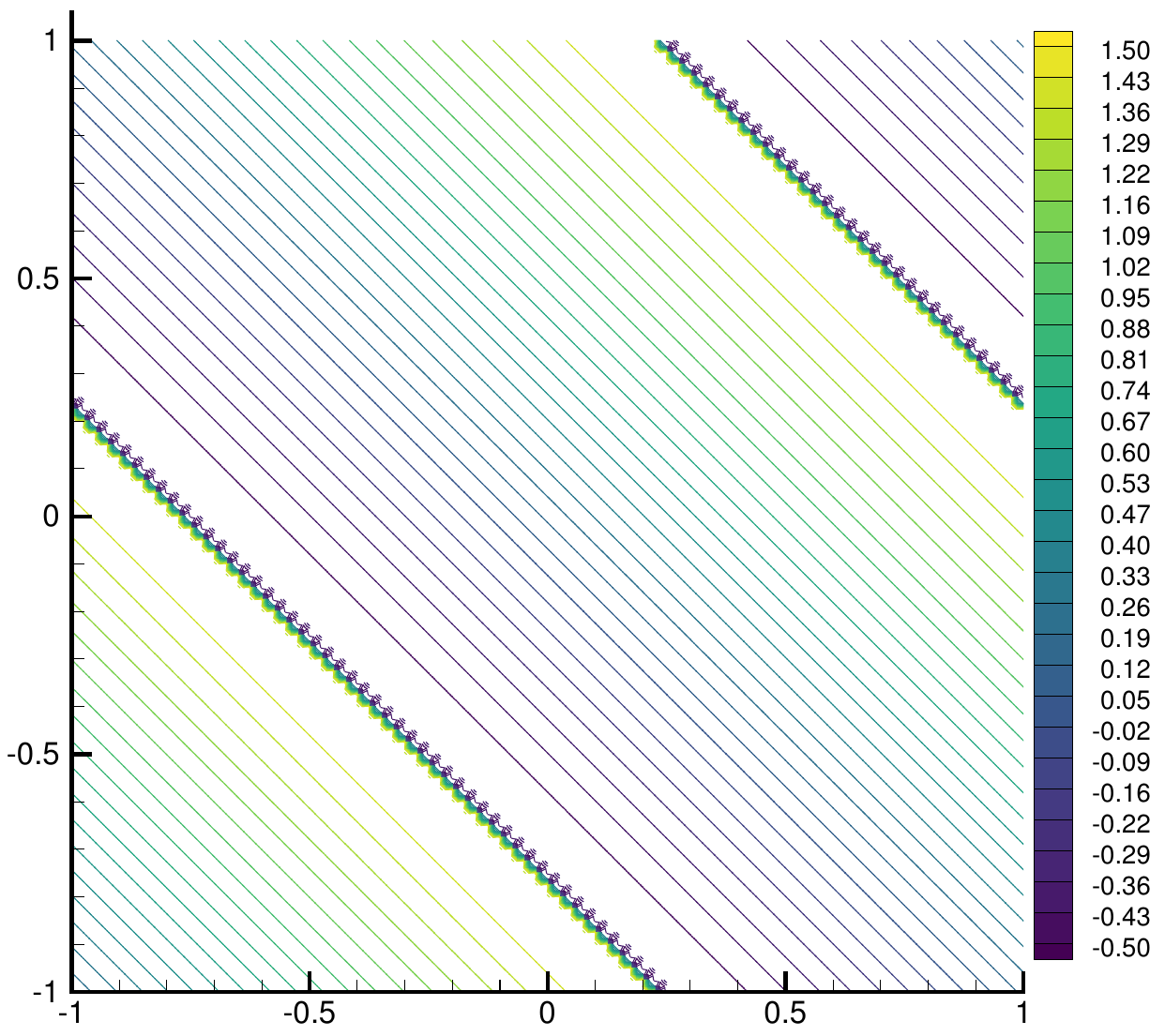}
		\caption{ADER3, \(T=0.23\)}
	\end{subfigure}
	\begin{subfigure}{0.3\linewidth}
		\centering
		\includegraphics[width=\linewidth]{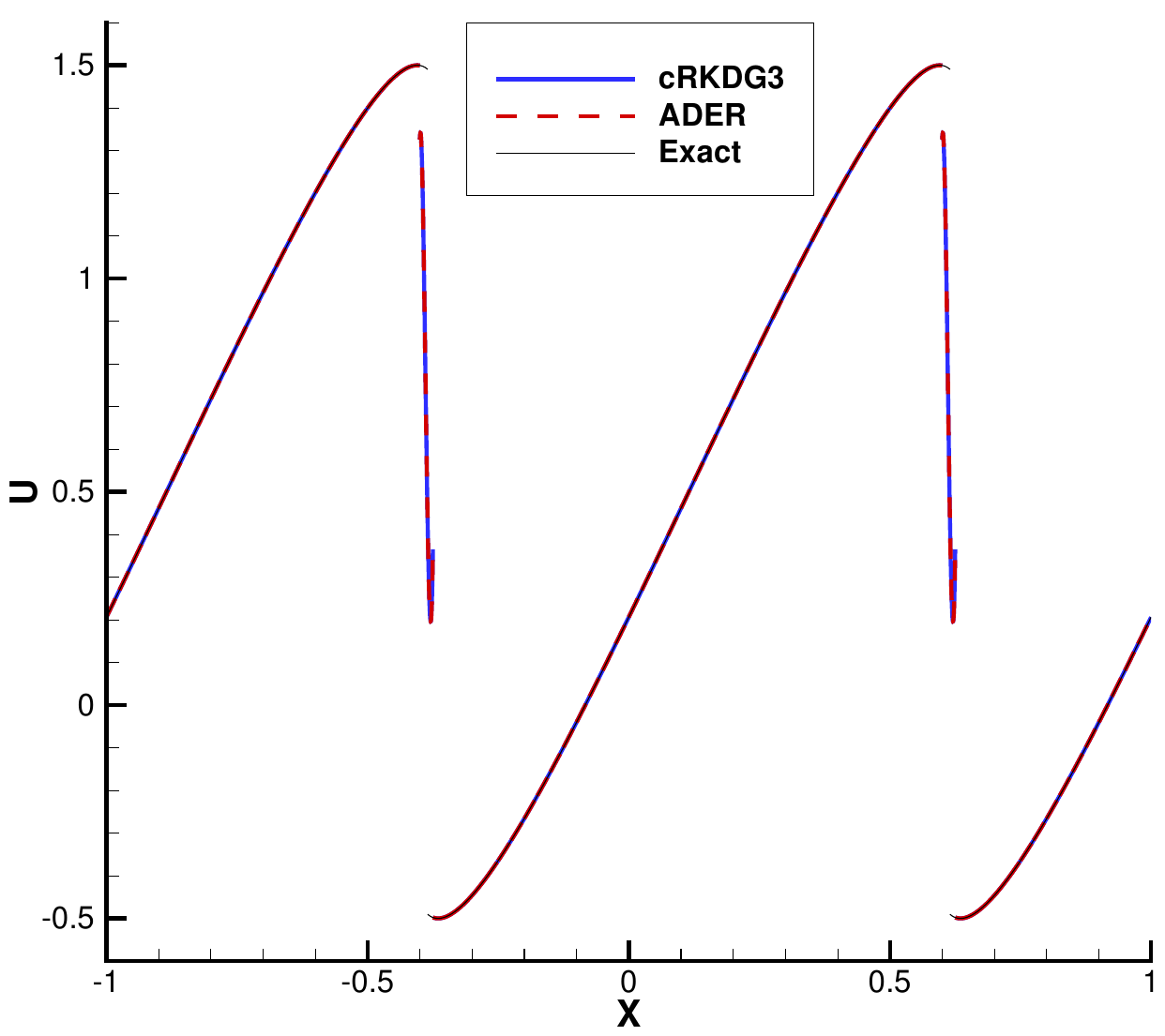}
		\caption{Diagonal slice, \(T=0.23\)}
	\end{subfigure}
	\\
	\begin{subfigure}{0.3\linewidth}
		\centering
		\includegraphics[width=\linewidth]{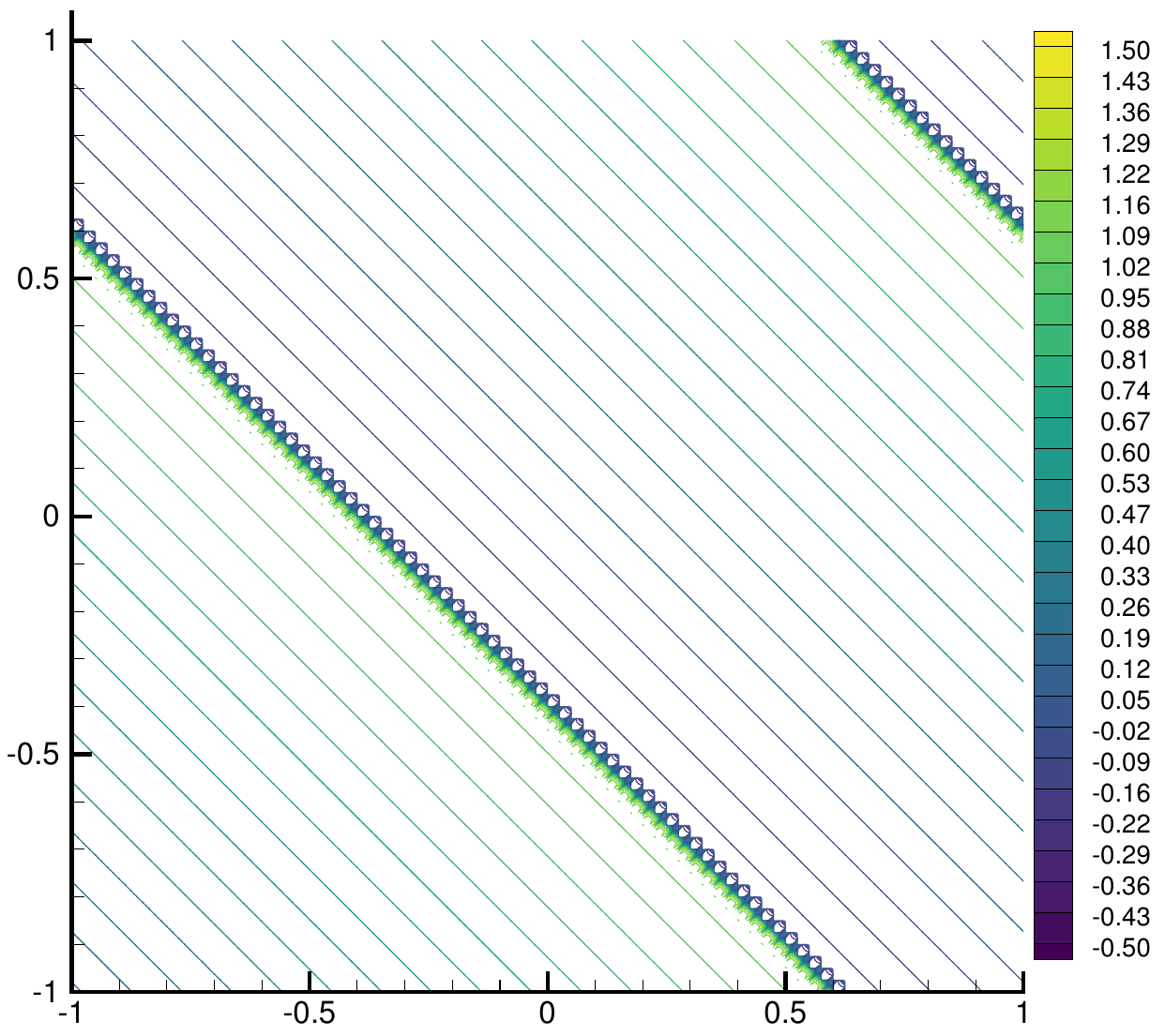}
		\caption{cRKDG3, \(T=0.6\)}
	\end{subfigure}
	\begin{subfigure}{0.3\linewidth}
		\centering
		\includegraphics[width=\linewidth]{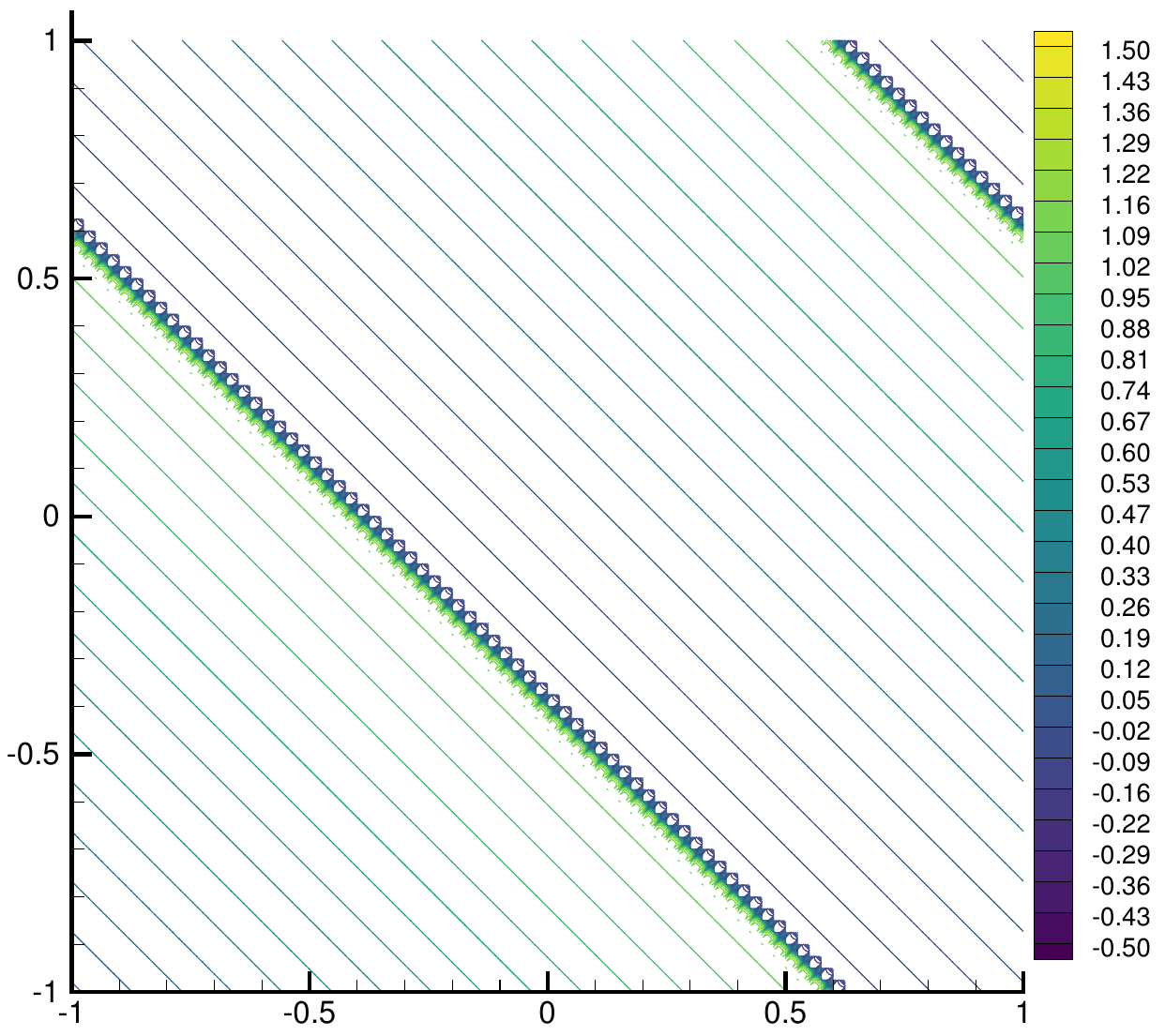}
		\caption{ADER3, \(T=0.6\)}
	\end{subfigure}
	\begin{subfigure}{0.3\linewidth}
		\centering
		\includegraphics[width=\linewidth]{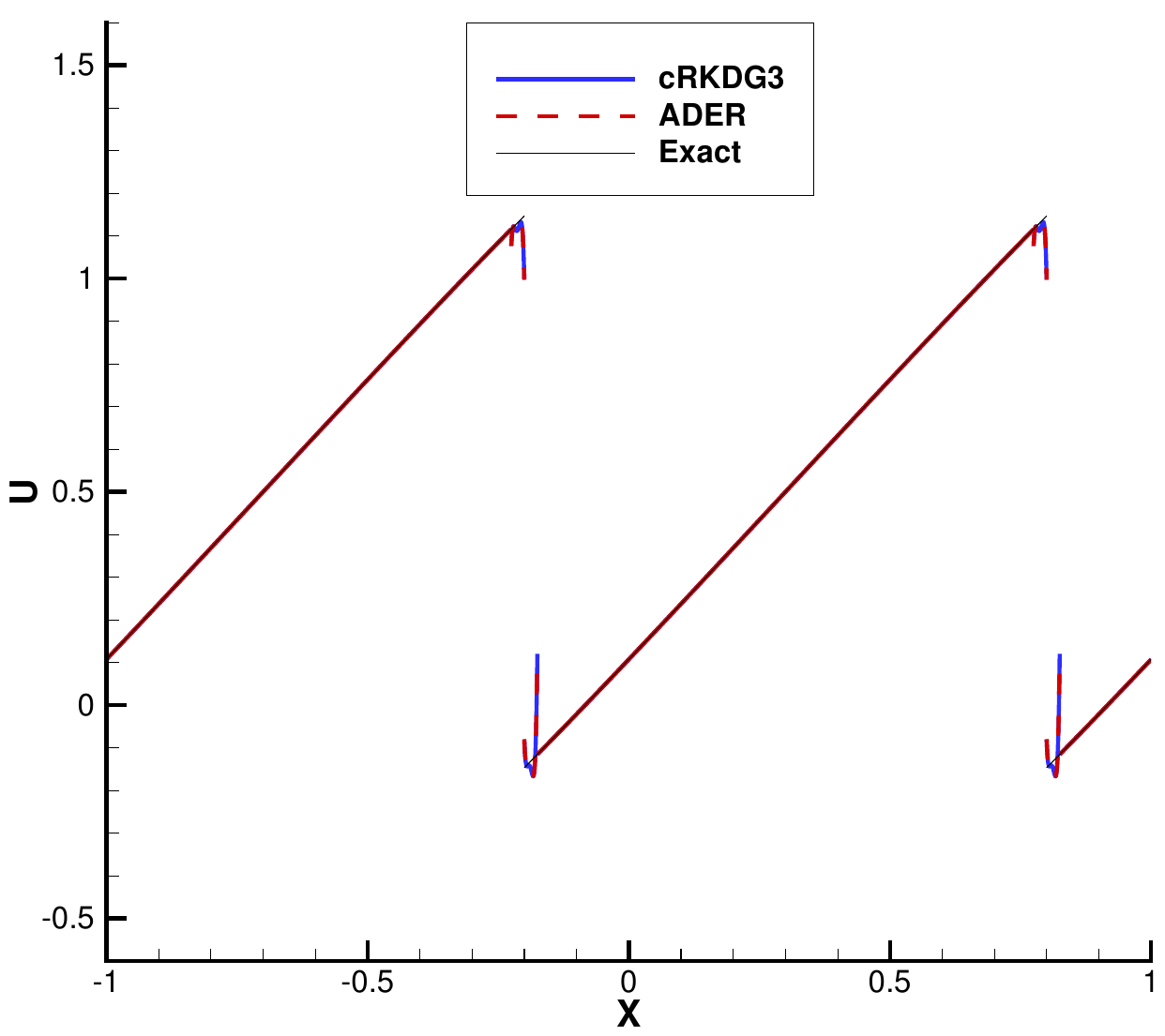}
		\caption{Diagonal slice, \(T=0.6\)}
	\end{subfigure}
	\caption{Two-dimensional Burgers, PACL without oscillation suppression, $N=80$.}
	\label{fig:num-burgers2d}
\end{figure}

\begin{example}[Two-dimensional Burgers shocks]\label{ex:num-burgers2d}
	We advance the Burgers initial data in \Cref{ex:num-burgers2d-acc} to $T=0.23$ and $0.6$ to examine the solution after shock formation \cite[Example~5.5]{10.1016/j.jcp.2009.12.030}. The mesh starts from an $N_x\times N_y$ rectangular grid with $N_x=N_y=80$, with the cells in half of the area split into right triangles. \Cref{fig:num-burgers2d} shows the solution field and a diagonal cut. These Burgers tests use no oscillation-suppressing limiter. The results align well with analytical solutions.
\end{example}

\begin{figure}[!ht]
	\centering
	\begin{subfigure}{0.3\linewidth}
		\centering
		\includegraphics[width=\linewidth]{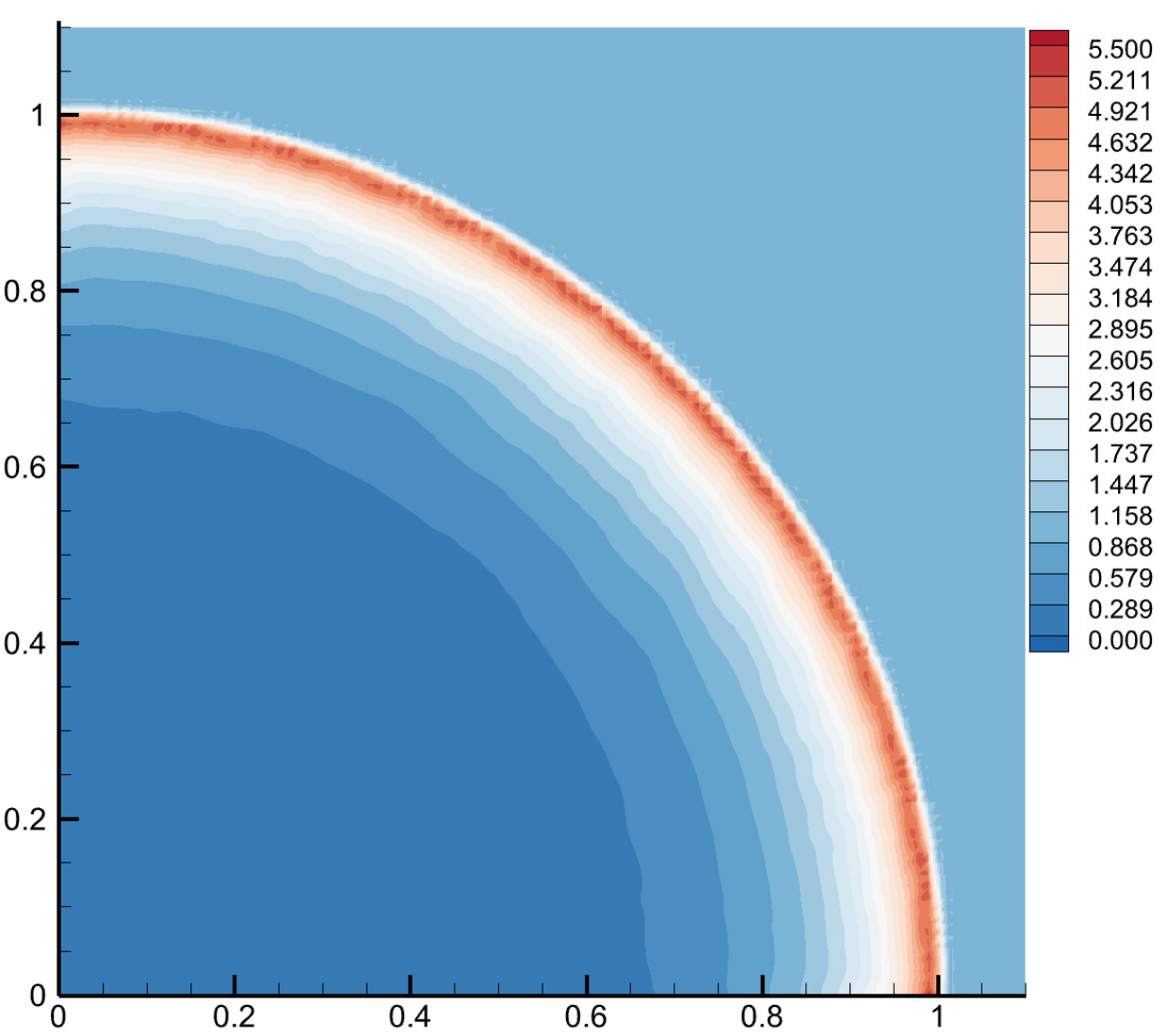}
		\caption{Density of cRKDG3}
	\end{subfigure}
	\begin{subfigure}{0.3\linewidth}
		\centering
		\includegraphics[width=\linewidth]{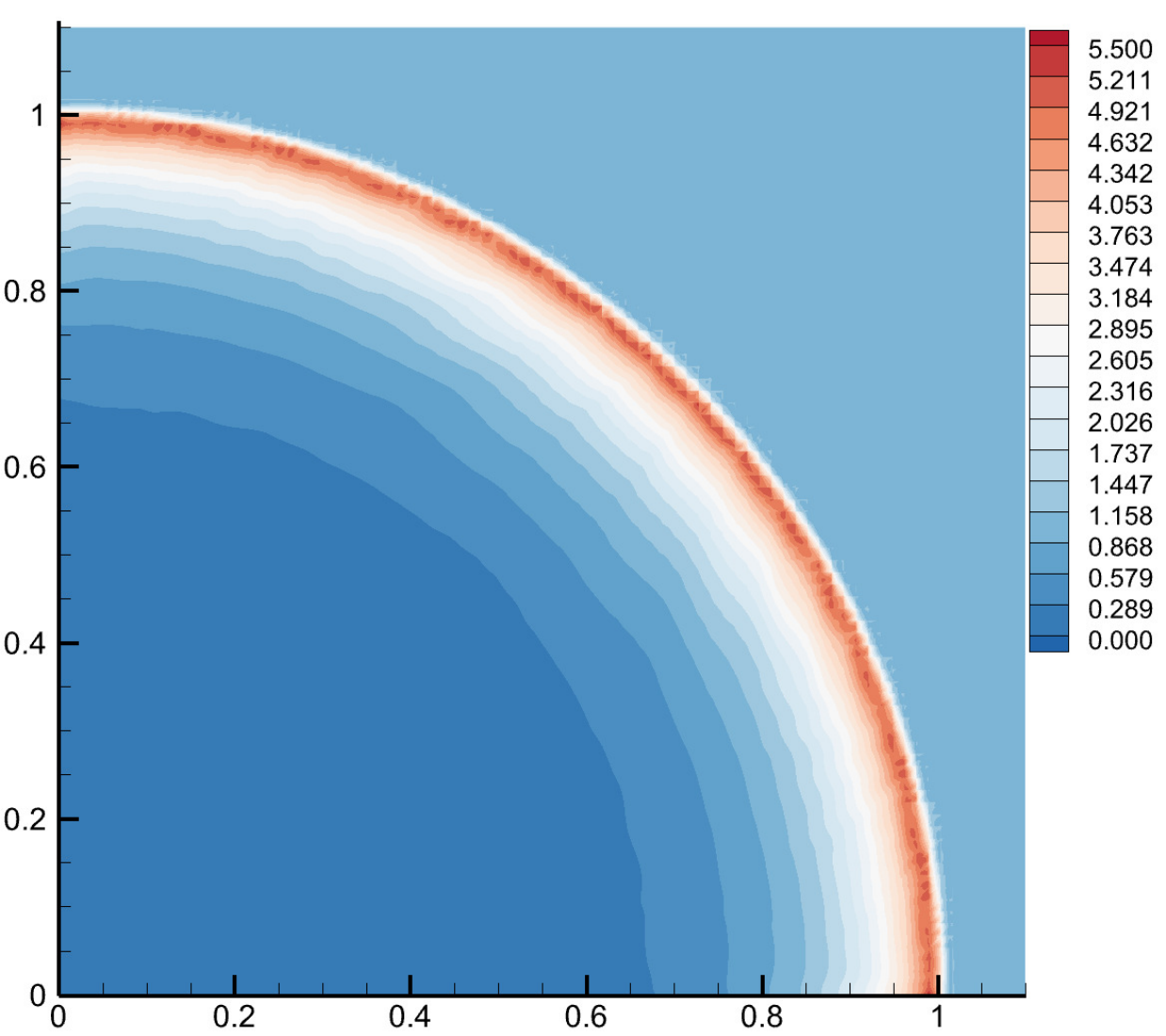}
		\caption{Density of ADER3}
	\end{subfigure}
	\begin{subfigure}{0.3\linewidth}
		\centering
		\includegraphics[width=\linewidth]{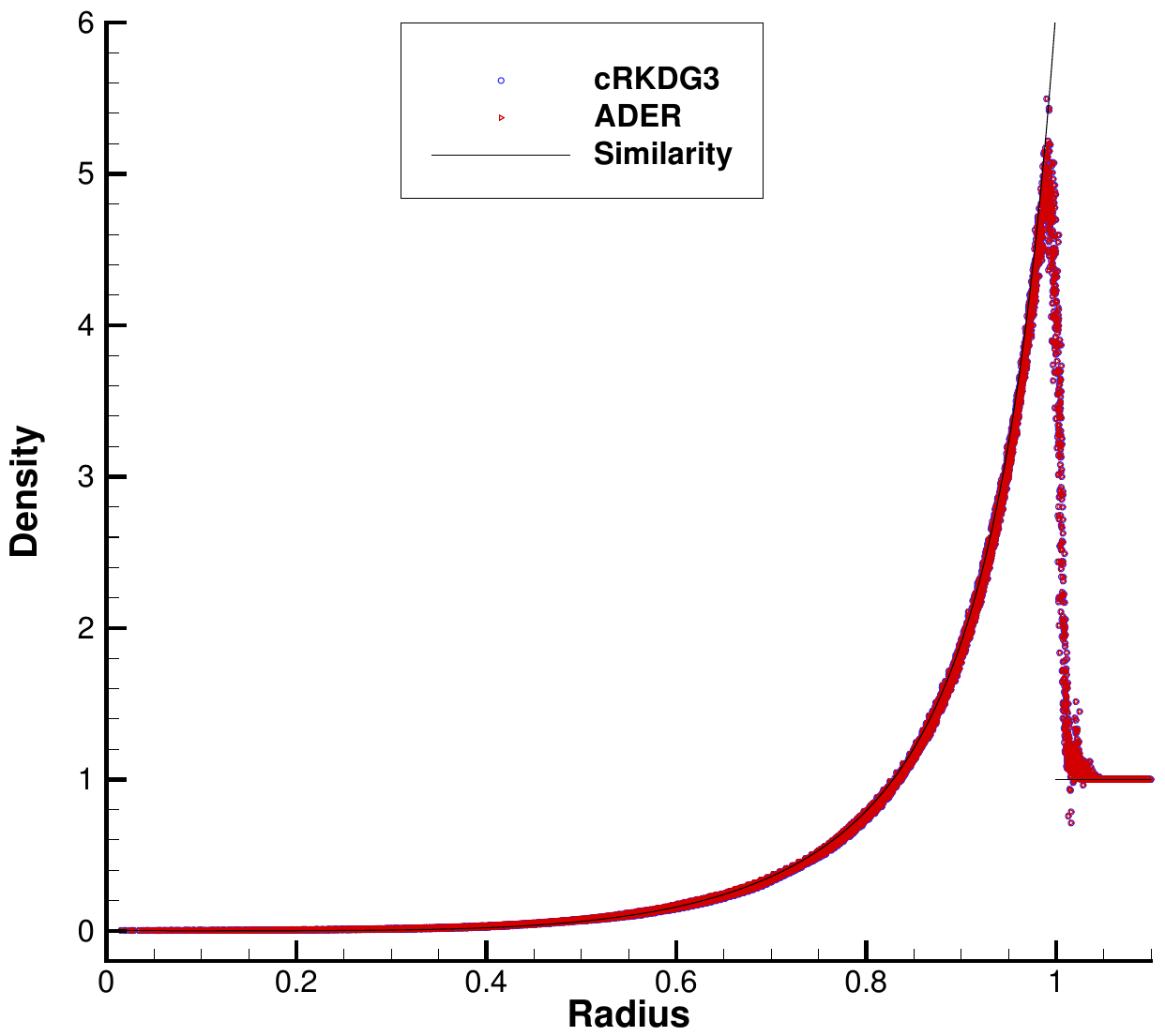}
		\caption{Radial plot of density}
	\end{subfigure}
	\caption{Two-dimensional Sedov problem, hybrid mesh of 9600 cells, $T=1$, PACL with RIOE.}
	\label{fig:num-sedov2d}
\end{figure}

\begin{example}[Two-dimensional Sedov blast]\label{ex:num-sedov2d}
	We simulate the two-dimensional Sedov blast on $[0, 1.1]^2$. The coarse mesh is an $N_x\times N_y$ rectangular grid with $N_x=N_y=20$, and the cells in half of its area are split into right triangles. The coarse mesh contains 400 triangles and 200 quadrilaterals. Two uniform refinements give 6400 triangles and 3200 quadrilaterals, for a total of 9600 cells. Initial density is one, velocity is zero and background energy density is $10^{-12}$ \cite{10.1016/j.jcp.2010.08.016}. Energy $0.244816$ is deposited uniformly in the lower-left coarse macrocell of side $1.1/20$. The left and bottom boundaries are reflective, while the other boundaries are outflow. Both refined computations reach $T=1$, and cRKDG3 and ADER3 agree quite well with each other. \Cref{fig:num-sedov2d} gives density plots and radial samples. Our numerical results show good radial pattern and robustness. 
\end{example}

\begin{figure}[!ht]
	\centering
	\begin{subfigure}[b]{0.4\linewidth}
		\centering
		\includegraphics[width=\linewidth]{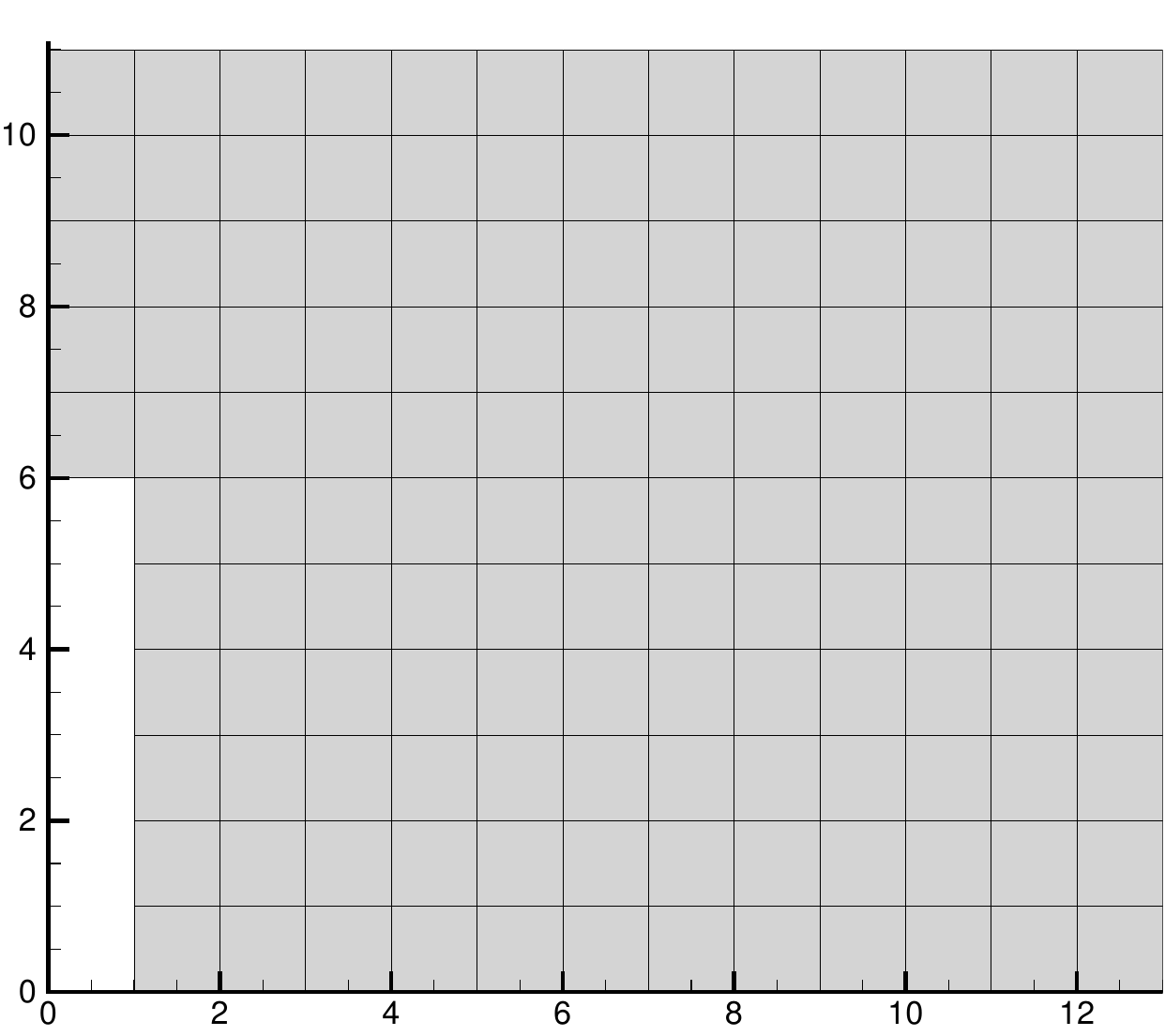}
		\caption{$90^\circ$ shock diffraction}
		\label{fig:mesh-diffraction90}
	\end{subfigure}
	\begin{subfigure}[b]{0.4\linewidth}
		\centering
		\includegraphics[width=\linewidth]{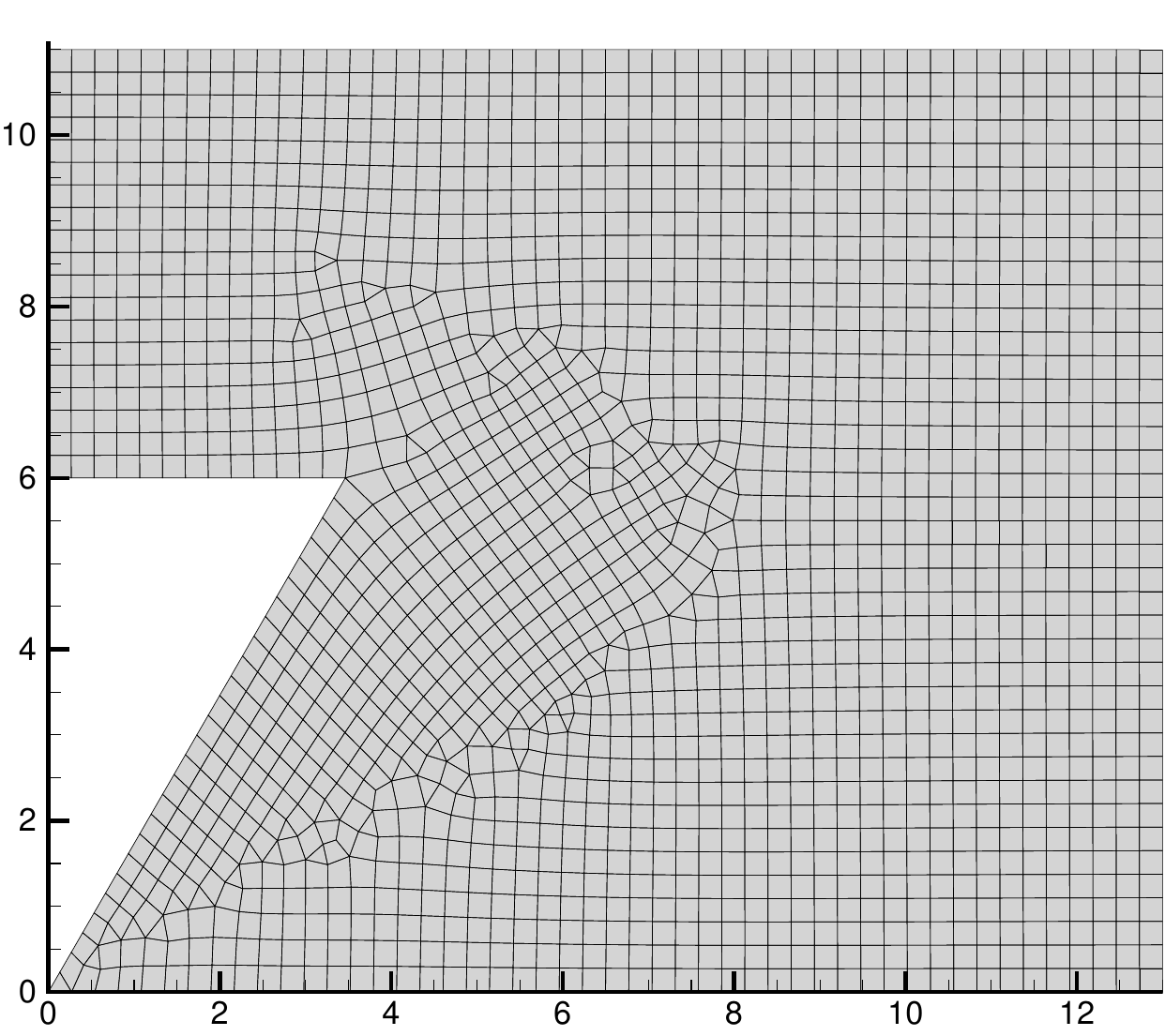}
		\caption{$120^\circ$ shock diffraction}
		\label{fig:mesh-diffraction120}
	\end{subfigure}
	\\
	\begin{subfigure}[b]{0.4\linewidth}
		\centering
		\includegraphics[width=\linewidth]{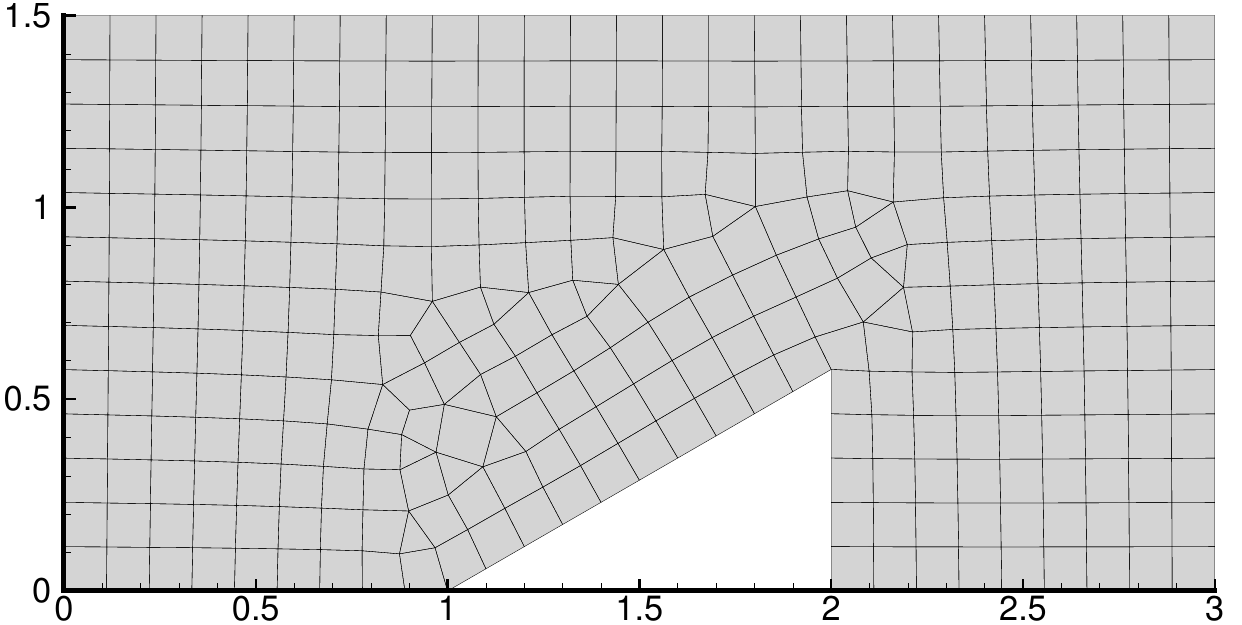}
		\caption{Schardin problem at Mach 1.34}
		\label{fig:mesh-schardin134}
	\end{subfigure}
	\begin{subfigure}[b]{0.4\linewidth}
		\centering
		\includegraphics[width=\linewidth]{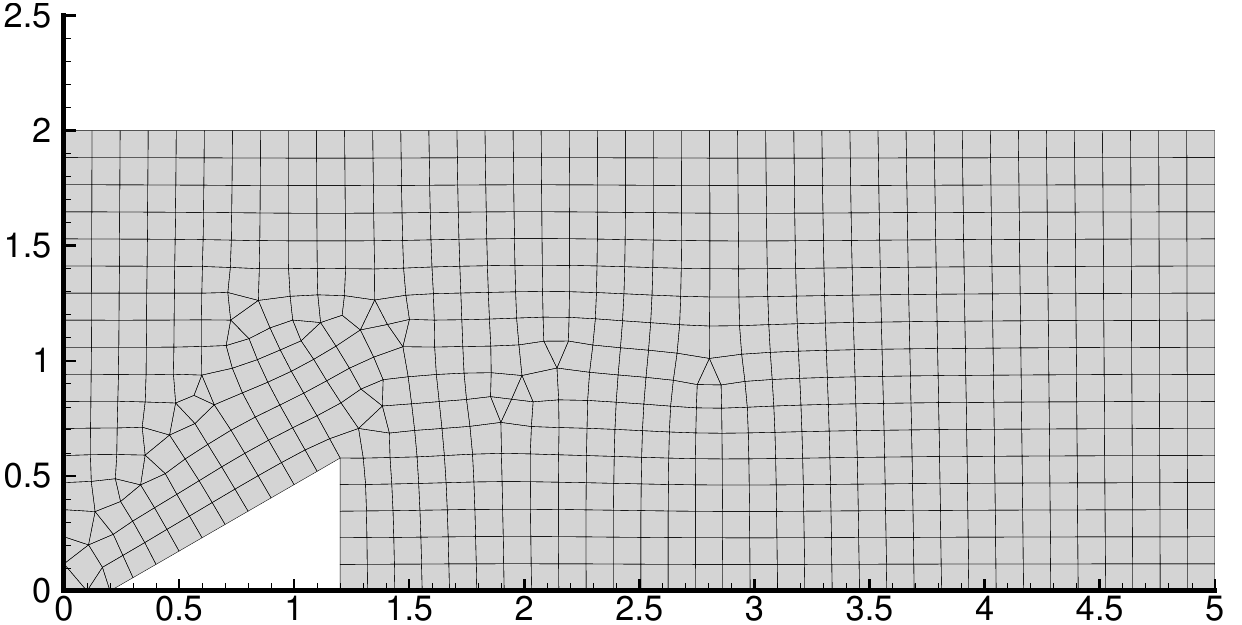}
		\caption{Schardin problem at Mach 10}
		\label{fig:mesh-schardin10}
	\end{subfigure}
	\caption{Coarse meshes for non-regular domains.}
\end{figure}

\begin{figure}[!ht]
	\centering
	\begin{subfigure}{0.4\linewidth}
		\centering
		\includegraphics[width=\linewidth]{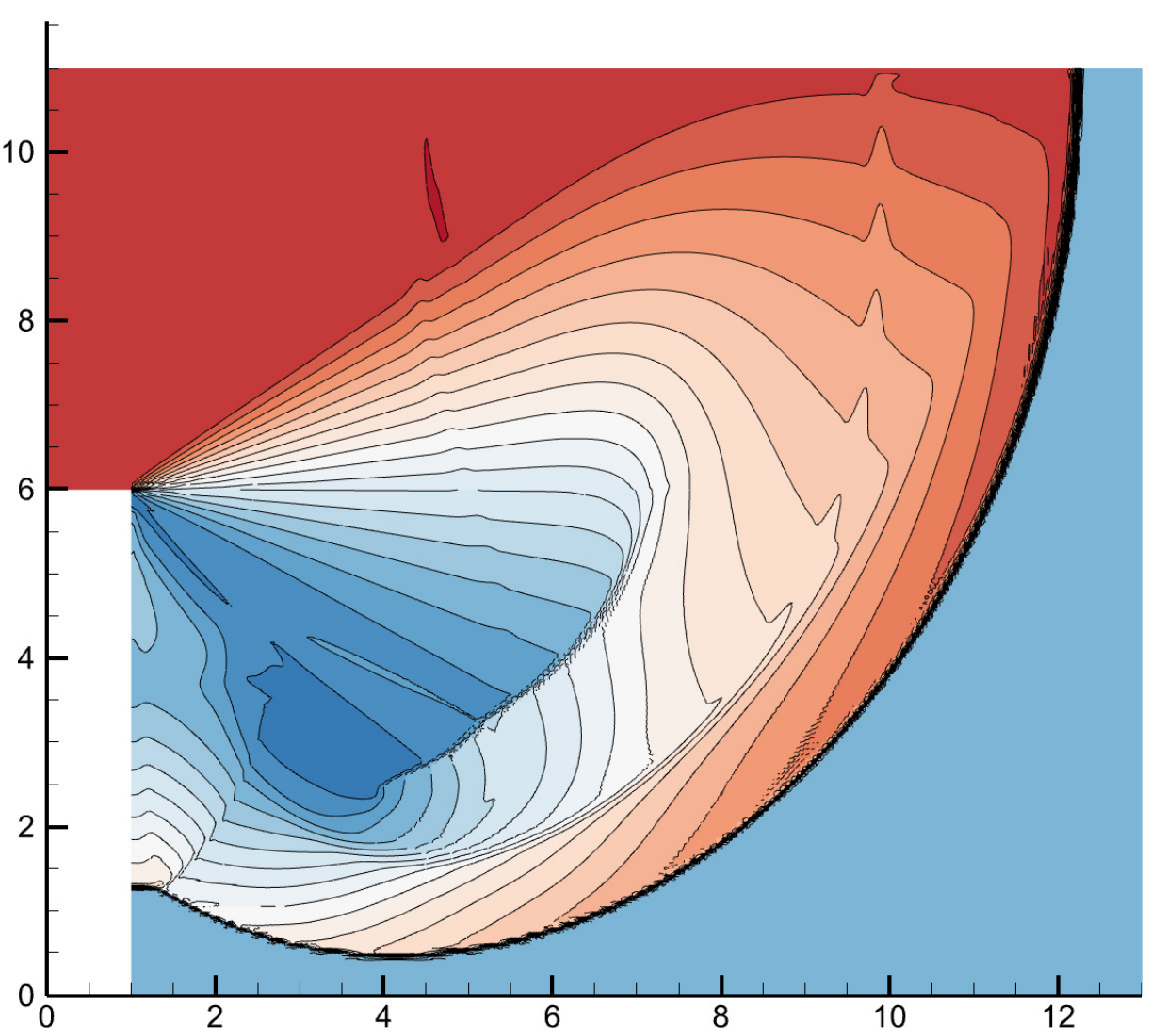}
		\caption{Density of cRKDG3: 20 equally spaced contours from 0.08 to 7.09}
	\end{subfigure}
	\begin{subfigure}{0.4\linewidth}
		\centering
		\includegraphics[width=\linewidth]{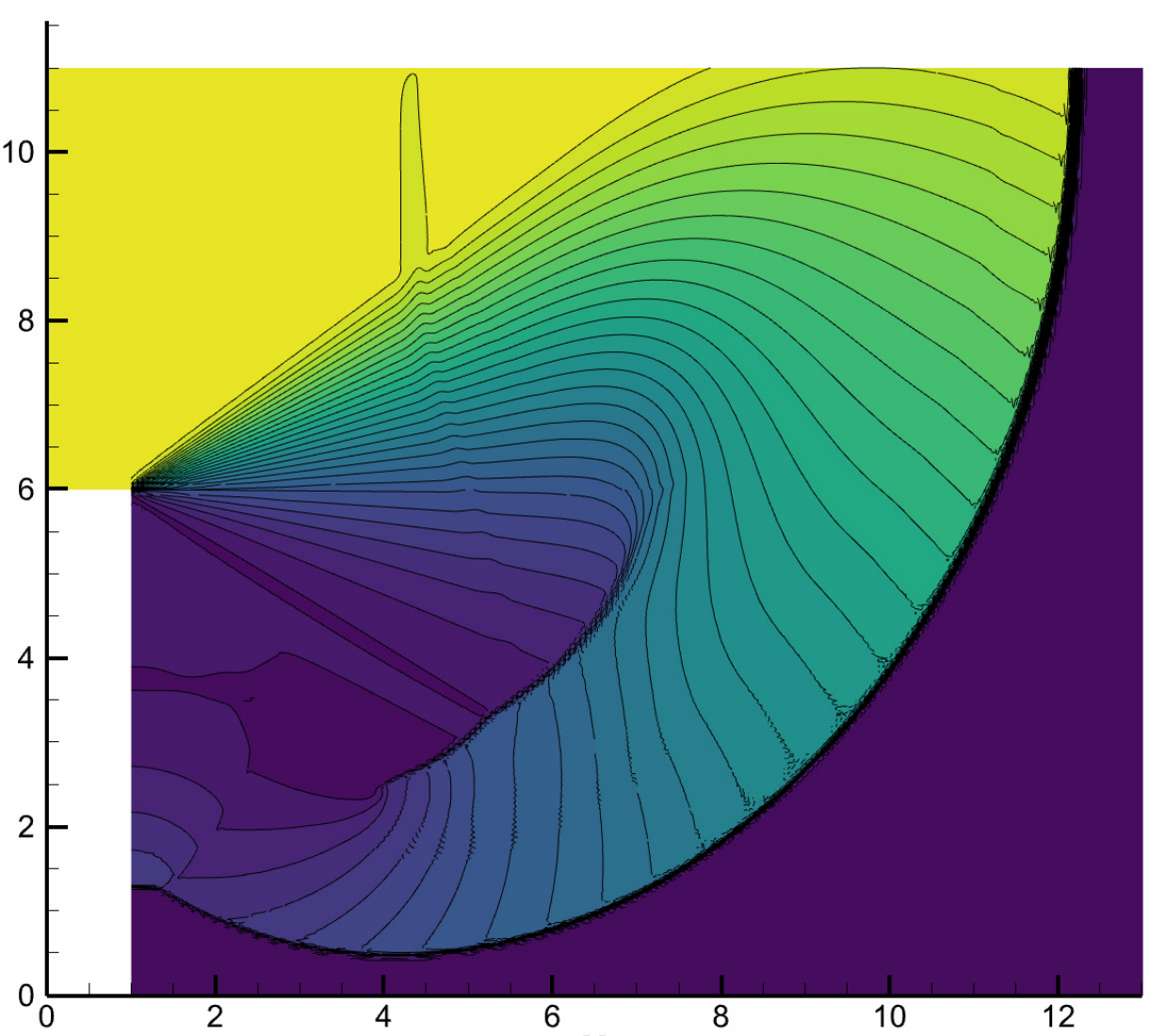}
		\caption{Pressure of cRKDG3: 30 equally spaced contours from 0.01 to 31}
	\end{subfigure}
	\\
	\begin{subfigure}{0.4\linewidth}
		\centering
		\includegraphics[width=\linewidth]{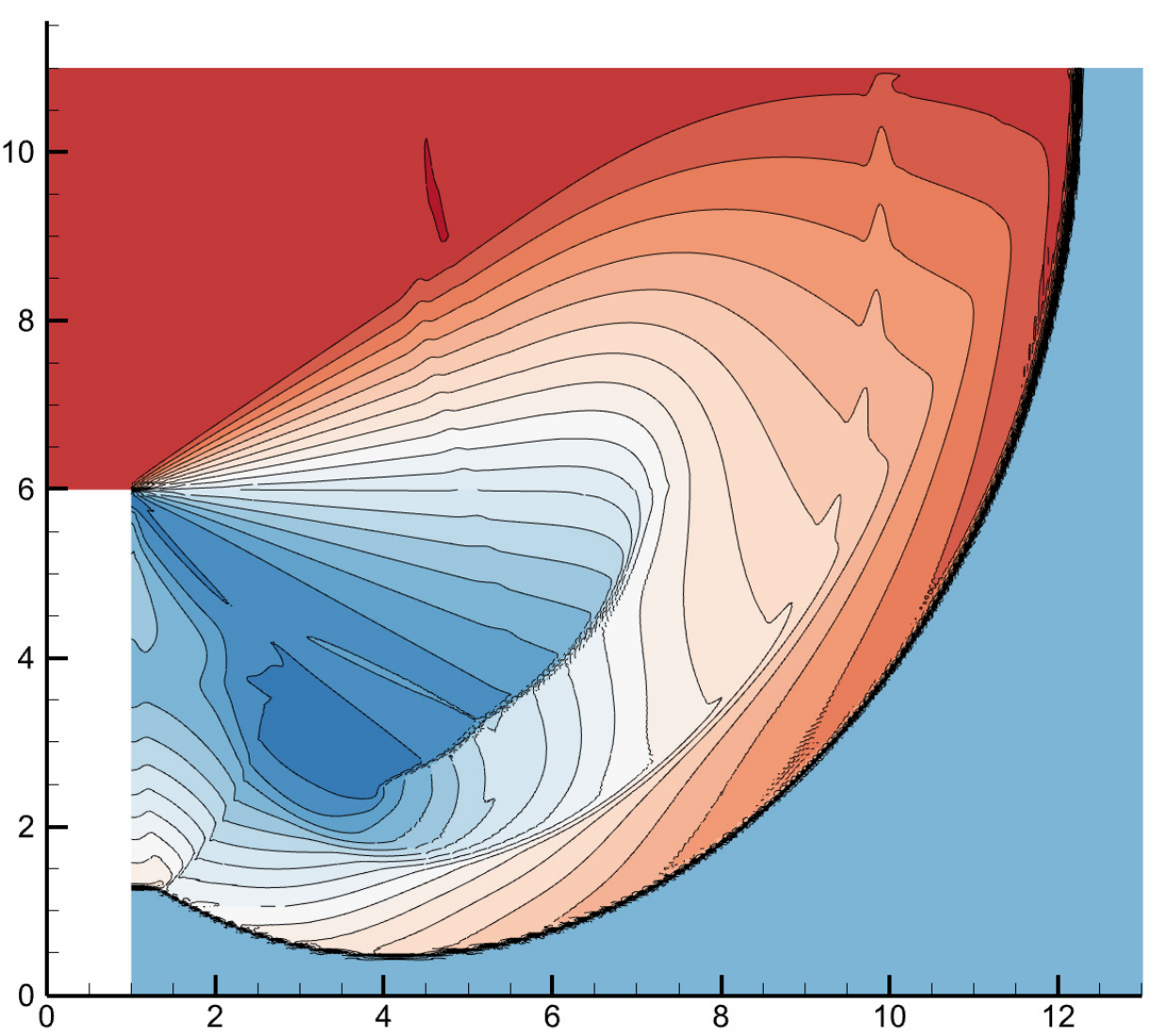}
		\caption{Density of ADER3: 20 equally spaced contours from 0.08 to 7.09}
	\end{subfigure}
	\begin{subfigure}{0.4\linewidth}
		\centering
		\includegraphics[width=\linewidth]{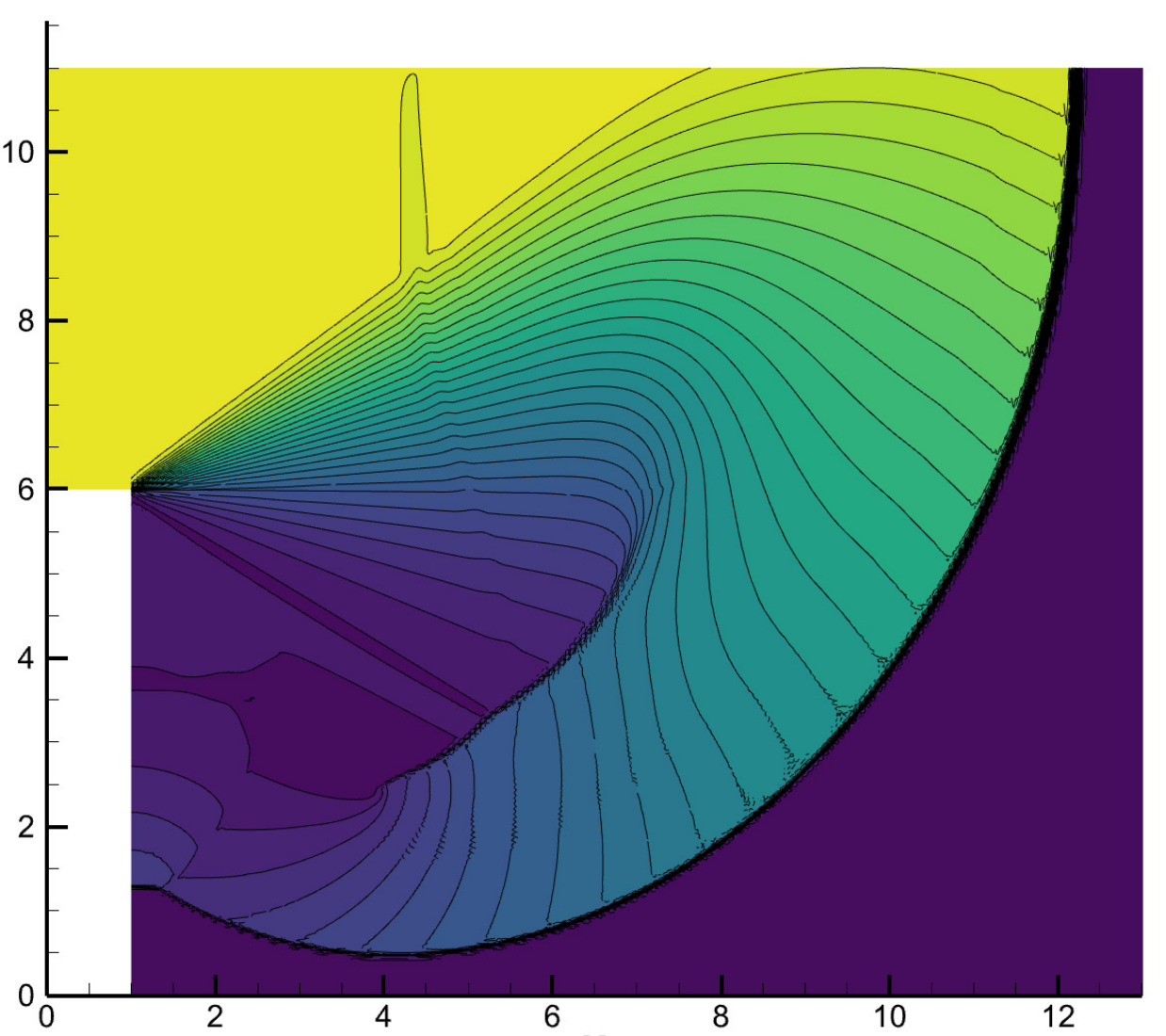}
		\caption{Pressure of ADER3: 30 equally spaced contours from 0.01 to 31}
	\end{subfigure}
	\caption{$90^\circ$ shock diffraction at $T=2.3$ on 35072 $\Qk^2$ cells, PACL with RIOE. Top: cRKDG3; bottom: ADER3. Left: density; right: pressure. The upper exterior state follows the exact moving shock.}
	\label{fig:num-diffraction90}
\end{figure}

\begin{example}[$90^\circ$ shock diffraction]\label{ex:num-diffraction90}
	A Mach-5.09 shock diffracts around a $90^\circ$ corner. Initially at $x=0.5$, the shock enters undisturbed gas with $(\rho,u,v,p)=(1.4,0,0,1)$ \cite[Section~4.5]{10.1016/j.jcp.2010.08.016}. The domain is $([0,1]\times[6,11])\cup([1,13]\times[0,11])$. The two corner walls are free-slip reflective, the left boundary prescribes the post-shock state, and the right and bottom boundaries are outflow. At the upper boundary, following \cite[Example 4.4]{10.1016/j.jcp.2022.111424}, the exterior state follows the exact moving shock $x_s(t) = 0.5 + 5.09t$, evaluated at each flux time node. Four uniform refinements of the coarse mesh (\Cref{fig:mesh-diffraction90}) with 137 unit square cells give 35072 $\Qk^2$ cells with side length $1/16$. \Cref{fig:num-diffraction90} shows density and pressure. This example is nearly vacuum near the diffraction corner and thus poses a challenge for pressure positivity.
\end{example}

\begin{figure}[!ht]
	\centering
	\begin{subfigure}{0.4\linewidth}
		\centering
		\includegraphics[width=\linewidth]{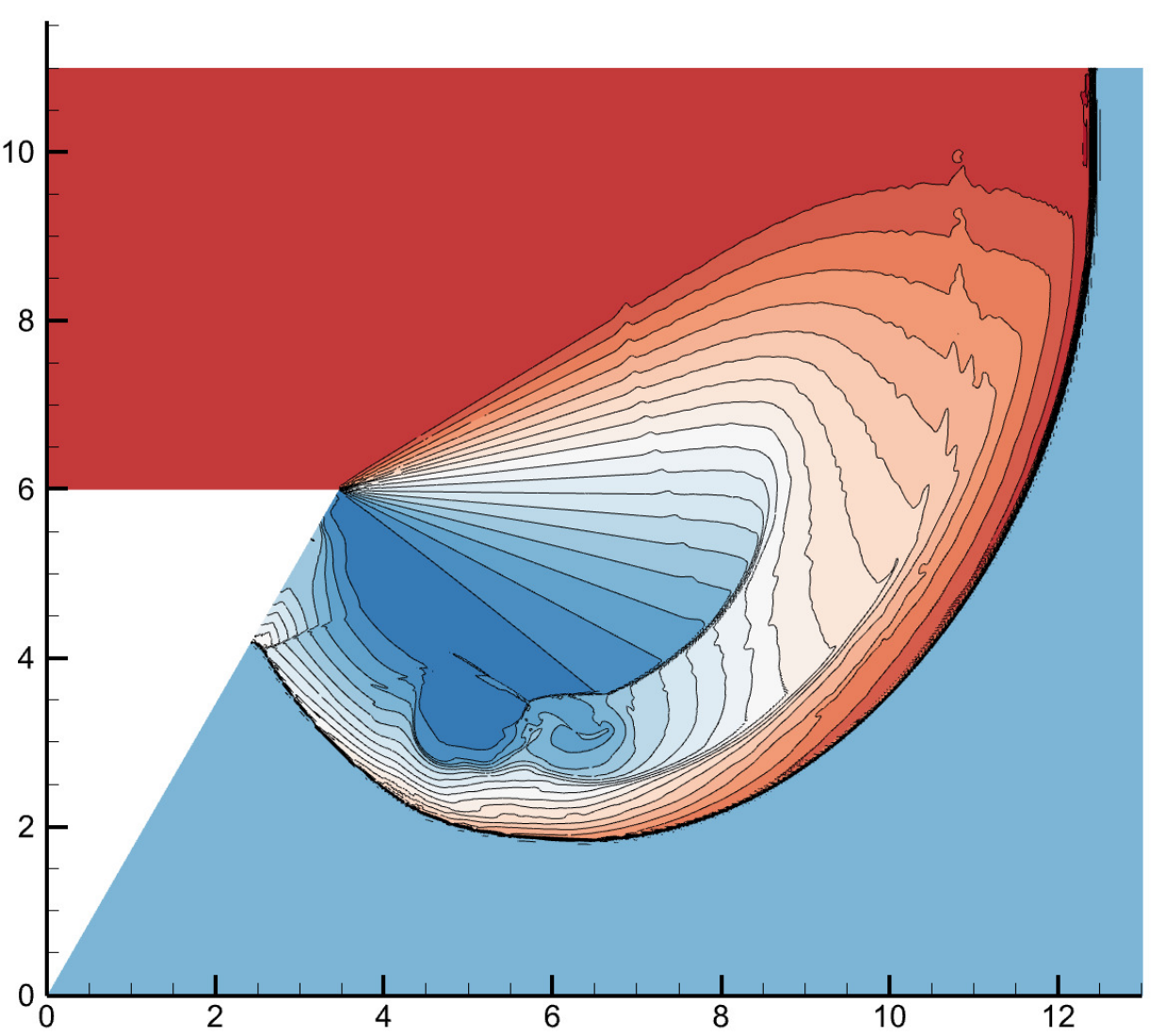}
		\caption{Density of cRKDG3: 20 equally spaced contours from 0.04 to 8.1}
	\end{subfigure}
	\begin{subfigure}{0.4\linewidth}
		\centering
		\includegraphics[width=\linewidth]{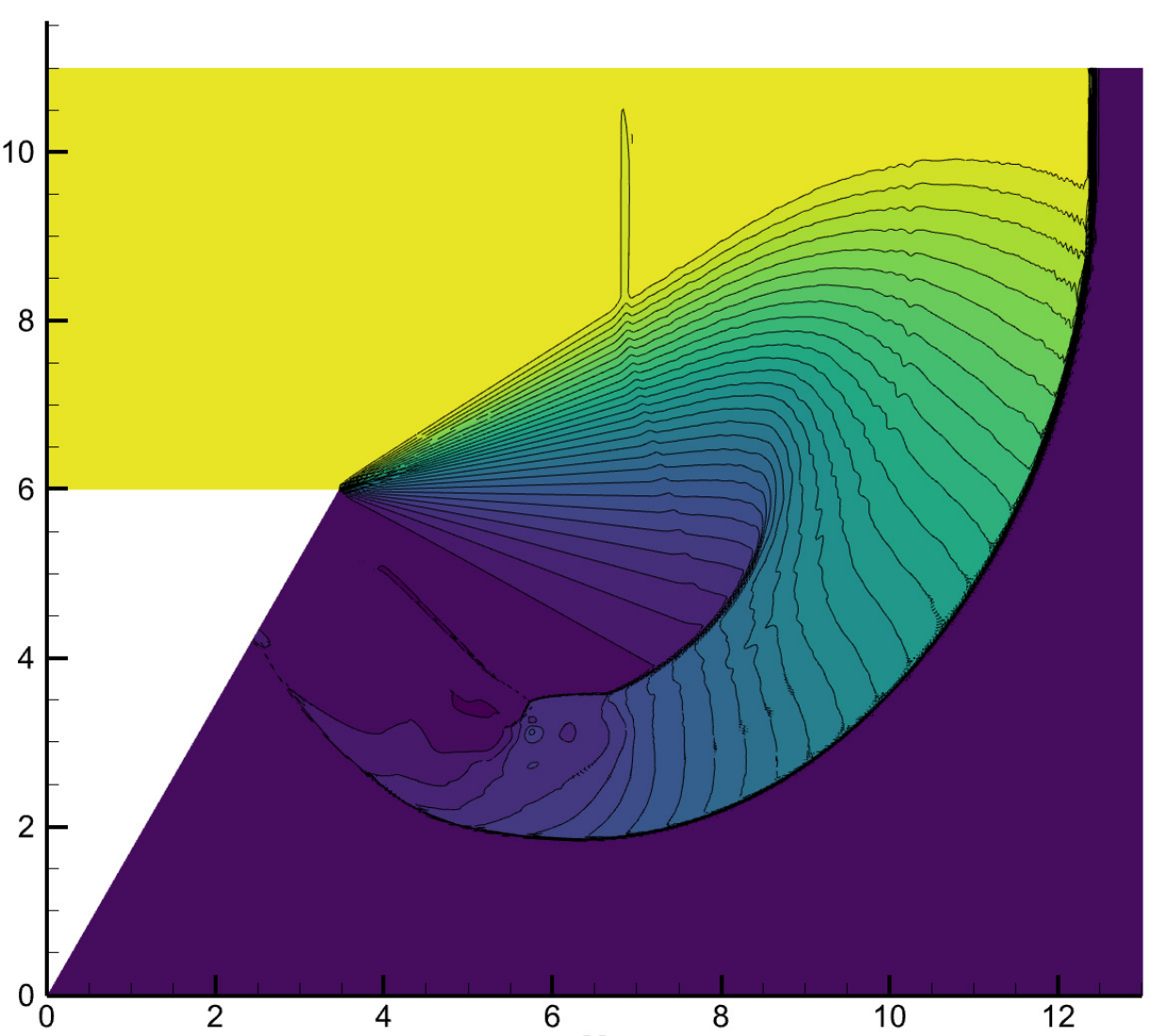}
		\caption{Pressure of cRKDG3: 30 equally spaced contours from 0.01 to 31}
	\end{subfigure}
	\\
	\begin{subfigure}{0.4\linewidth}
		\centering
		\includegraphics[width=\linewidth]{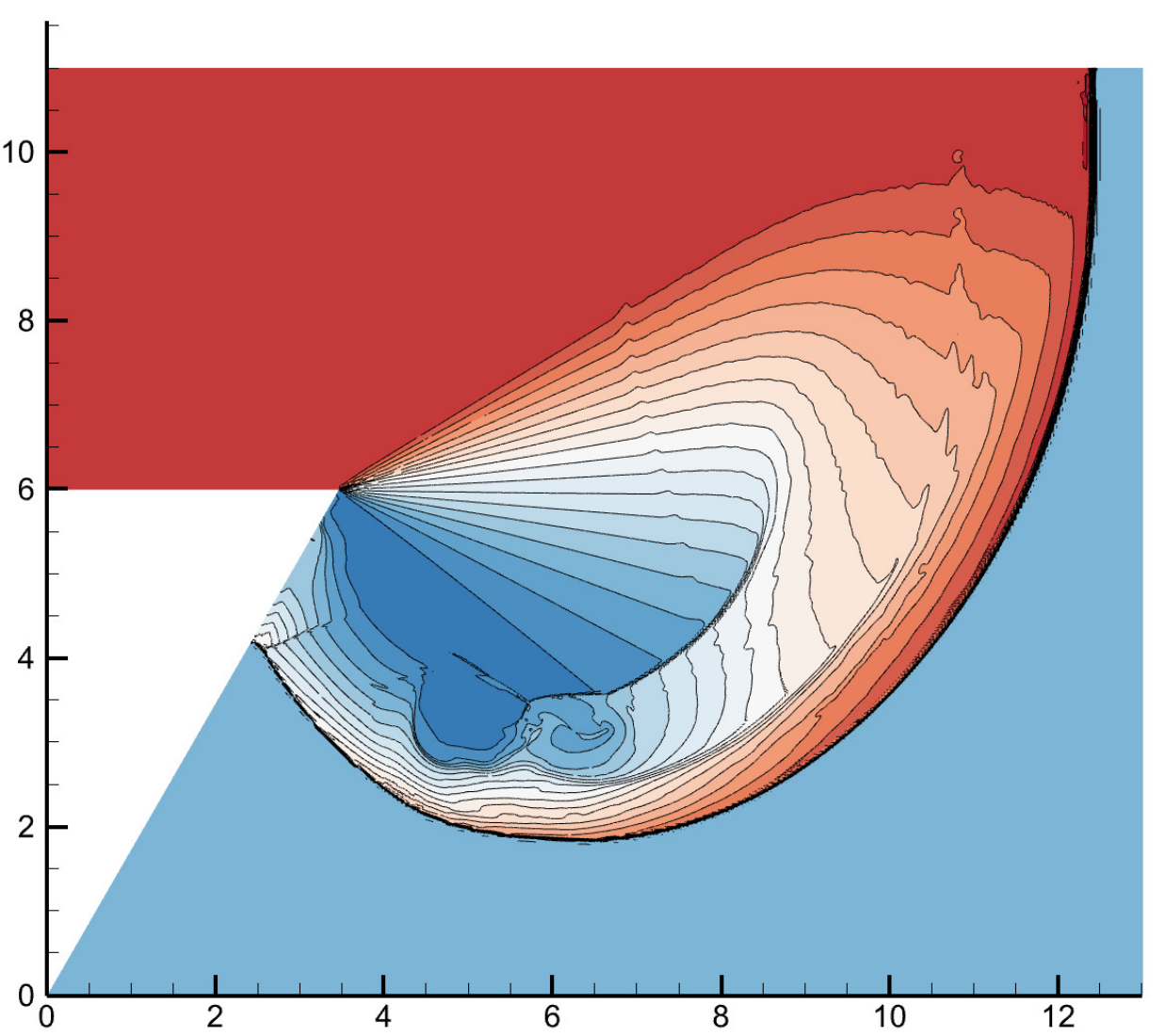}
		\caption{Density of ADER3: 20 equally spaced contours from 0.04 to 8.1}
	\end{subfigure}
	\begin{subfigure}{0.4\linewidth}
		\centering
		\includegraphics[width=\linewidth]{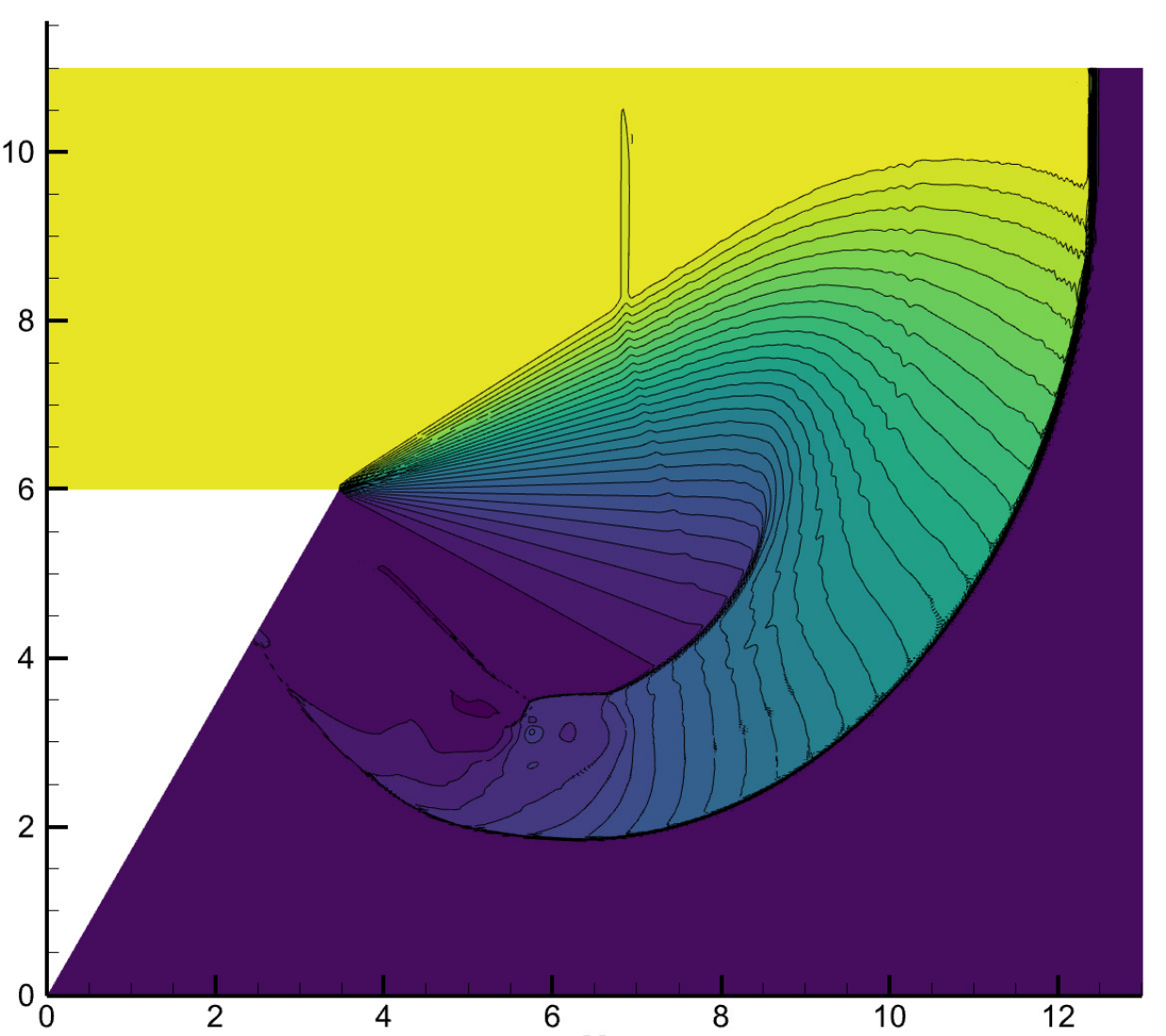}
		\caption{Pressure of ADER3: 30 equally spaced contours from 0.01 to 31}
	\end{subfigure}
	\caption{$120^\circ$ shock diffraction at $T=0.9$ on 119744 hybrid cells, PACL with RIOE. Top: cRKDG3; bottom: ADER3. Left: density; right: pressure.}
	\label{fig:num-diffraction120}
\end{figure}

\begin{example}[$120^\circ$ shock diffraction]\label{ex:num-diffraction120}
	We consider a Mach-10 shock diffracting around a $120^\circ$ corner. The initial shock is at $x=3.4$, $6\leq y\leq 11$, with ambient state $(\rho,v_1,v_2,p)=(1.4,0,0,1)$ \cite[Example~6.9]{10.1007/s10915-011-9472-8}. The straight ramp joins $(0,0)$ to $(2\sqrt3,6)$. The horizontal wall is $y=6$ to the left of the corner. The fluid fills the outer rectangle $[0,13]\times[0,11]$ above $y=6$ and lies to the right of the ramp below it. Both walls reflect. Following \cite[Example 6.2.6]{10.1016/j.jcp.2017.05.025}, the upper boundary prescribes the exact shock at $x_s(t)=3.4+10t$: $(8,8.25,0,116.5)$ behind the shock and $(1.4,0,0,1)$ ahead. The left boundary prescribes post-shock inflow, and the right and bottom boundaries are outflow. Three uniform refinements of the 1871-cell hybrid coarse grid (\Cref{fig:mesh-diffraction120}) give 119744 cells. \Cref{fig:num-diffraction120} shows the near-vertical upper incident shock and the curved diffracted front below it. Our result is comparable to those in the literature. 
\end{example}

\begin{figure}[!ht]
	\centering
	\begin{subfigure}{0.3\linewidth}
		\centering
		\includegraphics[width=\linewidth]{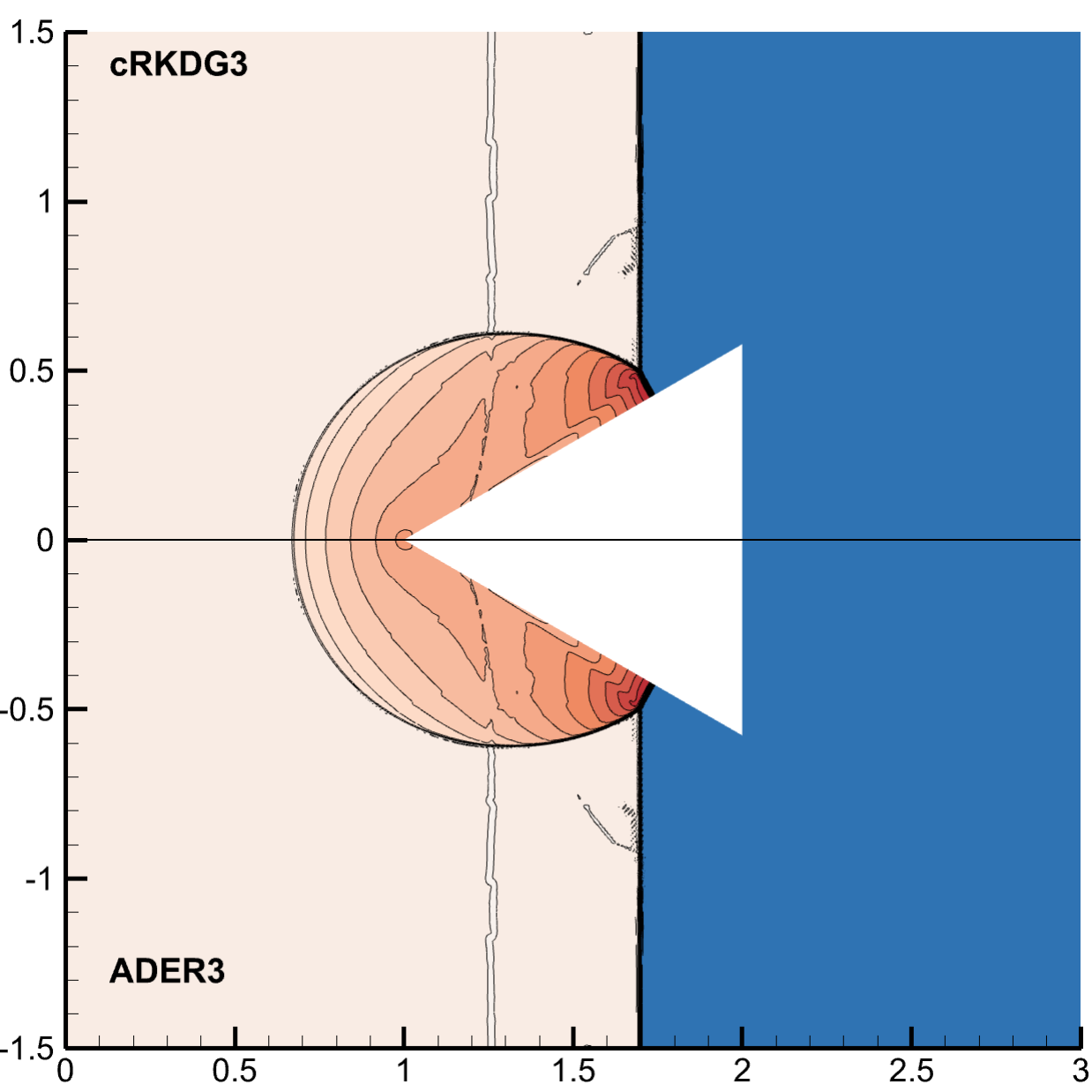}
		\caption{\(T = 0.521\), 30 equally spaced contours from 1.39 to 2.88}
	\end{subfigure}
	\begin{subfigure}{0.3\linewidth}
		\centering
		\includegraphics[width=\linewidth]{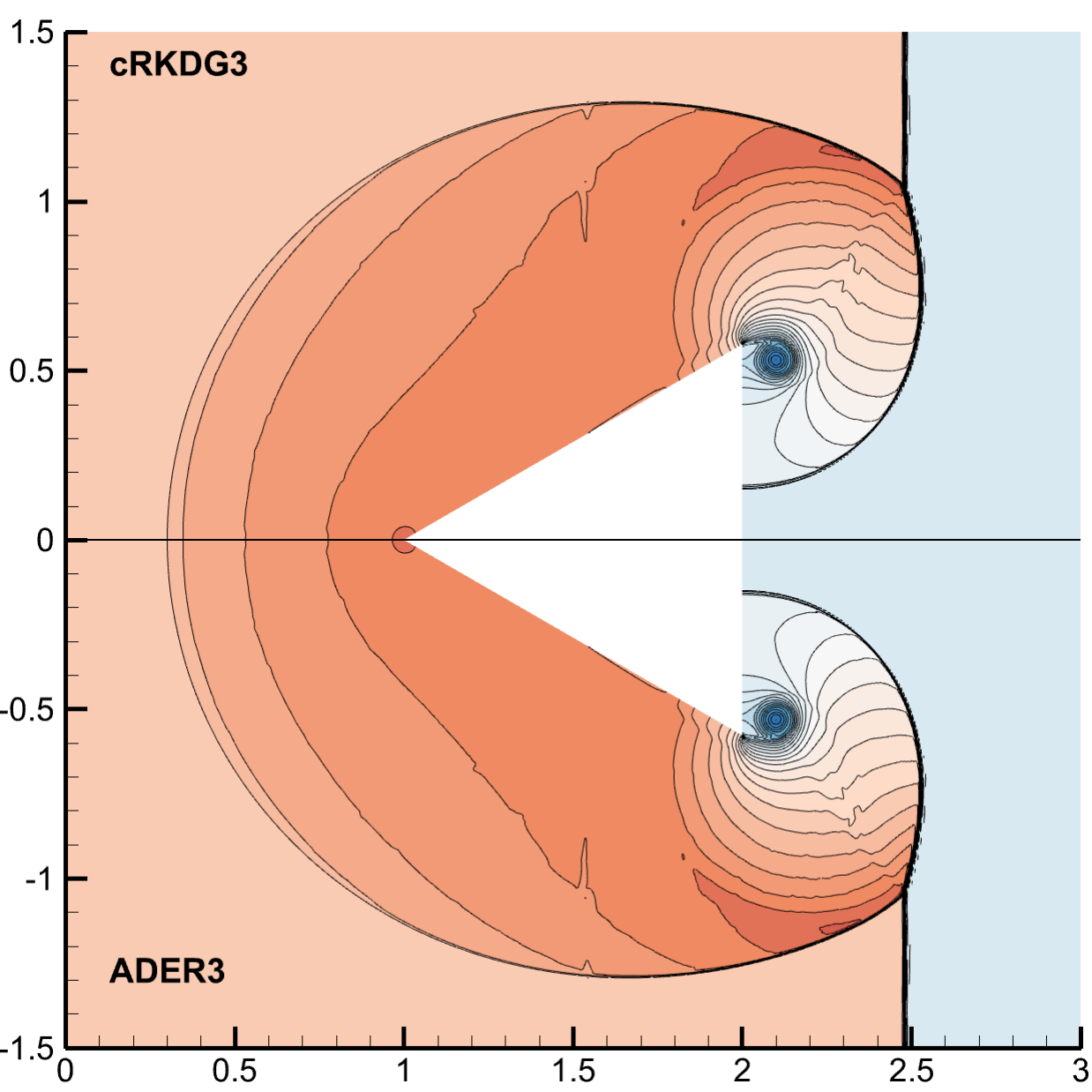}
		\caption{\(T = 1.104\), 30 equally spaced contours from 0.6 to 2.9}
	\end{subfigure}
	\begin{subfigure}{0.3\linewidth}
		\centering
		\includegraphics[width=\linewidth]{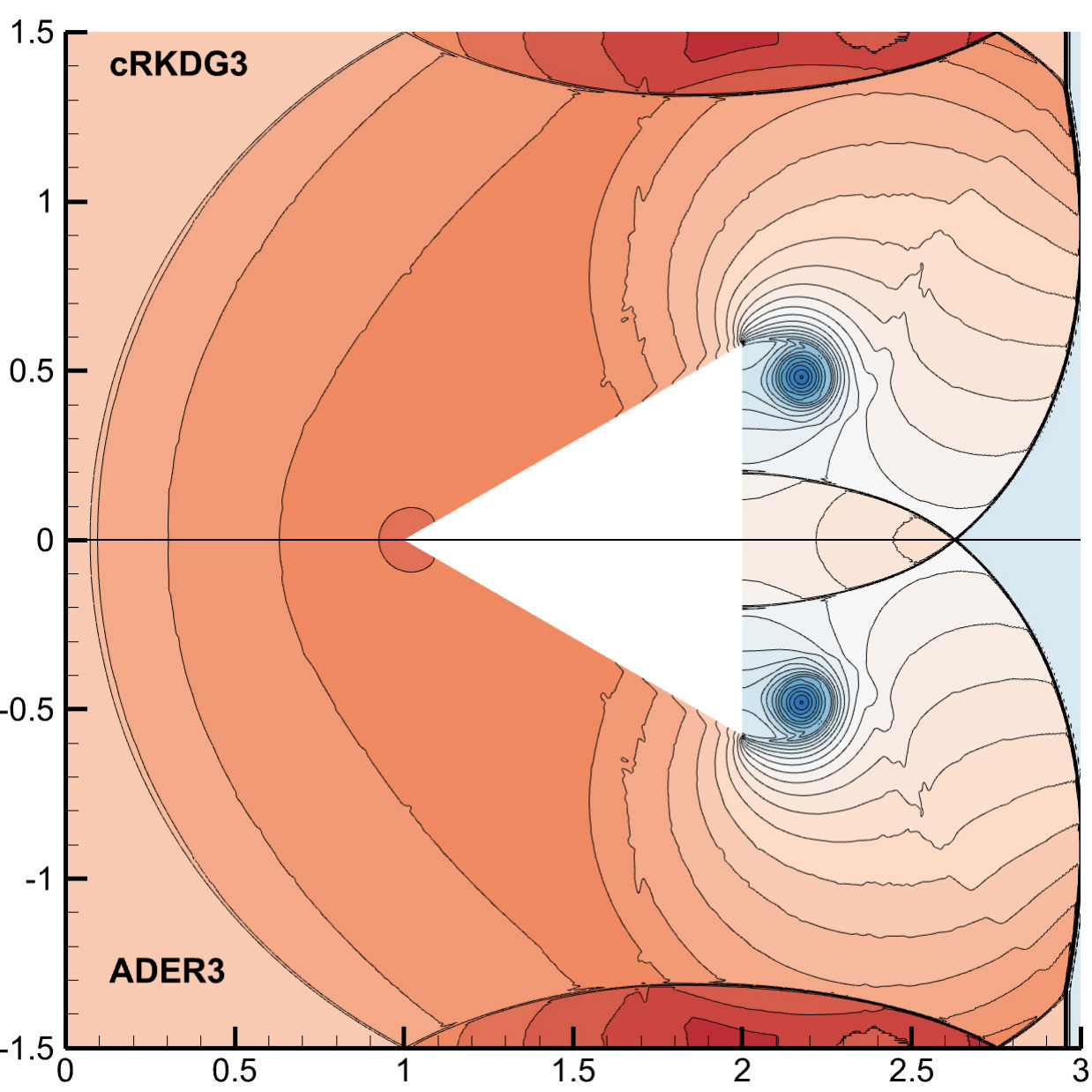}
		\caption{\(T = 1.46\), 30 equally spaced contours from 0.6 to 2.9}
	\end{subfigure}
	\caption{Schardin Mach 1.34, PACL with RIOE. The upper half of each plot is cRKDG3 and the lower half is ADER3.}
	\label{fig:num-schardin134}
\end{figure}

\begin{example}[Schardin problem at Mach 1.34]\label{ex:num-schardin134}
	The Schardin problem describes a planar shock passing an equilateral wedge. At Mach~1.34, we compute the upper half-domain $[0,3] \times[0,1.5]$ excluding the half-wedge with vertices $(1,0),(2,1/\sqrt3),(2,0)$. The initial shock is at $x=1$. The ambient state is $(\rho,v_1,v_2,p)=(1.4,0,0,1)$, and the post-shock state follows the normal-shock Rankine--Hugoniot relations \cite[Example~6.10]{10.1007/s10915-011-9472-8}. The left boundary prescribes post-shock inflow, the right boundary is outflow, and the upper boundary, wedge, and symmetry axis are free-slip reflective. The hybrid coarse mesh (\Cref{fig:mesh-schardin134}) has 321 cells. Four uniform refinements give 82176 cells. \Cref{fig:num-schardin134} shows density. Our numerical method preserves positive density and pressure near diffraction corners.
\end{example}

\begin{figure}[!ht]
	\centering
	\begin{subfigure}{0.3\linewidth}
		\centering
		\includegraphics[width=\linewidth]{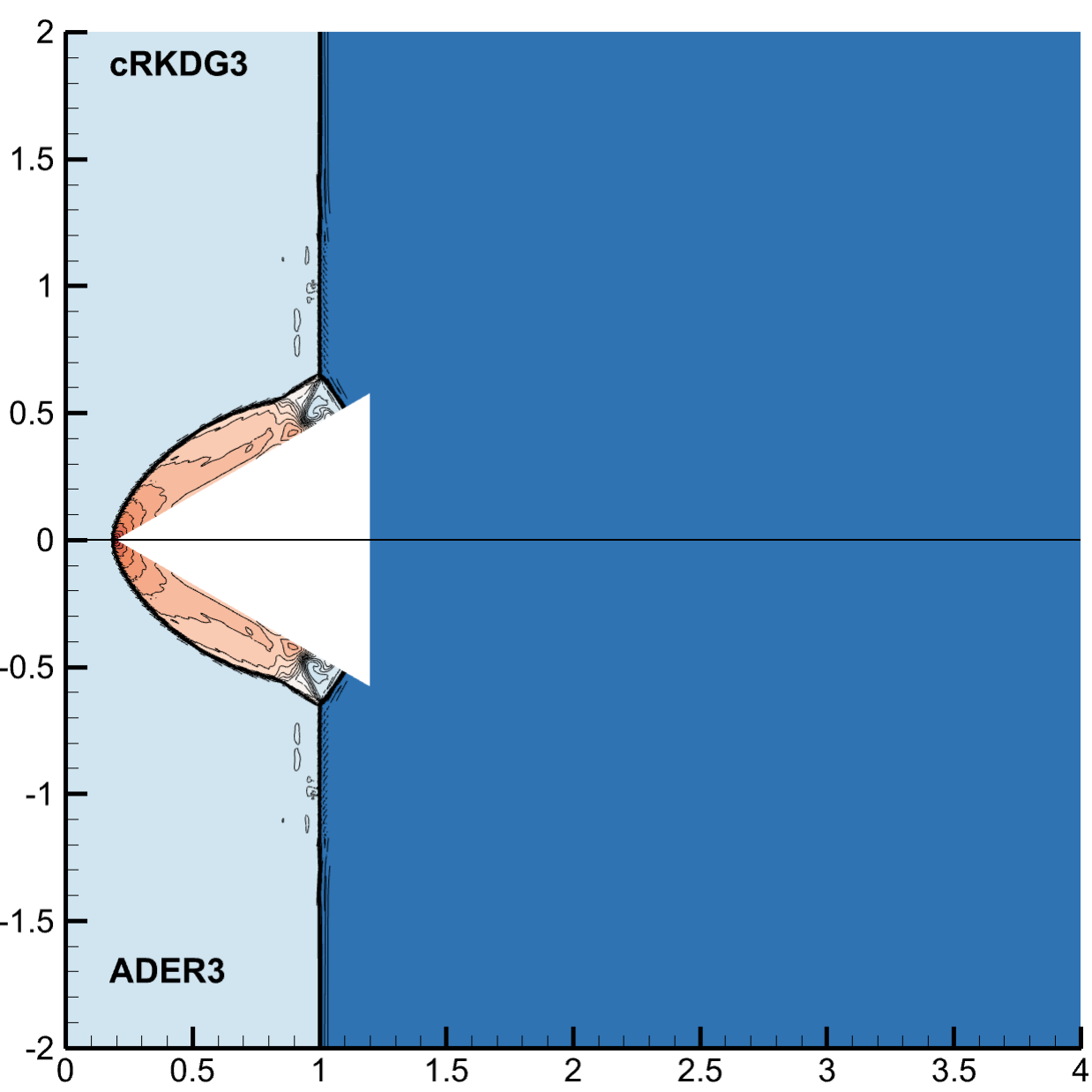}
		\caption{\(T = 0.08\), 30 equally spaced contours from 0.078 to 22}
	\end{subfigure}
	\begin{subfigure}{0.3\linewidth}
		\centering
		\includegraphics[width=\linewidth]{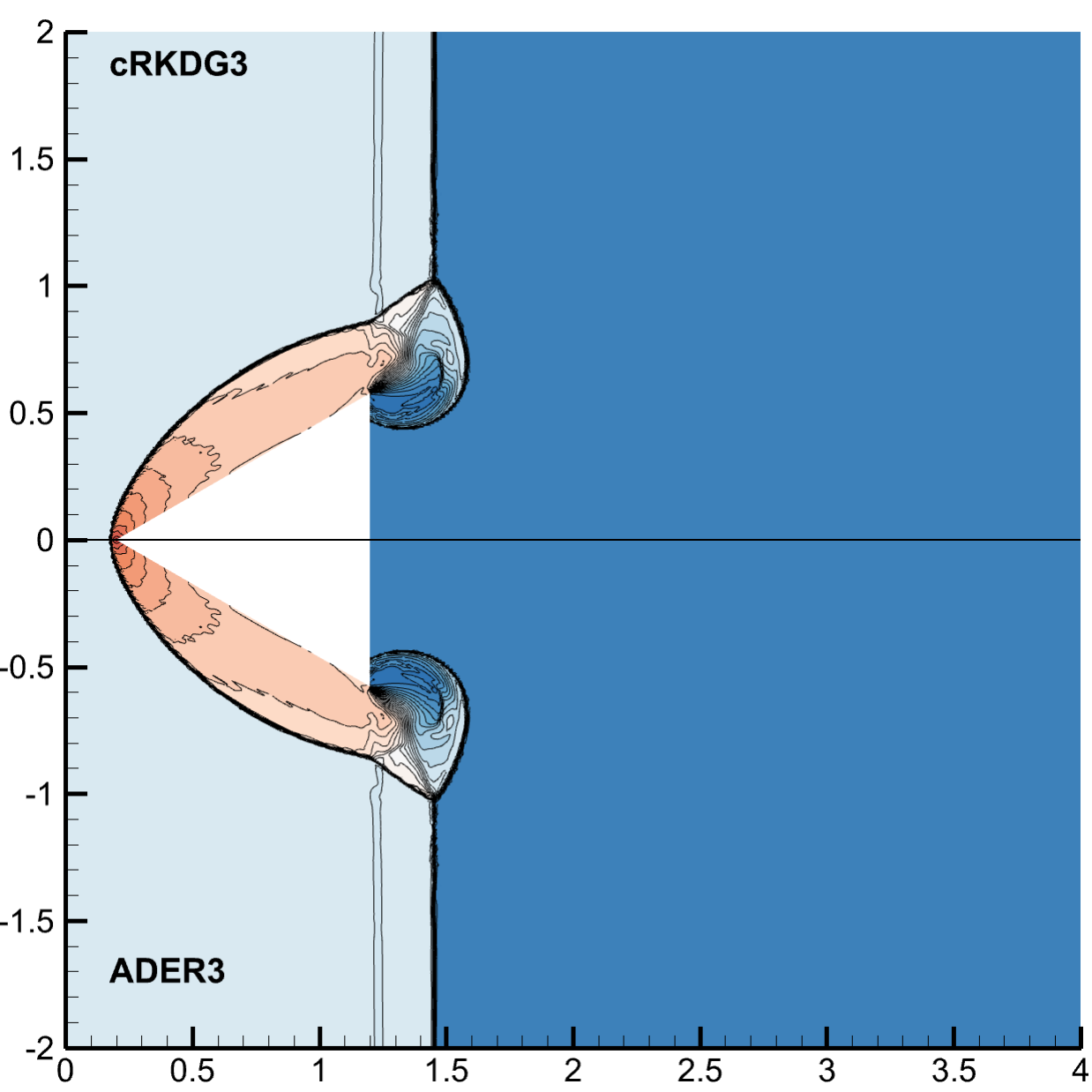}
		\caption{\(T = 0.125\), 30 equally spaced contours from 0.6 to 2.9}
	\end{subfigure}
	\begin{subfigure}{0.3\linewidth}
		\centering
		\includegraphics[width=\linewidth]{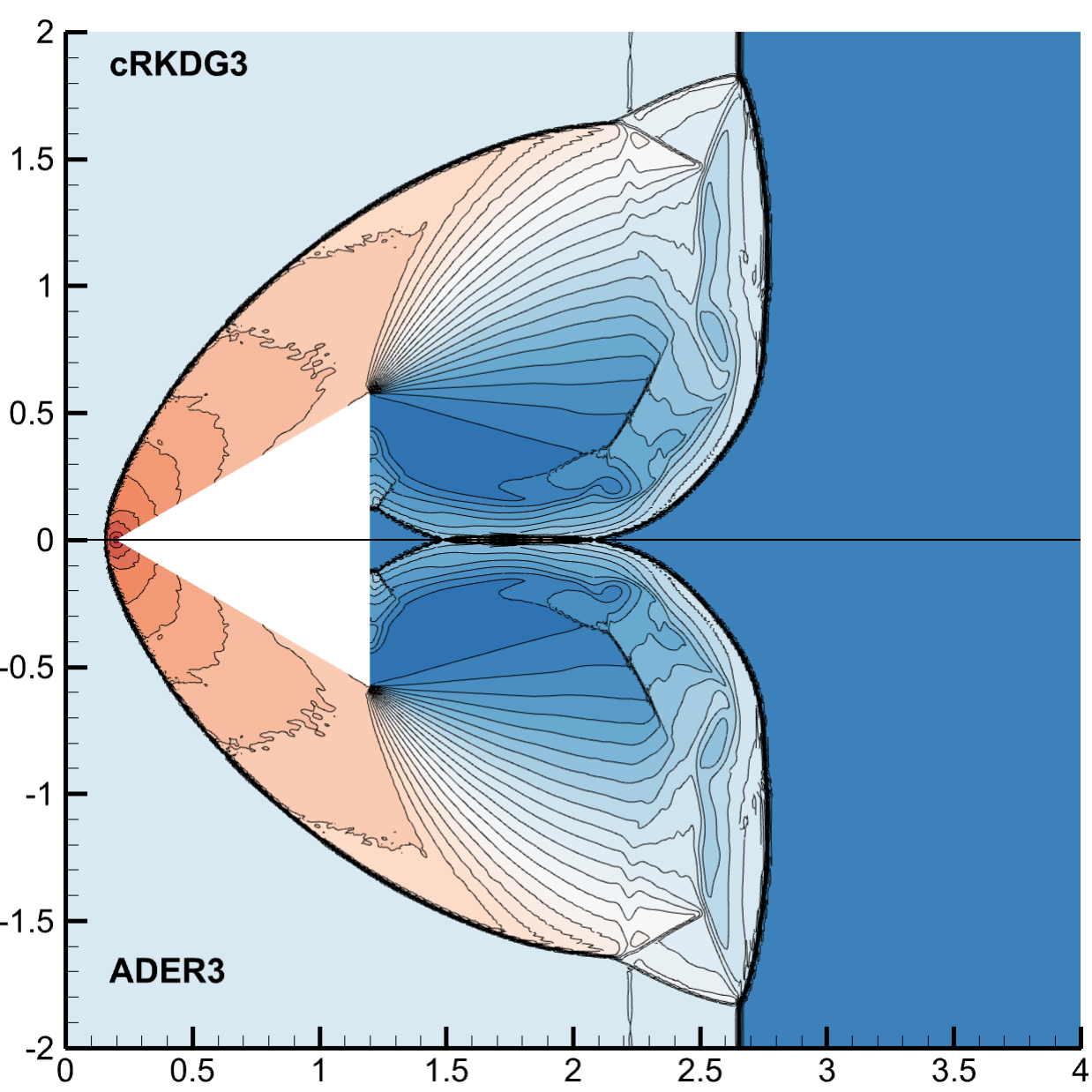}
		\caption{\(T = 0.245\), 30 equally spaced contours from 0.6 to 2.9}
	\end{subfigure}
	\caption{Schardin Mach 10, PACL with RIOE. The upper half of each plot is cRKDG3 and the lower half is ADER3.}
	\label{fig:num-schardin10}
\end{figure}

\begin{example}[Schardin problem at Mach 10]\label{ex:num-schardin10}
	The Schardin problem describes a planar shock passing an equilateral wedge. At Mach~10, we compute the upper half-domain $[0,5]\times[0,2]$ excluding the half-wedge with vertices $(0.2,0),(1.2,1/\sqrt3),(1.2,0)$. The initial shock is at $x=0.2$. The ambient state is $(\rho,v_1,v_2,p)=(1.4,0,0,1)$, and the post-shock state follows the normal-shock Rankine--Hugoniot relations \cite[Example~6.10]{10.1007/s10915-011-9472-8}. The left boundary prescribes post-shock inflow, the right boundary is outflow, and the upper boundary, wedge, and symmetry axis are free-slip reflective. The hybrid coarse mesh (\Cref{fig:mesh-schardin10}) has 702 cells. Three uniform refinements give 44928 cells. \Cref{fig:num-schardin10} shows density. The results show robustness and also ability to capture shock structure. 
\end{example}

\subsection{Accepted and rejected time step attempts}\label{subsec:num-rejections}

\Cref{tab:num-completed} reports accepted and rejected trial time step attempts for discontinuous tests. The counts correspond to the final output time listed for each case. These counts quantify the retry cost of the practical admissibility safeguard; the accuracy tables above compare the two methods at the same spatial resolution.

\begin{table}[h!]
	\centering
	\caption{Accepted steps and rejected trial attempts. ``Off'' denotes sole PACL; ``On'' includes OE/RIOE damping.}
	\label{tab:num-completed}
	\begin{tabular}{llrrrr}
		\toprule
		& & \multicolumn{2}{c}{cRKDG3} & \multicolumn{2}{c}{ADER3}\\
		\cmidrule(lr){3-4}\cmidrule(lr){5-6}
		Case & OE & Accepted & Rejected & Accepted & Rejected\\
		\midrule
		Burgers 2D, $T=0.600$ & Off & 8892 & 0 & 9822 & 0\\
		Buckley--Leverett & On & 1324 & 0 & 1987 & 0\\
		Double rarefaction & On & 1280 & 0 & 1921 & 0\\
		Traffic flow, $T=1.500$ & On & 78236 & 0 & 101468 & 0\\
		Sedov 2D & On & 104171 & 0 & 115096 & 0\\
		$90^\circ$ diffraction & On & 49290 & 0 & 54445 & 0\\
		$120^\circ$ diffraction & On & 87527 & 0 & 96533 & 0\\
		Schardin Mach 1.34 & On & 103460 & 0 & 114227 & 0\\
		Schardin Mach 10 & On & 69866 & 0 & 77062 & 0\\
		\bottomrule
	\end{tabular}
\end{table}

\section{Conclusions}\label{sec:conclusions}

For local polynomial degrees $k_K\geq 1$, we establish admissibility of the weighted patch averages of the completed DG update. Testing by continuous patch hat functions cancels the contribution of conservative interface numerical fluxes, and local volume-flux estimates provide admissibility conditions without first-order bound-preserving building blocks or SSP decompositions. The framework admits exact spatial integration and compatible quadrature and accommodates compact methods of arbitrarily high temporal order. Conservative patch limiting then supplies admissible physical cell averages for the final Zhang--Shu scaling. The limiter also accommodates cells with $k_K=0$ by forming auxiliary vertex moments for the patch correction and retaining the corrected cell average as the new cell value (\Cref{rem:degree-zero-cells}). The procedure is local and conservative in the sense that only the patches directly related to problematic cells are adjusted (\Cref{rem:locality}). Each patch correction preserves the mass assigned to its vertex, and the final cellwise scaling preserves the corrected cell averages (\Cref{rem:local-conservation}). 

The scalar analysis covers cRKDG, converged and fixed-Picard ADER-DG, and LWDG. For Euler, the cRKDG and fixed-Picard ADER results maintain positive predictor density while allowing nonpositive predictor pressure. The correction preserves each conserved quantity over the union of the active patches. Under the local mesh and scaling assumptions of \Cref{thm:pacl-accuracy}, one application of scalar PACL without optional oscillation control preserves approximation order, including at active bounds and flat extrema. The numerical experiments demonstrate the designed convergence rate and bound preservation in the reported near-vacuum and strong-shock tests.

For conservative systems, a further direction is to establish the patch admissibility estimates for predictors formulated in primitive variables. 
Future work will also include positivity preservation for nonconservative hyperbolic systems, including the Godunov--Powell formulation of ideal magnetohydrodynamics, the extension of PACL to nonconforming meshes, and bound-preserving high-order schemes for nonlinear convection--diffusion problems. These directions call for incorporating nonconservative and diffusive contributions into the patch evolution analysis and constructing compatible patch hat functions on nonconforming meshes.


\begin{appendices}

\section{Admissibility enforcement for Lobatto quadrature}\label{app:lobatto-enforcement}

Let $Q_K^{\mathrm L}$ and $Q_K^{\mathrm G}$ denote the tensor-product $(k+1)$-point Lobatto and Gauss volume rules. For affine tensor-product cells with $k\geq2$ and bilinear quadrilaterals with $k\geq3$, both rules integrate $p$ and $\psi_ap$ exactly. These restrictions follow from $2k-1\geq k+1$ and $2k-1\geq k+2$, respectively, with the latter accounting for the physical Jacobian.

If $\mathcal E_K$ contains both node sets, every nonzero local polynomial $p$ nonnegative on $\mathcal E_K$ satisfies
\[
	\frac{Q_K^{\mathrm L}(\psi_ap)}{Q_K^{\mathrm L}(|p|)}
	= \frac{Q_K^{\mathrm G}(\psi_ap)}{Q_K^{\mathrm G}(p)}
	\geq \min_{\bx\in\mathrm{GL}_{k+1}^{\otimes d}} \psi_a(\bx)>0.
\]
Here the Gauss nodes are understood as their physical images. Thus $\eta_{aK}^Q>0$.

The Lobatto nodes alone need not give positivity. On $K=[-1,1]$ with $k=2$, take
\[
	\mathcal E_K=\{-1,0,1\},
	\quad
	\psi_a(x)=\frac{1-x}{2},
	\quad
	p(x)=\frac{x(x+1)}2.
\]
Then $p\geq0$ on $\mathcal E_K$, whereas $Q_K^{\mathrm L}(\psi_ap)=0$ and $Q_K^{\mathrm L}(|p|)>0$, giving $\eta_{aK}^Q=0$.

For $k=1$, the polynomial $p(x)=(1+x)/2$ is nonnegative throughout $K$ and has zero weighted moment under the two-point Lobatto rule for the same hat. Additional enforcement points therefore cannot restore positivity in this case. The two-point rule also fails the required state--hat exactness.

\section{Physical boundary conditions}\label{app:physical-boundaries}

A physical boundary enters the patch update through the numerical flux in the corrector. The cell-local predictor uses the interior polynomial as in \Cref{sec:preliminaries}. Let $\mathcal B_K$ be the set of unpaired physical faces of $K$, with outward unit normal $\bn_e$ on each face $e$. Testing \eqref{eq:corrector} with $\psi_a$ and summing over the physical cells of the patch gives
\begin{equation}\label{eq:physical-boundary-patch}
	\calM_a(u_h^\star)-\calM_a(u_h^n)
	=\tau\sum_{K\in\mathcal T_a}\left[
		\int_K\FF_{K,\tau}(u_K^n)\cdot\nabla\psi_a\dd\bx
		-\sum_{e\in\mathcal B_K}\int_e
		\psi_a\widehat{\mathscr F}^{\,b}_{K,e}(u_K^n)\dd s
	\right],
\end{equation}
where $\widehat{\mathscr F}^{b}_{K,e}$ is the outward, time-averaged boundary flux. Contributions on shared and periodically paired faces cancel as in \Cref{lem:cancellation}. If $\psi_a$ vanishes on the physical boundary, \eqref{eq:weighted-identity} follows immediately. At a boundary vertex, each incident boundary face contributes its own flux, including at corners where two boundary conditions meet.

\subsection{Scalar equations: copied and fixed exterior traces.}

Write $z^-$ for the interior predictor trace at a flux evaluation point. For extrapolation at an outflow boundary, set $z^+=z^-$. The same formula applies to a scalar symmetry condition implemented by even reflection. On a Dirichlet face, set $z^+=g_e$, where $g_e\in[m,M]$ is prescribed and fixed during the step. In both cases the flux uses the same two arguments $\widehat f(z^-,z^+;\bn_e)$.

Write $D_{aK}(v)$ for the expression in brackets in \eqref{eq:physical-boundary-patch} evaluated with old cell polynomial $v$. The scalar comparison proof requires the local estimate
\begin{equation}\label{eq:scalar-physical-boundary-cfl}
	|D_{aK}(v)-D_{aK}(w)|
	\leq L^b_{aK}(\tau)\norm{v-w}_{1,K},
	\quad \tau L^b_{aK}(\tau)\leq\eta_{aK}
\end{equation}
for every active pair and all comparison inputs of \Cref{ass:weighted-lipschitz}. With the other hypotheses of \Cref{thm:main}, increasing one input by $p\ge0$ increases the updated patch mass by at least $(\eta_{aK}-\tau L^b_{aK}(\tau))\norm{p}_{1,K}$. Thus the updated patch average $H_a$ is monotone in each input.

It remains to evaluate $H_a(c)$ for constant interior data. For cRKDG and ADER, all predictor traces then equal $c$ and the corrector weights sum to one. By consistency, the copied face flux is $\widehat f(c,c;\bn_e)=\bfv(c)\cdot\bn_e$. Integration by parts in \eqref{eq:physical-boundary-patch} therefore leaves only the Dirichlet faces $\mathcal B_K^D$:
\begin{equation}\label{eq:physical-boundary-endpoints}
	H_a(c)-c
	=\frac{\tau}{\mathfrak m_a}
		\sum_{K\in\mathcal T_a}\sum_{e\in\mathcal B_K^D}
		\int_e\psi_a
		\left[\widehat f(c,c;\bn_e)-\widehat f(c,g_e;\bn_e)\right]\dd s .
\end{equation}
In particular, copied traces preserve constants. On a Dirichlet face, suppose the outward flux is \emph{nonincreasing} in its second argument on $[m,M]$. Since $m\leq g_e\leq M$, the integrand is nonnegative at $c=m$ and nonpositive at $c=M$. Comparing the old data with these two constants gives
\[
	m\leq H_a(m)\leq\avg{u_h^\star}_a\leq H_a(M)\leq M .
\]

For stage-based fluxes, \eqref{eq:scalar-physical-boundary-cfl} can be checked with the predictor estimates of \Cref{subsec:flux-estimates}. Let $P_{\nu K}$ bound the $L^1(K)$ Lipschitz constant of predictor stage $\nu$, with corrector weight $w_\nu$. Let $\alpha_b$ bound the Lipschitz constant of the composite boundary flux in $z^-$ for the chosen exterior trace on the predictor range. On a boundary cell the inverse trace estimate \eqref{eq:local-inverse-trace} gives the choice
\[
	L^b_{aK}(\tau)=L_{aK}(\tau)
		+\frac{C_{\rm tr}\alpha_b}{h_K}
		\max_{e\in\mathcal B_K}\norm{\psi_a}_{\infty,e}
		\sum_\nu |w_\nu|P_{\nu K}.
\]
For scalar LWDG the comparison uses the actual time-averaged boundary flux. It applies whenever that flux satisfies \eqref{eq:scalar-physical-boundary-cfl} and the endpoint inequalities $H_a(m)\ge m$ and $H_a(M)\leq M$.

\subsection{Euler equations: slip walls, outflow, and prescribed exterior states}

Use cRKDG or fixed-Picard ADER as in \Cref{eul:thm:exact}. Let $\vect U=(\rho,\vect m^\top,E)^\top$ denote the interior predictor trace on a physical face $e$. At each boundary flux evaluation, use $\widehat{\vect f}(\vect U,\vect U^+;\bn_e)$ with
\[
	\vect U^+=
	\begin{cases}
		(\rho,[\vect m-2(\vect m\cdot\bn_e)\bn_e]^\top,E)^\top,
		& \text{on a slip wall},\\
		\vect U,
		& \text{on an extrapolation outflow face},\\
		\vect U_{D,e}(\bx,t),
		& \text{on a face with prescribed exterior data}.
	\end{cases}
\]
On a prescribed-data face, let $\mathcal S_{D,e}\subset\GE$ be a compact set containing all exterior states used on $e$ during the current step. For fixed data this set is a singleton. For a prescribed moving shock it consists of the pre-shock and post-shock states.

For each interior state $\vect U$, let $\mathcal S_e(\vect U)$ contain the exterior trace specified above on a slip wall or an outflow face, and set $\mathcal S_e(\vect U)=\mathcal S_{D,e}$ on a prescribed-data face. On a slip wall or an outflow face, $\mathcal S_e(\vect U)$ is a singleton. Define
\[
	B_{K,e}
	=\sup_{\substack{\vect Z\in\mathcal D_K\\
		\vect U^+\in\mathcal S_e(\stateT_K^{-1}\vect Z)}}
	\left|\stateT_K\widehat{\vect f}
		(\stateT_K^{-1}\vect Z,\vect U^+;\bn_e)\right|_\infty,
	\quad
	\beta_{aK,e}
	=\frac{h_K}{\abs{K}}\int_e\psi_a\dd s .
\]
Assume $B_{K,e}$ is finite on each contributing face. Under the old-state assumptions and predictor condition \eqref{eul:eq:density-CFL} of \Cref{eul:thm:exact}, a sufficient condition for the completed update is
\begin{equation}\label{eq:euler-physical-boundary-cfl}
	C_d\,\omega_{\mathrm{abs}}\,\frac{\tau}{h_K}
	\left(
		C_{\gamma,d}\iota_K\alpha_K^n\beta_{aK}
		+\sum_{e\in\mathcal B_K}\beta_{aK,e}B_{K,e}
	\right)
	\leq \eta_{aK}/4
\end{equation}
for every active pair, with the sum taken over faces on which $\psi_a$ has nonzero trace.

Indeed, \eqref{eul:eq:Phi-bound} bounds the transformed volume flux by $C_{\gamma,d}\iota_K\alpha_K^n$. The predictor condition places each interior trace in $\mathcal D_K$. Applying $\stateT_K/\abs{K}$ to the face terms in \eqref{eq:physical-boundary-patch} bounds their total contribution to the maximum norm of the increment by
\[
	\omega_{\mathrm{abs}}\,\frac{\tau}{h_K}
	\sum_{e\in\mathcal B_K}\beta_{aK,e}B_{K,e}.
\]
Together with \eqref{eul:eq:simple-correction-bound}, \eqref{eq:euler-physical-boundary-cfl} gives $C_d\abs{\vect r_{aK}}_\infty \leq \eta_{aK}/4$ for the increment including the boundary faces. The conclusion follows from \Cref{eul:lem:budget} and the physical patch identity.

For an HLL flux, bounded signal speeds on the indicated trace pairs suffice for $B_{K,e}<\infty$, provided the flux also has a bounded value at coincident zero speeds. Indeed, in the branch $s_L<0<s_R$, the coefficients of the two physical fluxes are nonnegative and sum to one, while the jump coefficient has magnitude at most $\min(-s_L, s_R)$. The one-sided branches use a single physical flux. The trace pairs lie in a compact set with density bounded away from zero, so their physical fluxes and clipped acoustic speeds are bounded. The resulting boundary constant may depend on the exterior data and the old cell mean.

\subsection{Compatible quadrature}

Use $Q_K$ and $Q_e^f$ in \eqref{eq:physical-boundary-patch} under \Cref{ass:quadrature}. Exact state--hat pairings identify its left-hand side with the physical patch mass increment. For the scalar proof, replace $\norm{v-w}_{1,K}$ by $Q_K(|v-w|)$ and $\eta_{aK}$ by $\eta_{aK}^Q$. The discrete Green identity \eqref{eq:discrete-green} gives \eqref{eq:physical-boundary-endpoints} with $Q_e^f$ in place of the face integrals. Nonnegative Dirichlet face weights preserve the two endpoint signs. For other rules one can check the two discrete endpoint inequalities directly.

For Euler, impose the predictor condition in \eqref{eq:euler-Q-CFL}. 
For each contributing boundary face, define $B_{K,e}^Q$ as in the exact-integral case with $\mathcal D_K^Q$ in place of $\mathcal D_K$. 
If the face rule has nodes $\bx_{e,q}$ and physical weights $\omega_{e,q}^f$, set
\[
	\beta_{aK,e}^Q
	=\frac{h_K}{\abs{K}}
	\sum_q \abs{\omega_{e,q}^f}\,\psi_a(\bx_{e,q}).
\]
The nodal volume-flux estimate in \Cref{cor:euler-quadrature} and the face-flux bound $B_{K,e}^Q$ give the sufficient condition
\[
	C_d\,\omega_{\mathrm{abs}}\,\frac{\tau}{h_K}
	\left(
		C_{\gamma,d}\iota_K^Q\alpha_{K,Q}^n\beta_{aK}^Q
		+\sum_{e\in\mathcal B_K}\beta_{aK,e}^Q B_{K,e}^Q
	\right)
	\leq \eta_{aK}^Q/4
\]
for every active pair. The discrete patch identity then gives $\avg{\vect U_h^\star}_a\in\GE$ on the physical mesh. 
The patch limiter applies as in \Cref{sec:patch-limiting}.

\end{appendices}

\section*{Acknowledgements}

Part of this research was conducted using computational and visualization resources and services provided by the Center for Computation and Visualization at Brown University. The authors thank Professor Chi-Wang Shu of Brown University for his kind support and helpful comments. Both authors thank OpenAI for providing access to ChatGPT and Codex through the ChatGPT for Academic Researchers program. Zheng Sun also thanks Anthropic for providing access to Claude through the Claude Team plan for research labs.

\section*{Statement of Usage of Artificial Intelligence}

The authors used ChatGPT and Codex (OpenAI) and Claude (Anthropic) for mathematical discussions, manuscript drafting and revision, and \LaTeX{} preparation. These tools also assisted with the development, execution, and verification of computer code used for the numerical experiments and to produce the associated figures and tables. The authors reviewed and verified all AI-assisted content and computational results, made revisions as needed, and take full responsibility for the content of the manuscript.


\printbibliography[heading=bibintoc] 



\end{document}